%% file: mr.tex
\documentclass[article]{siamart}
\usepackage{amsmath,amssymb}
\usepackage{caption}
\usepackage{hyperref}
\usepackage[utf8]{inputenc}
\usepackage{mathtools}
\usepackage{color,graphicx}
\usepackage{thmtools, thm-restate}
\usepackage{enumitem}
\usepackage{algorithmicx}
\usepackage{algpseudocode}
\usepackage{tikz}
\usetikzlibrary{patterns}
\usetikzlibrary{arrows}
\usetikzlibrary{3d,calc}

\input{macros.tex}

\begin{document}

\ifpdf
\DeclareGraphicsExtensions{.pdf, .jpg, .tif}
\else
\DeclareGraphicsExtensions{.eps, .jpg}
\fi

\title{\TheTitle}

\author{Donsub Rim%
  \thanks{Department of Applied Physics and Applied Mathematics,  %
  Columbia University, New York, %
  NY 10027 (\email{dr2965@columbia.edu}, \email{kyle.mandli@columbia.edu}).}%
  \and %
  Kyle T. Mandli\footnotemark[1]  %
}
\maketitle

\begin{abstract} 
When approximating a function that depends on a parameter, one encounters many
practical examples where linear interpolation or linear approximation with
respect to the parameters prove ineffective.  This is particularly true for
responses from hyperbolic partial differential equations (PDEs) where linear,
low-dimensional bases are difficult to construct. We propose the use of
\emph{displacement interpolation} where the interpolation is done on the
optimal transport map between the functions at nearby parameters, to achieve an
effective dimensionality reduction of hyperbolic phenomena. We further propose
a multi-dimensional extension by using the intertwining property of the Radon
transform.  This extension is a generalization of the added{classical
translational representation of Lax-Philips [Lax and Philips, \emph{Bull.
Amer. Math. Soc.} 70 (1964), pp.130--142].}
\end{abstract}

\section{Introduction}

Linear interpolation or linear approximation is a concept that is ubiquitous in
computational mathematics.  Perhaps the most familiar setting is where one
approximates an arbitrary function by a linear sum of functions from a
carefully chosen basis $\cV$.  For example, the basis can be chosen as the
Fourier basis, Wavelet basis, or Chebyshev polynomial basis, depending on the
application at hand \cite{szego,trefethen,zygmund_fefferman,daubechies}.  The
approximation with respect to these well-studied bases are guaranteed to be
optimal in some sense, often meaning that the number of terms in the linear sum
that achieves a desired accuracy is small.  Put in other terms, an accurate
approximation is presumed to belong to some low-dimensional subspace of the
chosen basis $\cV$. 

A more challenging problem arises when the function to be approximated depends
on a parameter.  For one fixed parameter value a linear approximation may
indeed be optimal, accurately representing the approximand with a few members
of $\cV$, but a different set of members may be needed for a good approximation
at a different parameter value. In the worst case, one would need a large
subset of $\cV$ to accurately approximate the function over all parameter
values, and a low-dimensional representation would be possible only locally.
Mathematically speaking, we are describing the situation where the Kolmogorov
$N$-width decays slowly (e.g. linearly) with respect to the dimension $N$ of
the basis $\cV$ \cite{crb-book}.  Unfortunately, this is not a rare worst-case
scenario, but is the typical behavior for functions of interest in modeling
wave propagation: often the energy of the wave is concentrated at different
spatial locations for different parameter values, so that finding a global
representation of the wave profile using a fixed low-dimensional subspace is
not possible \added{\cite{marsden1,amsallem,Schulze15,Welper17p}. This
high-dimensional nature of wave phenomena can also be characterized terms of
the separability of the Green's function for the Helmholtz equation
\cite{engquist18}, however this difficulty arises even for the simplest of
examples \cite{rim18reversal}.}

\added{Naturally, this also has important implications for the field of
uncertainty quantification (UQ). Many numerical numerical methods devised for
propagating uncertainties rely on exploiting some low-dimensional linear
subspace for the random parameters. As a prime example, generalized polynomial
chaos (gPC) expansions use polynomial bases to efficiently compute the
statistics of some quantity of interest (QoI) \cite{xiu10}. However, the slow
decay of the Kolmogorov $N$-width with respect to the uncertain parameters also
implies slow convergence of the expansion \cite{welper17}. A closely related
subject is that of reduced order models (ROMs) \cite{crb-book,siamrev-survey}.
Once a ROM of a parametrized system is constructed, even na\"ive Monte Carlo
sampling of the full solution becomes feasible. This makes ROMs valuable in
various UQ problems, whenever they could be constructed. Nonetheless,
projection-based model reduction methods suffer from precisely the same issue
pointed out above, due to the lack of a low-dimensional linear basis.}

In this paper, we adopt an alternative approach in which we compute the
transport map between functions in order to find a low-dimensional structure in
the transport maps themselves.  Given two wave profiles, we will find the map
that will transport each unit mass in the first profile in the most optimal
manner to the second profile.  This is precisely an optimal transport problem,
and the interpolation we are proposing is referred to as \emph{displacement
interpolation} in the optimal transport literature \cite{mccann97,villani03,
villani2008}.

For our purposes, this transport map is found by using the simplest solution to
the Monge-Kantorovich \added{optimal transport} problem
\cite{Kantorovich,Benamou}  in a single spatial dimension (1D), called
\emph{monotone rearrangement}. \added{The problem minimizes the so-called
Wasserstein metric between two probability distributions \cite{dobrushin70},}
and we will show that a very simple computation will produce \added{the
minimizing} transport map.  Moreover, although the optimal transport problem
poses that the two profiles have to be non-negative, this restriction is easily
removed by implementing an integral formulation \eqref{eq:mr1d} rather than a
differential formulation \eqref{eq:ode} of the problem. \added{We also extend
the computation of transport maps to apply to functions of arbitrary sign}. 
While this may be satisfactory in 1D, wave propagation takes on a much more
complicated form in multiple spatial dimensions. We extend our displacement
interpolation procedure in 1D to multiple dimensions via a dimensional
splitting that exploits the intertwining property of the Radon transform
\cite{laxphilips,Helgason2011,radonsplit}.  That is, instead of dealing with
the multi-dimensional problem directly, we solve a collection of 1D
Monge-Kantorovich problems.  \added{This a natural extension of the
translational representation of the solution to the wave equation by Lax and
Philips \cite{laxphilips}.} As a result of this interpolation procedure, one
discovers a smooth map that interpolates between two wave profiles, even though
the profiles themselves can be nearly orthogonal to each other. 

\added{Individual components that form the interpolation techniques proposed in
this paper are closely related to recent works that have arisen in various
contexts. The subject matter spans different fields and its literature is
growing rapidly, so we will only provide a limited overview here.}  

\added{\emph{Model reduction}. The idea of finding a low-dimensional structure
by applying dimensionality reduction in transport maps was also proposed in
\cite{iollo14}, and these low-rank transport modes were called \emph{advection
modes}.  Nonetheless, the work suggested using linear programming to compute
the optimal transport map in place of the explicit monotone rearrangement, and
the multi-dimensional extension was also done by solving an equivalent
multi-dimensional optimal transport problem. This incurs a high computational
cost, and prohibits its use in practical PDE applications. We note that the
authors adapted the interpolation techniques here to build a ROM for
paramterized scalar conservation laws \cite{rim18}, motivated by an elegant
relationship between optimal transport and scalar conservation laws observed in
\cite{bolley05}.}

\added{\emph{Image Processing.} Displacement interpolants, or more generally
Wasserstein barycenters \cite{agueh11}, are useful in image processing for
tasks such as shape interpolation, warping, color transfer or texture mixing.
In these applications, the Wasserstein metric is commonly used along with an
optimization procedure but again, the computational cost was too high
\cite{bonneel11}.  In an effort to find cheaper approximations, both the
monotone rearrangement \cite{delon04}, as well as an approximation to the
Wasserstein metric called the \emph{sliced Wasserstein distance} were
considered in \cite{rabin12}, and the combination of using monotone
rearrangement along with the Radon transform was also proposed as an
approximation to the sliced Wasserstein barycenters themselves
\cite{pitie07,bonneel15}. Remarkably, the transform was proposed chiefly as a
way to reduce computational cost and the articles report appearance of
artifacts when Radon transforms were used to interpolate between images, which
is not unexpected since the most natural interpolants between common images
tend to be closer to the motion of a rigid body rather than that of waves. 
The link between the Radon transforms and hyperbolic PDEs through the
intertwining property \eqref{eq:intertwine} was not made, and the interpolant
with artifacts may in fact be the \emph{correct} interpolant we seek for our
purposes. On the other hand, breakthroughs in reducing computational costs of
the optimal transport problem via entropic-regularization
\cite{cuturi13,benamou15} has significantly alleviated the computational
burden, so direct computations on the multi-dimensional optimal transport
problem appears to be the more desired approach for this class of applications
\cite{solomon15}.}

\added{\emph{Wasserstein metric as a misfit function.} The Wasserstein metric
itself is of great interest apart from our application for interpolation or
dimensionality reduction. Also called the \emph{Mallows distance}
\cite{mallows72}, it has long been studied in the statistics literature in
relation to goodness-of-fit tests for non-normality or as measures of
similarity between histograms \cite{barrio99,dewet02,munk02}, due to its
robustness with respect to location-scale families of probability
distributions. More recently, the metric has been employed as a misfit function
in various applications. For example, it was used in seismic inverse problems
\cite{yang18,metivier16} to overcome common problems such as cycle-skipping
that entails the use of the usual $\ell^2$ misfit function. Its application is
also being explored in machine learning, for instance as the loss function in
supervised learning \cite{frogner15}.}

\added{Our main contribution is in the generalization of the classical
translational representation of Lax-Philips by combining monotone rearrangement
with the Radon transform, further extending the approach in conjunction with
more general transforms, and applying them for dimensionality reduction of
problems involving wave propagation exhibiting high-dimensional behavior.} The
paper is organized as follows. The displacement interpolant is derived in
\cref{sec:dinterp} through the integral formulation \eqref{eq:mr1d} of the
monotone rearrangement problem.  The simplest case of interpolating between two
strictly positive 1D profiles is considered first, then restrictions on the
signs are gradually removed.  Then, a multi-dimensional extension is defined
through the Radon transform.  Details on numerical implementation and numerical
examples illustrating various aspects of the interpolation scheme are presented
in \cref{sec:numerics}.

\section{Displacement interpolation}\label{sec:dinterp}

The most common approach to approximating a paramterized function is to apply
the method of separation of variables and then to perform linear approximation.
Once the parameter function is separated, one chooses a suitable linear basis,
then proceeds to find a linear combination of its members that approximates
this function. 

To make this more explicit, suppose that $u$ is a function that depends on a a
vector of parameters $\balpha$, and that we have evaluated the function at
various parameter values $\balpha_1, ..., \balpha_N$.  That is, we have
obtained the set of functions
\[
\cU = \{u_n = u(\balpha_n): n = 1, 2, ... N\}.
\]
The goal is to interpolate these function with respect to the parameters.

Throughout this paper, we will refer to each member in $\cU$ itself as a
function, and when referring to multiple members, we will refer to them in
plural as \emph{functions}, although they are merely the function $u$ evaluated
as given parameter values $\{\balpha_n\}$.  We will also suppress the
dependence of $u$ on its domain, for simplicity of notation. 

If there is a suitable basis $\cV$ of $\cU$ with small dimension $M$,
\[
    \cV = \{v_1, v_2, ... , v_M \},
\]
then a natural choice for approximating the solution $u$ at a new parameter
value $\balpha$, is to construct $\tilde{u}$ that is a linear combination of
the basis functions,
\begin{equation}
    \tilde{u}(\balpha) = \sum_{j=1}^M c_j(\balpha) v_j.
    \label{eq:linterp}
\end{equation}
There is a wealth of choices for the coefficients $c_j(\balpha)$, but they are
typically chosen to minimize the error in some sense, for example $\sum_n
\lVert \tilde{u}(\balpha_n) - u_n \rVert_2^2$.  When the coefficients are
chosen in a way that this error vanishes (that is, $\tilde{u}(\balpha_n) =
u_n$) we call the approximation $\tilde{u}$ an interpolant.

\begin{figure}
\centering
\begin{tabular}{cc}
    \includegraphics[width=0.35\textwidth]{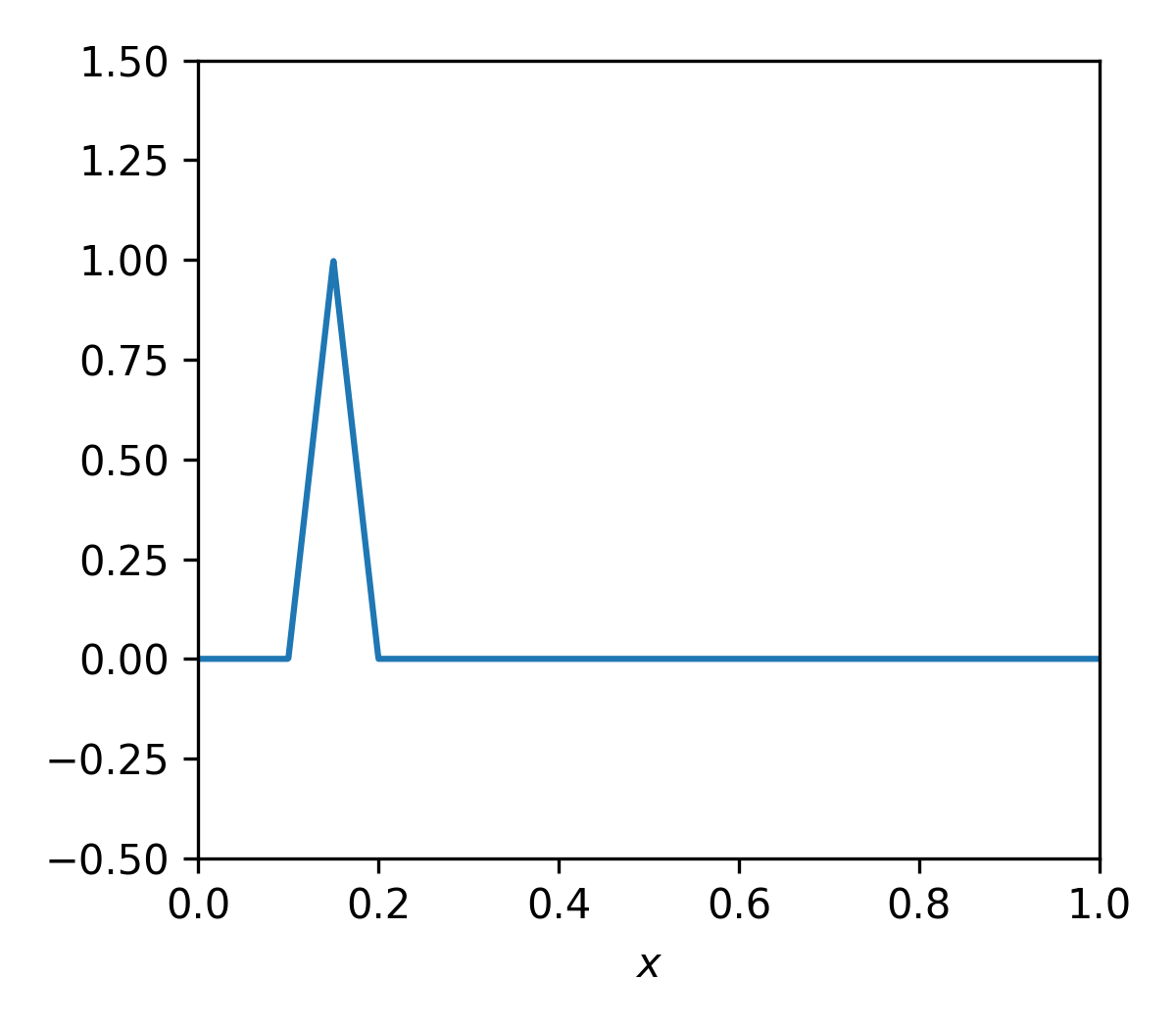}
    &
    \includegraphics[width=0.35\textwidth]{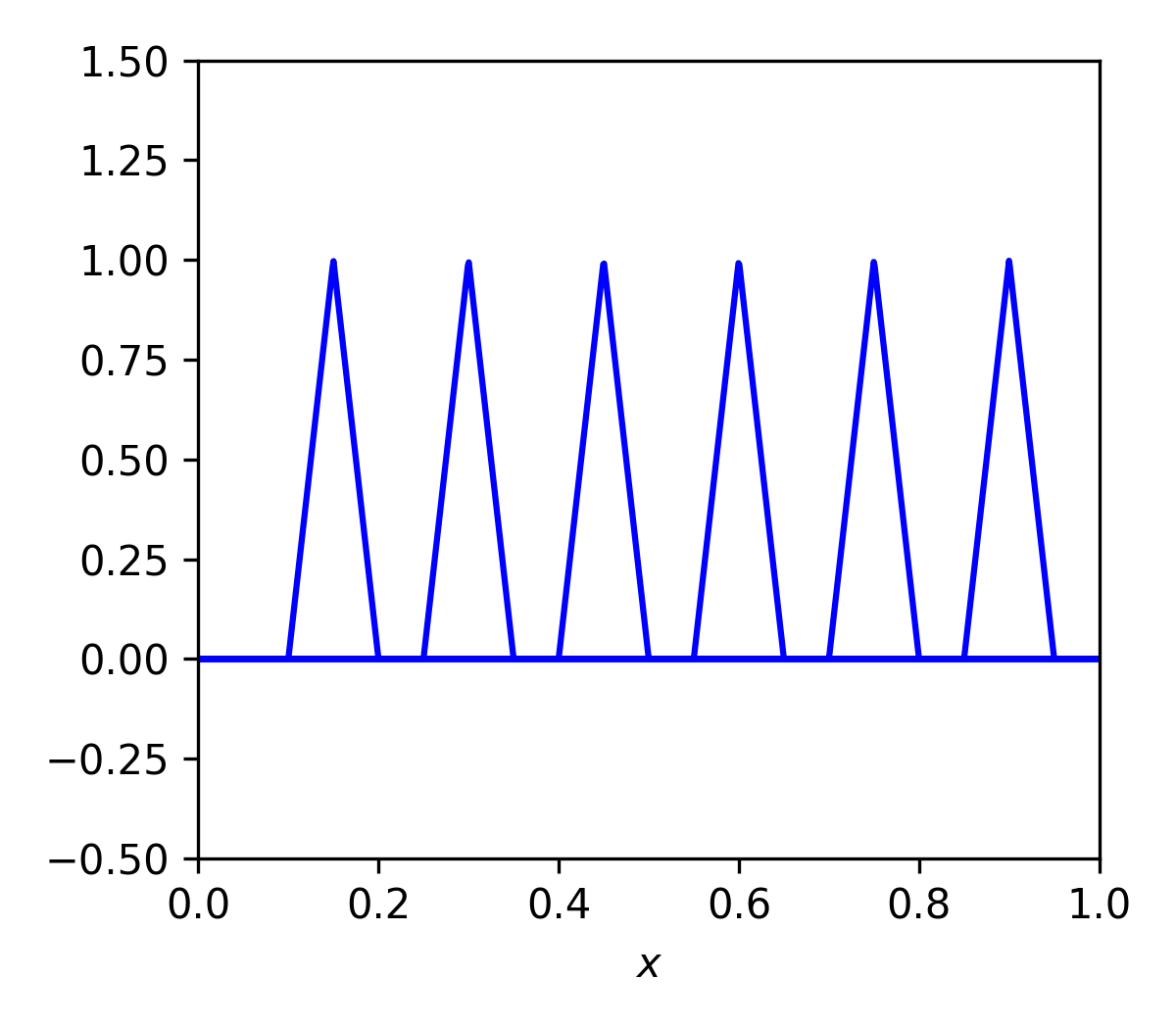}\\
\end{tabular}
\caption{The initial condition in \eqref{eq:transport1d}, the hat function
$\phi(x-3w)$ with $w =0.05$ (left). The solution of the advection equation
\eqref{eq:transport1d} at different times $t_n = 3(n-1)w$ overlayed on top of
each other, for $n=1, 2, ..., 6$ (right).}
\label{fig:hats}
\end{figure}

Although linear approximation (or interpolation) is a very powerful tool, there
are examples where linear approximation with respect to the parameters fail to
be optimal, due to the fact that there exists no low-dimensional linear
subspace $\cV$, even though $\cU$ may be low-dimensional in some other sense.
An example of particular interest is when $u(\balpha)$ is the solution to a
hyperbolic PDE.  To demonstrate this difficulty, suppose the parameter is the
time variable, $\balpha = t$, and $u(t)$ is the solution to a 1-dimensional
transport equation
\begin{equation}
    \left\{
    \begin{aligned}
        u_t + u_x &= 0 \quad \text{ in } (0,1),\\
           u(x,0) &= \phi(x - 3w), \\
           u(0,t) &= u(1,t),
    \end{aligned}
    \right.
    \label{eq:transport1d}
\end{equation}
where $w = 0.05$, and $\phi(x)$ is a hat function centered at $0$ of width
$2w$, 
\begin{equation}
    \phi(x;w) = 
    \begin{cases}
        \frac{1}{w} (x + w) &  \text{ if } -w \le x \le 0,\\
      - \frac{1}{w} (x - w) & \text{ if } 0 \le x \le w,\\
        0 & \text{ otherwise.}
    \end{cases}
    \label{eq:hatf}
\end{equation}
Suppose that we are given $u(t)$ at various times $t_1 < t_2 < \cdots < t_M$ as
shown in \cref{fig:hats}.  These functions $u(t_n)$ are all orthogonal to each
other since $\supp u(t_n) \, \cap \, \supp u(t_m) = \emptyset$ for $n \ne m$,
hence no low-dimensional basis $\cV$ can be found for $\cU =\{u(t_n)\}$.
Nonetheless, the functions are merely translates of each other.

Naturally, our goal is to construct an interpolant $\tilde{u}(t)$ that
satisfies $\tilde{u}(t_j) = u(t_j)$ that exploits this translation symmetry
without having to compute additional solutions (that is, adding to $\cU$) or
adaptively refining the linear basis (adding to $\cV$).  We propose to
construct an interpolant using \emph{displacement interpolation}
\added{\cite{mccann97}}, by solving a simple optimal transport problem between
the functions $\{u(t_n)\}$ then find the low-rank structure in the computed
transport maps. To keep the operations as simple as possible, we make use of a
simple optimal transport solution, called \emph{monotone rearrangement}
\added{\cite{mallows72,villani03}}. The resulting solution is straightforward
to compute, making it a potentially useful tool in a wide variety of contexts.

\subsection{Monotone rearrangement} \label{ssec:mr}

Suppose we wish to interpolate between two functions $u_1, u_2 > 0$ in
$\added{\cL:=}L^2(\Omega) \cap L^\infty(\Omega)$ that satisfies
\begin{equation}
        \int_\Omega u_1 \, dx = \int_\Omega u_2 \, dx = 1.
        \label{eq:massbalance}
\end{equation}
In our setting, the two functions can be considered to be the parametrized
function $u(\alpha)$ evaluated at two parameter values $\alpha_1$ and
$\alpha_2$, that is, $u_1 = u(\alpha_1)$ and $u_2 = u(\alpha_2)$.  In the
optimal transport literature these two densities are referred to as probability
distribution functions (PDFs), and we will follow this convention as long as
the functions we are referring to are non-negative.

We define the optimal transport problem as:
\beq
    \begin{aligned}
    \text{find } M \text{ that minimizes } 
    \int_\Omega \lvert x - M(x) \rvert u_1(x) \, dx, \\
    \quad
    \text{ subject to } u_1(M(x)) = u_2(x).
    \end{aligned}
    \label{eq:ot_prob}
\eeq
This is a version of the mass transfer problem posed by Monge
\cite{Kantorovich,villani03}\added{, and the minimizing functional is called
the Wasserstein metric \cite{dobrushin70}.} \remove{} The problem is well-known
to be degenerate even in the simple 1-dimensional setting, as illustrated by
the book-shifting example \cite{gangbo96}.  However, there is a particularly
simple explicit formulation of the solution when $\Omega = \RR$: we seek a
non-decreasing function $U(x)$ such that
\begin{equation}
    \int_{-\infty}^x u_1 \, dx = \int_{-\infty}^{U(x)} u_2 \, dx.
    \label{eq:mr1d}
\end{equation}
This solution is called monotone rearrangement, as it rearranges the mass in a
monotone manner. Physically, this implies that particle trajectories do not
cross, when they are transported across the domain by the transport map.
Without loss of generality, we will first consider the case of approximating
between two such 1-dimensional functions $u_1$ and $u_2$ that are positive in
the domain $\RR$. The restrictions on the dimension, the number of functions,
and the strict positivity, will all be removed as we progress.

\remove{}
\added{Upon differentiating \cref{eq:mr1d},} we obtain the ordinary
differential equation
\begin{equation}
    U'(x) = \left. u_1(x) \right/ u_2(U(x)),
    \label{eq:ode}
\end{equation}
valid when $u_1, u_2 > 0$ with regularity conditions, for example if $u_1,u_2$
are smooth. 

As an alternative to solving the differential equation, we may compare the
cumulative distribution functions (CDFs) of $u_1$ and $u_2$:
\begin{equation}
    U_1(x) := \int_{-\infty}^x u_1(y) \, dy
    \quad \text{ and } \quad
    U_2(x) := \int_{-\infty}^x u_2(y) \, dy.
    \label{eq:CDFs}
\end{equation}
Since $u_1,u_2$ are both positive $U_1,U_2$ are both strictly monotone
increasing. Moreover, the CDFs are continuous hence the inverses of $U_1$ and
$U_2$ exist. The CDF $\widetilde{U}_\lambda$ of the displacement interpolant
(for given weight $\lambda$) is defined by letting its inverse equal the linear
combination of $U_1^{-1}$ and $U_2^{-1} $\cite{barrio99},
\begin{equation}
    \widetilde{U}_\lambda^{-1}(y) 
    := (1-\lambda) U_1^{-1}(y) + \lambda U_2^{-1}(y)
    \quad \text{ where } 0 \le \lambda \le 1.
    \label{eq:interp}
\end{equation}
Then, letting $\lambda(\alpha)$ be the barycentric coordinate\footnote{%
    Here, if $\alpha_1 < \alpha_2$ then %
    $\lambda(\alpha) = (\alpha - \alpha_1)/(\alpha_2 - \alpha_1)$ %
    for $\alpha_1 \le \alpha \le \alpha_2$.}%
of $\alpha$ with respect to the nodes $\{\alpha_1, \alpha_2\}$, the
displacement interpolant is obtained by taking the derivative of its CDF,
\begin{equation}
    \tilde{u}_\alpha := 
    \cI(\alpha;u_1,u_2) := \frac{d}{dx} \widetilde{U}_{\lambda(\alpha)}.
    \label{eq:interp1def}
\end{equation}
It is easy to see that $\tilde{u}_\alpha$ is indeed an interpolant; when
$\alpha = \alpha_1$ ($\alpha = \alpha_2$) one recovers $\tilde{u}_\alpha = u_1$
($\tilde{u}_\alpha = u_2$, respectively). We will also denote this
interpolation by $\cI(\alpha;u_1,u_2)$, as indicated above.

Nonetheless, this interpolant cannot be directly used for the 1D transport
equation example above \eqref{eq:transport1d} because the solution $u$ vanishes
in some portions of the domain. Therefore, we next generalize the definition
for the interpolant $\tilde{u}_\alpha$ \eqref{eq:interp1def}, so that $u_1,u_2$
are allowed to be zero in some parts of the domain.

\subsection{Two functions with non-negative values}

In this section, we remove the requirement that $u_1,u_2$ both have to be
strictly positive, and define the interpolation procedure for two functions
with non-negative values. The intervals where these functions vanish appear as
intervals in which the respective CDFs is constant. Now the CDFs are no longer
invertible, but their graphs certainly can be reflected across the line
$\{(x,y):x = y\}$ and these intervals can be conveniently represented as
Heaviside jump functions.

Suppose $u_1, u_2$ are allowed to vanish, that is, $u_1, u_2 \ge 0$.  This
implies that the CDFs $U_1,U_2$ may have intervals in which they are constant.
We will restrict our discussion to functions $u_1,u_2$ for which there are only
finitely many such intervals.  Now, let us define a set of values for which
$U_n(x)$ is constant,
\begin{equation}
    \cY_n := \{y: U_n^{-1}(y) \text{ is not a singleton} \},
    \quad \text{ for } n = 1,2,
    \label{eq:cY}
\end{equation}
then enumerate the members of the set $\cY_n$ in increasing order,
\[
    y_1^n < y_2^n < \cdots < y_k^n < \cdots < y_K^n.
\]
We define the end-points of the intervals as follows,
\[
    a_k^n := \min\{ U^{-1}_n(y_k^n) \}
    \quad \text{ and } \quad
    b_k^n := \max\{ U^{-1}_n(y_k^n) \},
\]
and let $\ell^n_k = b^n_k - a^n_k$ denote the length of these intervals.

On the other hand, when $y \notin \cY_n$, the inverse $U^{-1}_n(y)$ is
well-defined. Let us define $W_n$ as the left limit of $U_n^{-1}$ where it
exists,
\[
        W_n(y) := \lim_{z \to y^+} U_n^{-1}(z).
\]
We define $U^\dagger_n$ to be a \emph{pseudo-inverse} of $U_n$ in which the
intervals where $U_n$ is constant are represented using the Heaviside jump
function $H$,
\begin{equation}
    U^\dagger_n (y):= W_n(y) + \sum_{k=1}^K \ell_k^n \, H(y - y_k^n).
    \label{eq:pinv}
\end{equation}
Note that the pseudo-inverse encodes all the necessary information to recover
$U_n$, so we can extend the interpolation previously defined in
\eqref{eq:interp} for nonnegative profiles.  The CDF $\widetilde{U}_\lambda$ of
the displacement interpolant (for given weight $\lambda$) is naturally given in
terms of its pseudo-inverse,
\begin{equation}
\widetilde{U}_\lambda^{\dagger} := (1-\lambda) U_1^{\dagger}(y) 
    + \lambda U_2^{\dagger}(y)
    \quad \text{ where } 0 \le \lambda \le 1.
    \label{eq:interp2}
\end{equation}
Although $\widetilde{U}_\lambda^{\dagger}$ may involve Heaviside functions, one
can uniquely determine $\widetilde{U}_{\lambda}$ from its pseudo-inverse
defined here.  Letting $\lambda(\alpha)$ denote the barycentric coordinate of
$\alpha$ with respect to the nodes $\{\alpha_1, \alpha_2\}$, the interpolant is
given by the derivative of the CDF,
\begin{equation}
    \tilde{u}_\alpha := 
    \cI(\alpha;u_1,u_2) := \frac{d}{dx} \widetilde{U}_{\lambda(\alpha)}.
    \label{eq:interp2def}
\end{equation}
This removes the restriction that the two functions have to be strictly
positive. 

\subsection{Two functions of arbitrary total mass}

Suppose the two functions $u_1, u_2 \ge 0$ do not integrate to one. We can
incorporate this additional generalization into our definition by simply
scaling \eqref{eq:interp2def} by a multiplicative constant. We simply let
\begin{equation}
    M_n := \int_{-\infty}^\infty u_n \,dx,
    \quad 
    \text{ for } n = 1,2,
\end{equation}
to normalize $U_1,U_2$, so that the linear combination
\eqref{eq:interp2} now becomes
\begin{equation}
    \widetilde{U}_\lambda^\dagger(y) := (1-\lambda) U_1^{\dagger}(M_1 y) 
                                          + \lambda U_2^{\dagger}(M_2 y)
    \quad \text{ where } 0 \le \lambda \le 1.
    \label{eq:interp2scaled}
\end{equation}
Letting $\lambda(\alpha)$ denote the barycentric coordinate, we scale the
interpolant linearly in the final step,
\begin{equation}
    \tilde{u}_\alpha := 
    \cI(\alpha;u_1,u_2) := 
                 \left[ (1-\lambda(\alpha))M_1 + \lambda(\alpha) M_2  \right]
                         \frac{d}{dx} \widetilde{U}_{\lambda(\alpha)},
    \label{eq:sinterp}
\end{equation}
This guarantees that although the total integral of $u_1$ and $u_2$ are not the
same, that of $\tilde{u}_\alpha$ will interpolate between the two values
linearly, that is,
\[
    \int_{-\infty}^\infty \tilde{u}_\alpha \, dx = 
    (1- \lambda(\alpha)) \int_{-\infty}^\infty u_1\, dx
      +  \lambda(\alpha) \int_{-\infty}^\infty u_2\, dx.
\]

\begin{figure}
\centering
\begin{tikzpicture}[scale=0.8]
\begin{scope}
    \draw (5.5,1.8) node[anchor=east]{$\cI(\alpha; u_1,u_2)$};
\draw[black] (1.25, 1.25) node {$u_1$};
\draw[white] (0,-.6) -- (0, -.4);
\draw[white] (0,-.5) -- (5.5,-.5);
\draw[white] (5.5,-.6) -- (5.5, -.4);
\draw[black,thick] (0.,.42) -- (4.,.42);
\draw[black,thick] (5.,.42) -- (5.5,.42);
\draw [black,dashed,thick,domain=0:180] 
plot ({0.5 + .5*cos(\x)}, {.42 + 1.5*sin(\x)});
\draw [black,dashed,thick,domain=0:180] 
plot ({1.5 - .5*cos(\x)}, {.42 - 0.8*sin(\x)});
\draw [black,thick,domain=0:180] plot ({4.5 + .5*cos(\x)}, {.42 + .85*sin(\x)});
    \draw[black] (4.5, .8) node {$u_2$};
\end{scope}
\begin{scope}[shift={(7,0)}]
\draw (0.90,1.8) node[anchor=west]{$u_1 = u_1^+ - u_1^-$};
\draw (4.0,1.8) node[anchor=west]{$u_2 = u_2^+$};
\draw[white] (0,-.6) -- (0, -.4);
\draw[white] (0,-.5) -- (5.5,-.5);
\draw[white] (5.5,-.6) -- (5.5, -.4);
\draw[black,thick] (0.,.42) -- (4.,.42);
\draw[black,thick] (5.,.42) -- (5.5,.42);
\draw [black,dashed,thick,domain=0:180,pattern=north west lines] 
plot ({0.5 + .5*cos(\x)}, {.42 + 1.5*sin(\x)});
\draw [black,dashed,thick,domain=0:180,pattern=north east lines] 
plot ({1.5 - .5*cos(\x)}, {.42 - 0.8*sin(\x)});
\draw[black] (0.5, 1.25) node[fill=white,inner sep=1pt] {$u_1^+$};
\draw[black] (1.5, -.008) node[fill=white,inner sep=-0.3pt] {$-u_1^-$};
\draw [black,thick,domain=0:180] plot ({4.5 + .5*cos(\x)}, {.42 + .85*sin(\x)});
    \draw[black] (4.5, .8) node {$u_2^+$};
\end{scope}
\begin{scope}[shift={(0,-3)}]
\draw (5.5,1.8) node[anchor=east]{$\cI(\alpha; u_1^+,u_2^+)$};
\draw[white] (0,-.6) -- (0, -.4);
\draw[white] (0,-.5) -- (5.5,-.5);
\draw[white] (5.5,-.6) -- (5.5, -.4);
\draw[black,thick] (0.,.42) -- (4.,.42);
\draw[black,thick] (5.,.42) -- (5.5,.42);
\draw [black,dashed,thick,domain=0:180,pattern=north west lines] 
plot ({0.5 + .5*cos(\x)}, {.42 + 1.5*sin(\x)});
\draw[black] (0.5, 1.25) node[fill=white,inner sep=1pt] {$u_1^+$};
\draw [black,thick,domain=0:180] plot ({4.5 + .5*cos(\x)}, {.42 + .85*sin(\x)});
    \draw[black] (4.5, .8) node {$u_2^+$};
\end{scope}
\begin{scope}[shift={(7,-3)}]
\draw (5.5,1.8) node[anchor=east]{$-\cI(\alpha; u_1^-,u_2^+)$};
\draw[white] (0,-.6) -- (0, -.4);
\draw[white] (0,-.5) -- (5.5,-.5);
\draw[white] (5.5,-.6) -- (5.5, -.4);
\draw[black,thick] (0.,.42) -- (4.,.42);
\draw[black,thick] (5.,.42) -- (5.5,.42);
\draw [black,dashed,thick,domain=0:180,pattern=north east lines] 
plot ({1.5 - .5*cos(\x)}, {.42 - 0.8*sin(\x)});
\draw[black] (1.5, -.008) node[fill=white,inner sep=-0.3pt] {$-u_1^-$};
\draw [black,thick,domain=0:180] plot ({4.5 + .5*cos(\x)}, {.42 - .85*sin(\x)});
 \draw[black] (4.5, -.008) node {$-u_2^+$};
\end{scope}
    \draw[<-,thick] (1.05   , .42 + 1.5 - 3.25) to[out=30,in=-130] %
                    (0.0 + 6.85, 1. - .25 ) ;
    \draw[<-,thick] (1.5 + 7, .42 + 1.5 - 4.25) -- (1.5 + 7, - .5);
\end{tikzpicture}
\caption{Illustration for displacement interpolation of two functions $u_1,
u_2$ of arbitrary sign when $u_2^- = 0$ \eqref{eq:spmzinterp}, defined in terms
of two interpolations of the non-negative cases \eqref{eq:zerointerp}.  One
treats the positive and negative parts (top row) by applying the interpolation
for non-negative functions separately to both pairs (bottom row).}
\label{fig:fourprobs}
\end{figure}
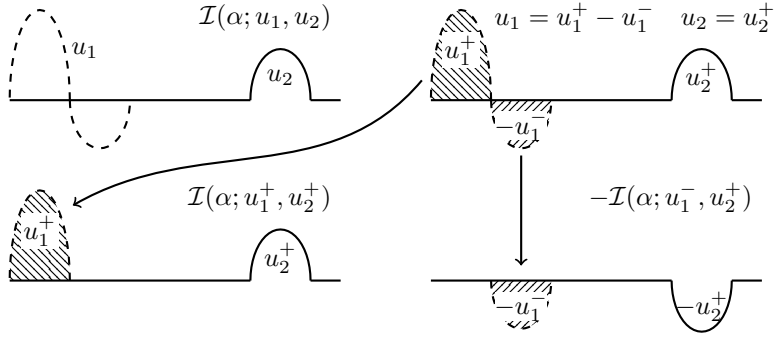

\subsection{Two functions with values of arbitrary sign}\label{sec:arb_sgn}

Now we remove the restriction that the two functions $u_1$ and $u_2$ have to be
non-negative. \added{We will separate the positive and negative parts and treat
them separately. This marks a point of departure from the optimal transport
problem \eqref{eq:ot_prob}, and our extension will not be a genuine extension
of the optimal transport solution in many aspects. However, the extension here
are intuitive and more suitable in our context of dimensionality reduction of
responses with large convective effects. This straightforward strategy was also
considered for defining the Wasserstein metric between such functions in
\cite{bonneel11,yang18}, but neither applies when any of the parts vanish, the
case considered in \cref{eq:spmzinterp} below. While this case may be
negligible in other contexts, in hyperbolic problems wave profiles often
reflect and propagate with negated sign, so it would be meaningful to construct
a well-defined interpolation even for that case. In \cite{yang18}, it was also
suggested that a fixed constant $c$ be added to $\{u_n\}$ to achieve
positivity, but this addition of low-frequency content can cause the low-rank
structures to be lost; see the concluding remarks of \cref{sec:2pfctns}.}

\begin{figure}
\centering
\begin{tabular}{ll}
    (a) & (b) \\
    \includegraphics[width=0.45\textwidth]{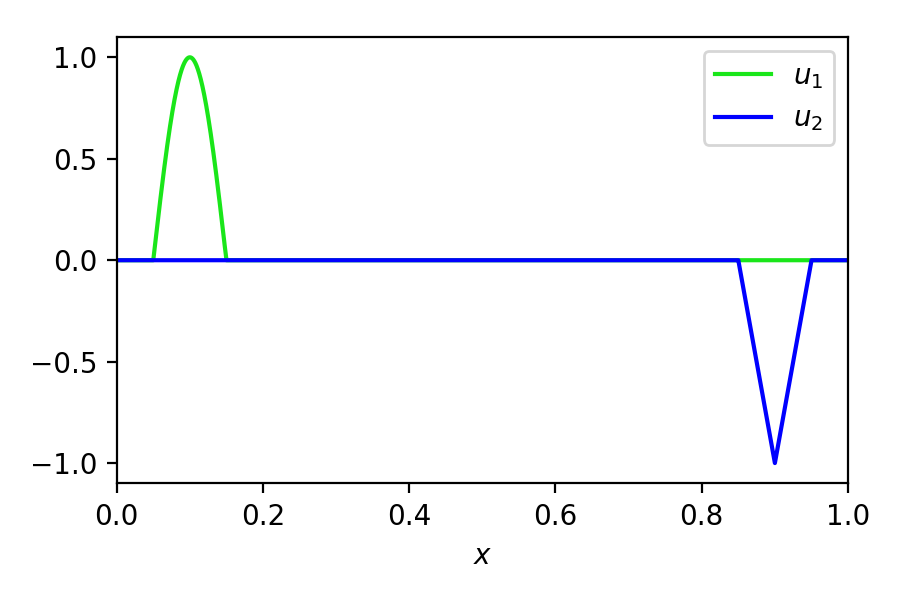}
    &
    \includegraphics[width=0.45\textwidth]{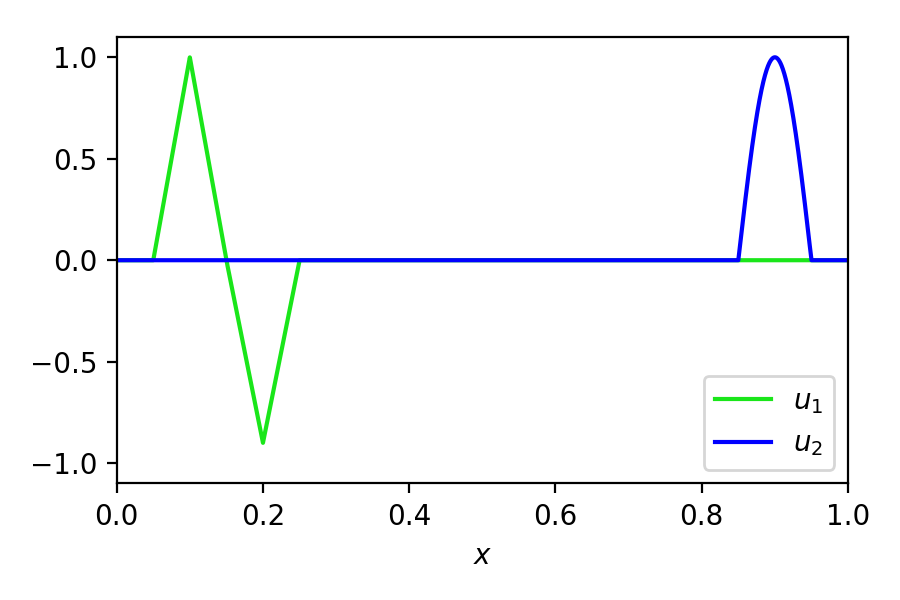}
    \\
    (c) & (d) \\
    \includegraphics[width=0.45\textwidth]{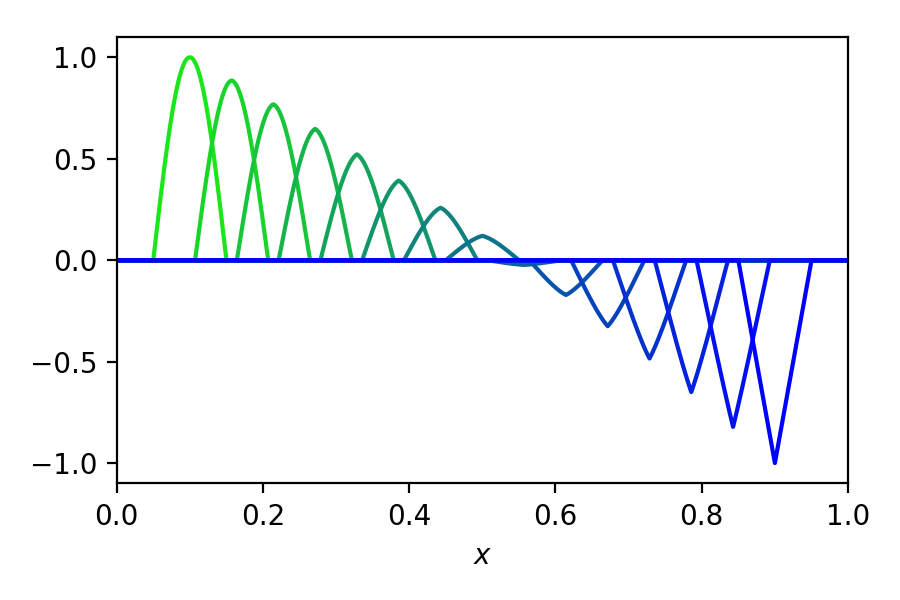}
    &
    \includegraphics[width=0.45\textwidth]{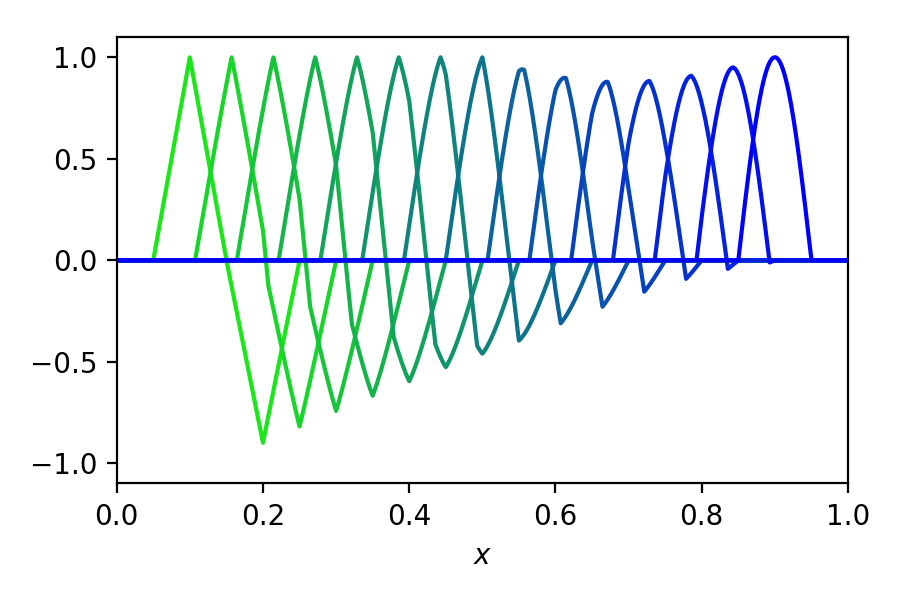}
\end{tabular}
\caption{\added{Two functions $u_1$ and $u_2$ and their displacement
interpolant defined by the interpolation procedure \cref{eq:spmzinterp}. The
case when $u_1^-$ and $u_2^+$ vanish (a,c), and the case only $u_2^-$ vanishes
(b,d).}}
\label{fig:dinterp_sgn}
\end{figure}

Let us denote the postive and negative parts of $u_1,u_2$ as follows,
\begin{equation}
    u_n^+ := \max\{u_n, 0\}, 
    \quad 
    u_n^- := \min\{u_n, 0\},
    \quad 
    \text{ for }
    n = 1,2.
\end{equation}
If all of $u_1^+, u_1^-, u_2^+, u_2^-$ are nontrivial, we interpolate as above
for $u_n^+$ and $u_n^-$ separately as in \eqref{eq:sinterp} to obtain
$\tilde{u}_\alpha^+$ and $\tilde{u}_\alpha^-$,
\begin{equation}
    \tilde{u}_\alpha^+ := \cI(\alpha;u_1^+, u_2^+)
    \quad \text{ and } \quad
    \tilde{u}_\alpha^- := \cI(\alpha;u_1^-, u_2^-),
    \label{eq:zerointerp}
\end{equation}
and then take their difference to be the interpolant,
\begin{equation}
    \tilde{u}_\alpha
    := \cI(\alpha;u_1,u_2) := \tilde{u}_\alpha^+ - \tilde{u}_\alpha^-.
    \label{eq:spminterp}
\end{equation}
If any of $u_1^+, u_1^-, u_2^+, u_2^-$ vanishes, one uses the part of opposite
sign to interpolate. As an example, suppose $u^{-}_2 = 0$ but all the other 
parts are nonzero, we compute $\tilde{u}_\alpha^-$
as the interpolant between $u_1^-$ and $-u_2^+$,
\begin{equation}
    \tilde{u}_\alpha^+ := \cI(\alpha;u_1^+, u_2^+)
    \quad \text{ and } \quad
    \tilde{u}_\alpha^- := \cI(\alpha;u_1^-, u_2^+),
\end{equation}
but we combine them so that the integral of $\tilde{u}_\alpha$ is a linear
interpolant,
\begin{equation}
    \tilde{u}_\alpha
    := \cI(\alpha;u_1,u_2) 
    := \tilde{u}_\alpha^+ - \beta \tilde{u}_\alpha^-.
    \label{eq:spmzinterp}
\end{equation}
The coefficient $\beta$ is determined by imposing that the integral of the
interpolant should be a linear interpolant between the integral of $u_1$ and
$u_2$ (that is, $M_1$ and $M_2$, respectively) as was done in
\eqref{eq:sinterp},
\begin{equation}
\int_{-\infty}^\infty \tilde{u}_\alpha\, dx
= \int_{-\infty}^\infty (\tilde{u}_\alpha^+ - \beta \tilde{u}_\alpha^-) \, dx
    = (1 - \lambda(\alpha)) M_1 + \lambda(\alpha) M_2.
\end{equation}
Then one may easily compute $\beta$,
\begin{equation}
    \beta = \frac{(1-\lambda(\alpha)) M_1^- + \lambda(\alpha) M_2^-}
                  {(1-\lambda(\alpha)) M_1^- + \lambda(\alpha) M_2^+}
    \quad \text{ where } \quad
    M_n^{\pm} := \int_{-\infty}^\infty u^{\pm}_n \,dx
    \quad \text{ for } n = 1,2.
\end{equation}
\added{
Since $M_2^- = 0$ in this case, we have that
\begin{equation}
    \beta = \frac{(1-\lambda(\alpha)) M_1^-}
                  {(1-\lambda(\alpha)) M_1^- + \lambda(\alpha) M_2^+}.
\end{equation}
}
Analogous definitions follow for other cases (when more of $u_n^{\pm}$ vanish.)
\added{Two representative examples are shown in \cref{fig:dinterp_sgn}}.

\added{The extension suggested here appears natural for our applications, but
it is by no means the only possible one. We also note that unlike the
monotone rearrangement, the resulting interpolant defined here for arbitrary
signs will not preserve the monotonicity in general. However, since we are
simply interested in finding a low-dimensional structure in the transport maps
thusly defined, the lack of monotonicity is not necessarily a serious concern.

On the other hand, it may be advantageous to adapt the interpolation procedure
above by decomposing the derivatives of $u_1$ and $u_2$ into individual
\emph{pieces} with connected supports and apply the interpolation above between
individual pairs, in relation to scalar conservation laws \cite{bolley05}. This
interpolation procedure was called \emph{displacement interpolation by pieces}
and used in the construction of ROMs in \cite{rim18}. This map will preserve
monotonicity if the given functions satisfy the \emph{signature condition}
defined therein. We will show in \cref{sec:transform} that these extensions are
merely special cases of a more general formulation.}

\subsection{More than two functions}\label{sec:2pfctns}

Consider the case when we are given multiple functions, corresponding to the
parametrized function evaluated at multiple parameter values,
\[
\cU = \{u(\balpha_1), u(\balpha_2), ..., u(\balpha_n)\}.
\]
\added{For simplicity of exposition, let us assume that all functions in $\cU$
are non-negative}.  We wish to interpolate these functions with respect to the
parameters $\balpha$.  Let us compute $U_n^\dagger$ for each $u_n =
u(\balpha_n)$ by the formula \eqref{eq:pinv}, and denote the set of these CDFs
by $\cU^\dagger = \{U_n^\dagger: n = 1, ... , N\}$.  Then we can extend the
interpolation procedure, by generalizing the definition of the CDF of the
interpolant \eqref{eq:interp2} in a piece-wise linear manner,
\begin{equation}
    \widetilde{U}_{\balpha}^\dagger 
    := \sum_{n=1}^N \lambda_n(\balpha) \, U^\dagger_n.
    \label{eq:bary_dinterp}
\end{equation}
There is a freedom in choosing $\lambda_n(\balpha)$ but we will choose it to be
the barycentric coordinates for conceptual simplicity. \added{For more general
discussions on Wasserstein barycenters, see \added{\cite{agueh11}}.} That is,
$\{\balpha_n\}$ will serve as nodes (or vertices) of a tessellation, and
$\lambda_n(\balpha)$ will yield the barycentric coordinate with respect to the
nodes of the polytope $\balpha$ belongs to.

\begin{figure}
\centering
\includegraphics[width=0.45\textwidth]{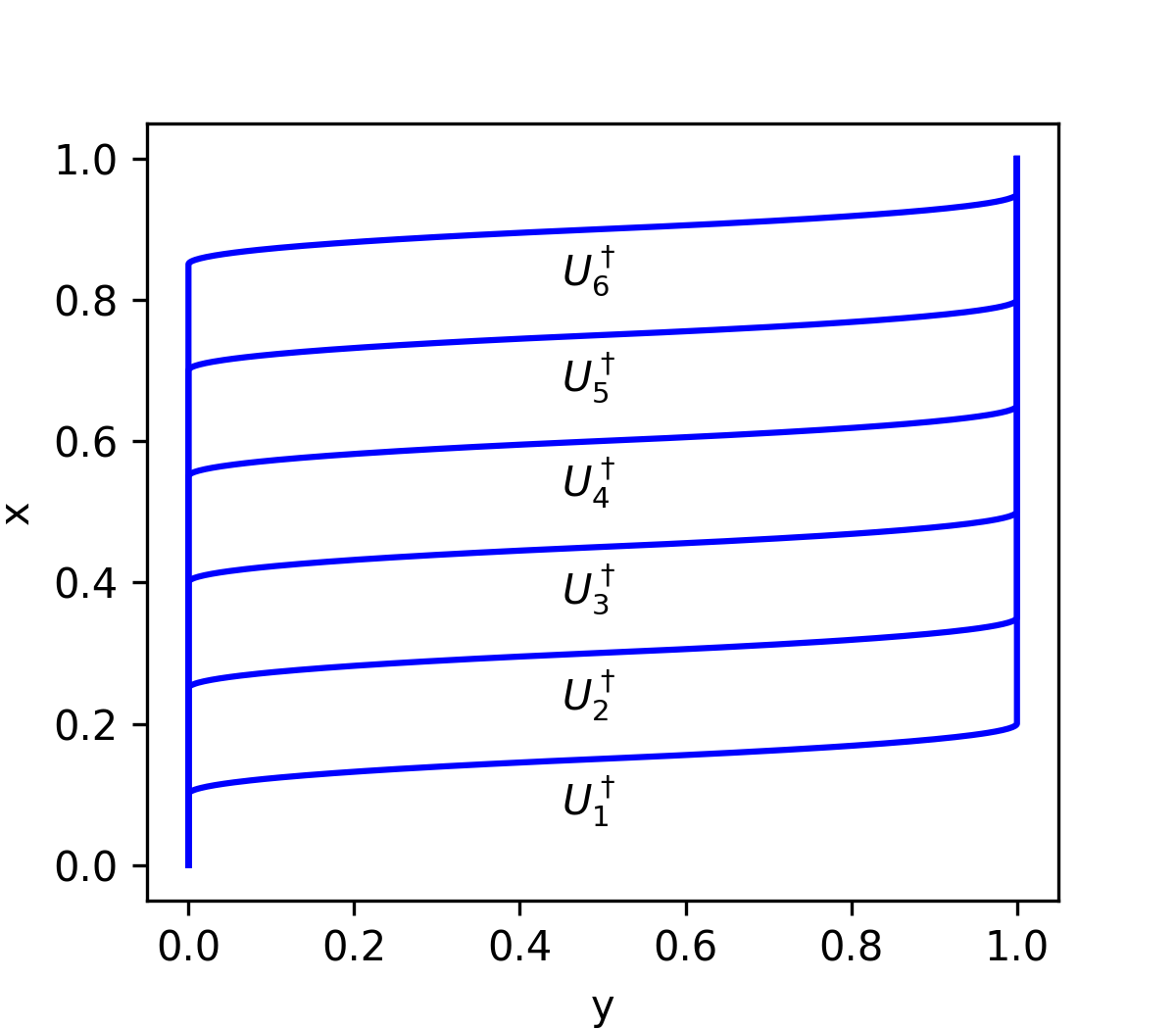}
\caption{Plot of $U_n^\dagger$ corresponding to the six hat functions to the
right of \cref{fig:hats}.}
\label{fig:hats_CDFs}
\end{figure}

Further suppose that $\cU^\dagger$ has a low-rank representation with the
corresponding low-dimensional basis $\cV^\dagger = \{V^\dagger_n: n = 1, ...
M\}$.  Then, the interpolant can be computed in a similar manner as was done
for the linear approximation \eqref{eq:linterp},
\begin{equation}
    \widetilde{U}_{\textrm{LR},\balpha}^\dagger 
    := \sum_{n=1}^N \nu_n(\balpha) \, V^\dagger_n.
    \label{eq:dinterp_lr}
\end{equation}
for coefficients $\nu_n$ that depend on $\balpha$.  Therefore we obtain a
low-dimensional representation $\widetilde{U}_{\textrm{LR},\balpha}^\dagger$ in
terms of the optimal transport map.  Then $\tilde{u}_{\balpha}$ is given by the
analogue of \eqref{eq:sinterp},
\begin{equation}
    \tilde{u}_{\balpha} :=
        \cI(\balpha; u_1, \cdots, u_N)
     := \left[ \sum_{n=1}^N \nu_n(\balpha) M_n \right]
        \frac{d}{dx} \tilde{U}_{\textrm{LR},\balpha}
    \label{eq:interpn}
\end{equation}
\added{Straightforward adjustments to \cref{eq:bary_dinterp,eq:dinterp_lr} can
be made for the case when the functions take on arbitrary signs. The basis
functions $\{V_n\}$ in \cref{eq:dinterp_lr} were named \emph{advection modes}
in \cite{iollo14}, although the modes here can be computed explicitly without
any optimization.} In order to demonstrate that the interpolant
$\tilde{u}_{\balpha}$ is useful, let us revisit the 1D transport example
\eqref{eq:transport1d}.  Suppose we were given the set of solutions to the
transport equation,
\[
\cU = \{u(t_n): n=1, ... , 6\}
\quad \text{ where } t_n = 3w (n-1), \quad (w = 0.05),
\]
as was shown in \cref{fig:hats}.  If one computes the CDFs $\{ U_n \}$ for each
$u_n = u(t_n)$, as shown in \cref{fig:hats_CDFs}, one discovers that their
pseudo-inverses are 
\beq
    U_n^{\dagger}(y) = 3w(n-1) H(y) + U_1^{\dagger}(y).
\eeq
This shows that translation is a low-dimensional operation when viewed in the
form \eqref{eq:dinterp_lr}, representable by the addition of a Heavside
function.  That is, although direct application of low-rank approximations such
as the singular value decomposition (SVD) will not succeed in finding a
low-rank approximation directly with  $\{u_n\}$, it can be successfully applied
to the CDFs $\{U^\dagger_n\}$. \added{See \cref{sec:hats} for a related
example.}

We remark that an addition of a constant value to $u(t_j)$ can affect this
low-rank property for this problem. Just as linear operation such as the
addition by a constant in the $y$ variable result in rank-increasing operations
(translation) in the $x$ variable, linear operations in $x$ can result in an
increase of rank as functions in the $y$ variable. This implies that some
simple preprocessing of the functions may be necessary for this approach to be
successful in general, e.g, the application of these methods to the derivative
$\partial u_n / \partial x$ instead of $u_n$, \added{as is done when applying
displacement interpolation by pieces \cite{rim18}.  This was also observed in
\cite{bonneel11} where low-band modes were treated independently.} This being
said, the issue appears to be easy to circumvent and did not significantly
affect the applicability of the interpolation method in our examples.

\subsection{Multiple spatial dimensions} \label{sec:multid}

We have so far defined the displacement interpolation $\cI$
(\ref{eq:interp1def}, \ref{eq:sinterp}, \ref{eq:spmzinterp}) in a single
spatial dimension. This definition can be naturally extended to multiple
dimensions, and one approach is to consider the multi-dimensional
Monge-Kantorovich optimal transport problem \cite{Benamou}.  However, the
multi-dimensional problem itself requires a solution of a nontrivial problem,
e.g., the Monge-Amp\`ere equation or related optimization problem which makes
the resulting interpolant implicit and expensive to compute. To avoid such
hurdles we instead propose a multi-dimensional extension which makes use of the
intertwining property of the Radon transform
\cite{Helgason2011,natterer,laxphilips,radonsplit}.  \added{Note that this
extension was also explored in the image processing literature
\cite{pitie07,bonneel15} as a method for avoiding the high computational cost
of solving the multi-dimensional optimal transport problem.} 

Suppose $u$ is a function of multiple spatial dimensions (e.g. $d=2,3$). The
Radon transform of $u$ is a function of $s \in \RR$ and $\boldsymbol{\omega}
\in S^{d-1}$, defined as the integral over the hyperplane oriented by $s$ and
$\boldsymbol{\omega}$,
\begin{equation}
    \cR \left[ u \right](s,\boldsymbol{\omega}) 
    :=  \int_{\mathbf{x} \cdot \boldsymbol{\omega} = s}  
    u(\mathbf{x}) \,dm(\mathbf{x}),
    \label{eq:radon}
\end{equation}
where $m(\mathbf{x})$ is the Euclidean measure over the hyperplane
$\{\mathbf{x} \in \RR^d: \mathbf{x} \cdot \boldsymbol{\omega} = s\}$.

The main property of this transform we will be interested in is the
intertwining property \cite{Helgason2011,natterer}, 
\beq
    \cR \left[ \frac{\partial u}{\partial x_i }\right] (s,\boldsymbol{\omega})
    = \omega_i \left(\frac{\partial }{\partial s}   \cR \left[u \right]
    \right)
    (s,\boldsymbol{\omega}),
    \label{eq:intertwine}
\eeq
where $\omega_i$ is the $i$-th component of $\boldsymbol{\omega}$. 
For example, the multi-dimensional transport equation
\[
    u_t + \boldsymbol{\theta} \cdot \nabla u = 0
    \quad \text{ where }  \boldsymbol{\theta} \in S^{d-1},
\]
for $u$ is transformed into a collection of 1D transport equations for $\cR
[u]$, where the transport speed for each 1D problem now depends on
$\boldsymbol{\omega}$,
\begin{equation}
    \cR\left[ u \right]_t 
    + (\boldsymbol{\theta} \cdot \boldsymbol{\omega})\, \cR\left[ u \right]_s = 0.
\end{equation}
This property is useful not only for the transport equation, but other
hyperbolic PDEs such as the acoustics equations. We refer to \cite{radonsplit}
for further details.

A natural extension of the displacement interpolation procedure $\cI$ in 1D
defined in the previous sections is to perform the interpolation for each fixed
$\boldsymbol{\omega}$ in the transformed variable $s$.  That is, given two
multi-dimensional PDFs $u_1, u_2 \ge 0$, we define the displacement interpolant
in terms of its Radon transform,
\begin{equation}
\cR \left[\tilde{u}_\alpha\right](\cdot,\boldsymbol{\omega}) 
    := \cI\left(\alpha; \cR \left[ u_1 \right](\cdot,\boldsymbol{\omega}), 
                    \cR \left[ u_2 \right](\cdot,\boldsymbol{\omega})\right)
\label{eq:dinterp-slice}
\end{equation}
\added{We will also denote the interpolation that appears on the RHS of
\cref{eq:dinterp-slice} by
\beq
  \cI_{\otimes}
        \left(\alpha; \cR \left[ u_1 \right], \cR \left[ u_2 \right]\right),
\eeq
where $\cI_\otimes$ acts on the $s$-variable of the transform individually.
We now} invert the transform to obtain $\tilde{u}_\alpha$,
\begin{equation}
    \tilde{u}_\alpha 
    := \cI_d(\alpha; u_1,u_2)
    := \cR^{-1} \left[ \cR[\tilde{u}_\alpha ]\right].
    \label{eq:ndinterp}
\end{equation}
\added{The inversion $\cR^{-1}$ is ill-posed, but unlike in the more common
tomography setting one can reduce the inversion error simply by over-sampling
the forward data. For a more detailed discussion, we refer the reader to
\cite{radonsplit,rim17thesis}}. In \cref{sec:numerics}, we will see that this
explicit construction yields reasonable numerical results. We will call this
\added{interpolation} operator the \emph{generalized Lax-Philips operator}, as
it generalizes the translation representation of \added{the one-parameter
family of solutions to the wave equation by} Lax-Philips \cite{laxphilips} to
representations using monotone rearrangement.

To the best of our knowledge, it has not yet been rigorously investigated
whether this particular construction $\tilde{u}_\alpha$ is a solution to a
multi-dimensional optimal transport problem, perhaps a variant of
\eqref{eq:ot_prob}. We will not pursue such an undertaking here, and leave the
task to a future work.

\added{
\subsection{Composition with general transforms} \label{sec:transform} The
application of the Radon transform in \cref{sec:multid} suggests that the
interpolant can achieve a rich set of behavior when applied together with the
right transform.  We formally define such a displacement interpolation
procedure in this section. Let us denote an invertible transform by $\cT$.  We
define $\cI_\cT:\cL \times \cL \to \cL$ by a conjugation of $\cI_\otimes$,
\beq
\cI_\cT := \cT^{-1} \cI_\otimes \cT.
\label{eq:dinterp_transform}
\eeq
\added{Note that the dependence on the interpolation parameter $\alpha$ is
omitted here.} If the operands are 1D functions, let $\cI_\otimes = \cI$.

As an elementary example, the operation of separating the positive and negative
parts as in \cref{sec:arb_sgn} can be formulated differently. For the case when
neither the positive nor the negative parts vanish, the transform $\cT$ may be
given by
\beq
\cT [u] = [\max\{u,0\},\max\{-u,0\}]
\quad \text{ and } \quad
\cT^{-1} [u^+,u^-] = u^+ - u^-.
\eeq
Similar definitions could be made for the other cases discussed in
\cref{sec:arb_sgn}. Displacement interpolation by pieces in \cite{rim18} could
also be expressed in this form. 

We will remark on how \cref{eq:dinterp_transform} serves to unify the various
interpolants based on the monotone rearrangement in 1D that appear here and
elsewhere under one framework. The generalized Lax-Philips operator
\cref{eq:ndinterp} is a special case of the above when one lets $\cT = \cR$. We
also note that the operator $\cI_\otimes$ could have been applied to the
$\boldsymbol{\omega}$-variable as an alternative. One obvious choice for $\cT$
could be the derivative operator; see \cref{sec:1dacoustics}. In image
processing, choosing $\cT$ to be the Fourier transform or wavelet transform was
shown to be useful for color or texture mixing \cite{bonneel15,rabin12}. The
treatment of low-frequency components in \cite{bonneel11} can also be related
to by \cref{eq:dinterp_transform}.

In our numerical example in \cref{sec:osc}, we will show that $\cT= \cR \otimes
\cP \cF$ for some permutation $\cP$ can capture the low-rank structure of
oscillatory functions.}

\section{Numerical implementation and examples} \label{sec:numerics}

In this section, we discuss some issues related to implementation and provide
some numerical examples that illustrate the behavior of the displacement
interpolant introduced in the previous sections.

\subsection{Implementation}\label{sec:implement}
In the discretized setting, the function $u_n$ in the previous sections will be
represented as a vector $\bu_n$ in $\RR^D$.  We will choose a piece-wise
constant representation of the function $u_n$.  The reason for this is
simplicity, as the exact cumulative distribution function $U_n$ defined in
\eqref{eq:CDFs} will then be piece-wise linear.  This is convenient since we
plan to perform linear operations on $U_n^\dagger$, e.g., to interpolate
between two profiles in \eqref{eq:interp2}, or to compute low-rank linear bases
for multiple profiles in \eqref{eq:dinterp_lr}.  A higher-order representation
of $u_n$ will cause the exact representation of $U_n^\dagger$ to involve
fractional powers of $y$ (such as the square root) which is cumbersome to work
with, especially when the grid-points are not uniform in terms of the variable
$U_n$.

To be more precise, let us suppose that $\bu_n$ represents $u_n$ by its
integral over cells of uniform width $h$, i.e., $x_{j+1/2} = jh$,
\[
    (\bu_n)_j = \frac{1}{h}\int_{x_{j-1/2}}^{x_{j+1/2}} u_n\, dx.
\]
Now, let us denote by $\bfU_n$ the cumulative distribution of $\bu_n$ when
viewed in its piece-wise constant representation. The result will be a
piece-wise linear function and will serve as a discretization of $U_n$, a
linear interpolant through the points $\sum_{k=1}^j(\bu_n)_k$ on the same
uniform grid $\{x_j\}$. 

To compute the pseudo-inverse $U^\dagger_n$ one simply exchanges the two
entries of the interpolating points. That is 
\[
(x_{j+1/2},(\bfU_n)_j) 
\quad \to \quad 
((\bfU_n)_j, x_{j+1/2}). 
\]
We will define $\bfU_n^\dagger$ as the linear interpolant of these points.
When represented as a 2D array, for example, we can write
\[
    (\bfU_n^\dagger)_{j,1} = (\bfU_n)_j
    \quad \text{ and } \quad
    (\bfU_n^\dagger)_{j,2} = x_{j+1/2}.
\]
However, although $(\bfU_n)_j$ is non-increasing, it does not form a uniform
grid in general, and there will be redundant entries if $\bu_n$ vanishes.  So
some care must be taken when computing linear operations such as
\eqref{eq:interp2} on  $\bfU_n^\dagger$. One simple solution would be to merge
the grid $(\bfU_n)_j$ whenever two pseudo-inverses are added, and in the worst
case the size of the array will grow linearly with the number of operations. We
will make this choice for the implementation presented here as the
interpolation procedure will be exact and free from approximation errors, but
one may opt for a fixed grid interpolation away from the points in $\cY$ where
$U_n$ is constant \eqref{eq:cY} for efficiency.  \added{The implementation we
have used is available in a public online repository \cite{dinterp}.}

To compute the Radon transform \eqref{eq:radon}, we will make use of a  fast
approximate algorithm called the Discrete Radon Transform (DRT) %
\cite{drt,pressdrt,radonsplit}. To control the inversion error, the inverse is
computed using the conjugate gradient algorithm \cite{greenbaum} on the
prolongated transform. While the computational cost for the inversion could be
high, it is still conjectured to be $\cO(D^{5/2}\log D)$ when the 2-dimensional
functions $u_n$ are each represented on the grid of size $D \times D$ (so that
$\bu_n$ is in $\RR^{D \times D}$).  Further improvements in the inversion
algorithms may reduce this computational burden, but this topic is also left as
future work.

\subsection{1D interpolation with 2D parameters}\label{sec:2param}

\begin{figure}
    \centering
\begin{tikzpicture}[scale=1.9]
\node[inner sep=0pt] (u1) at (-1.05, -1.10)
    {\includegraphics[width=.23\textwidth]{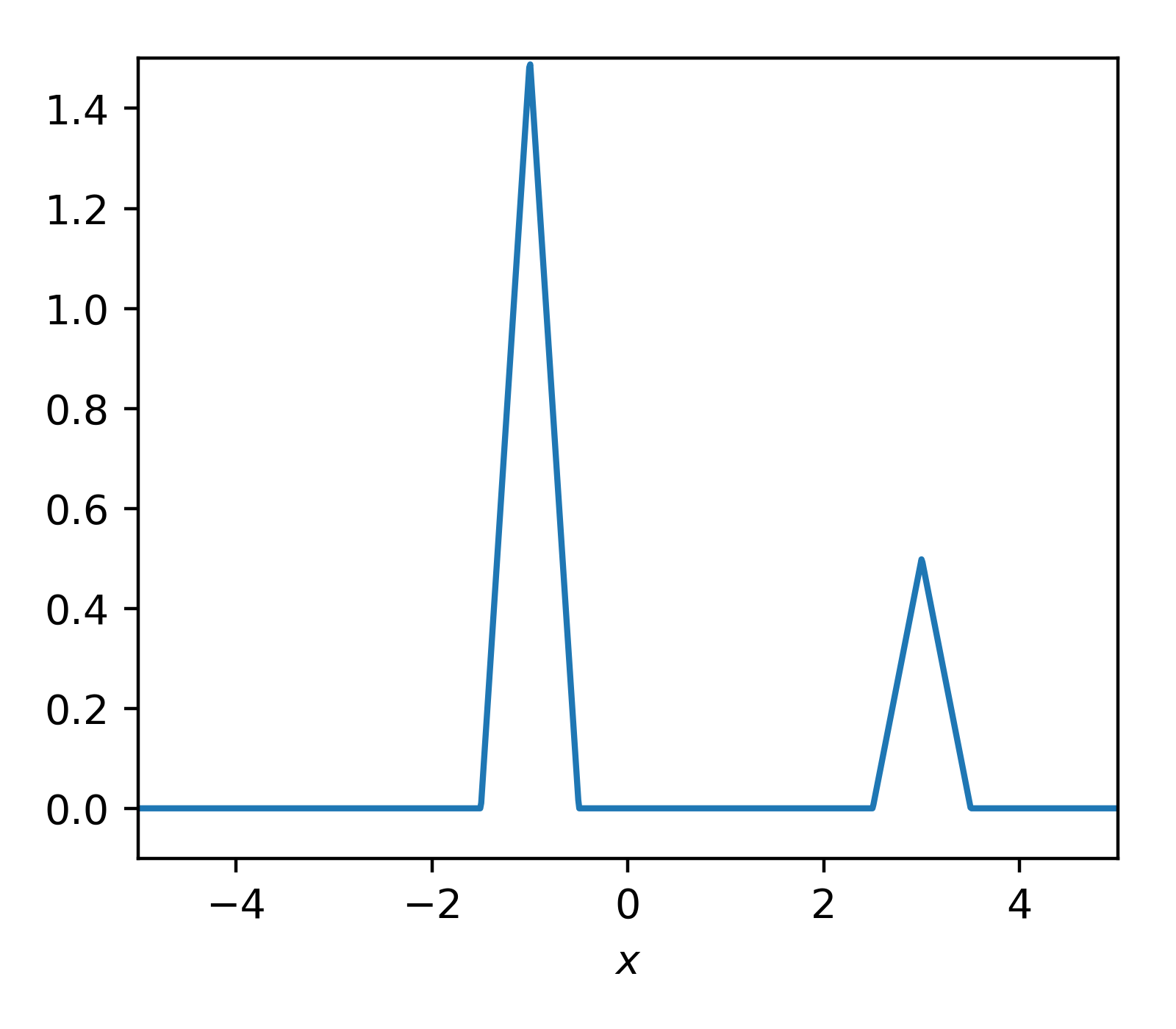}};
\node[inner sep=0pt] (u2) at ( 1.30, -1.10)
    {\includegraphics[width=.23\textwidth]{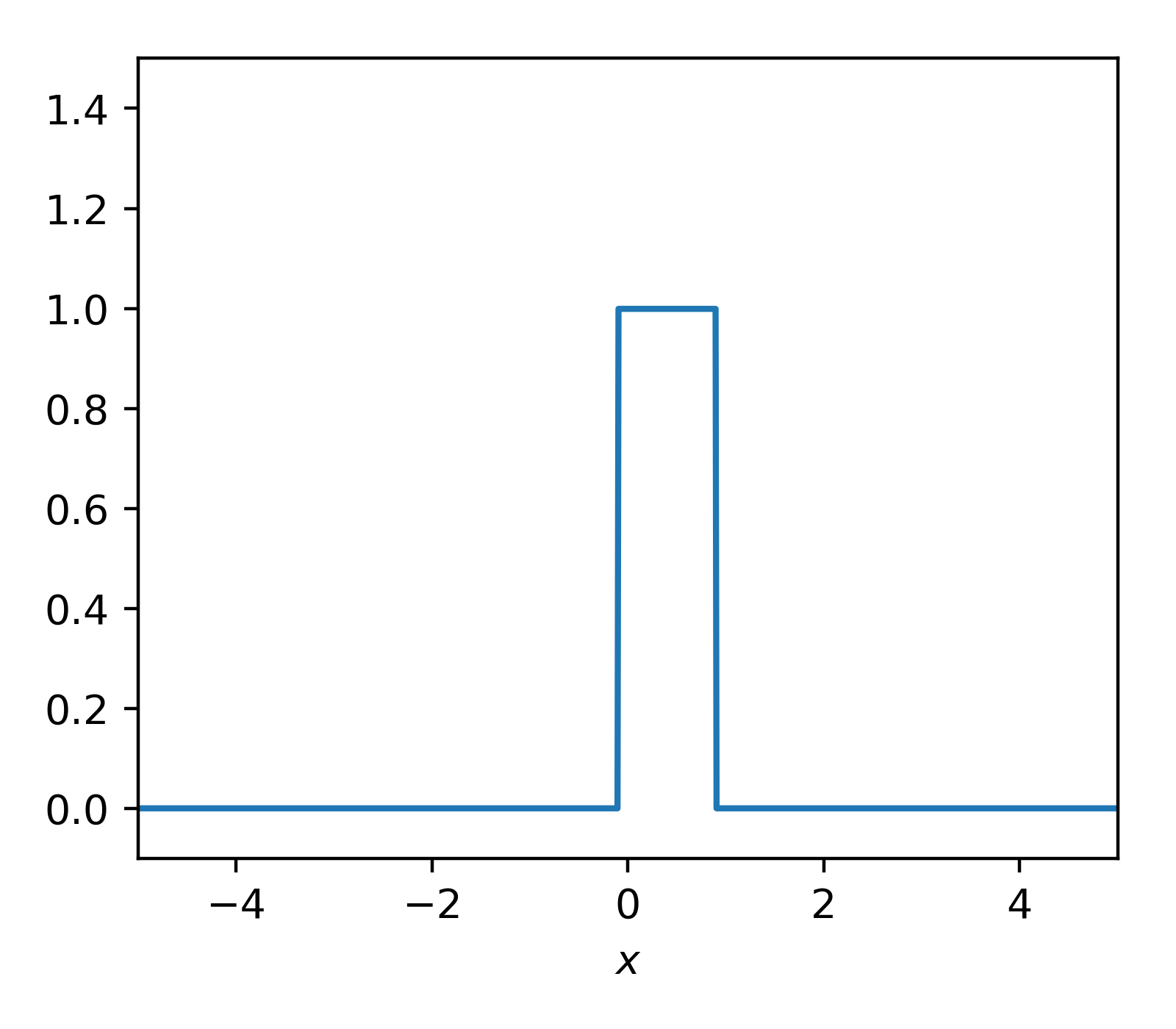}};
\node[inner sep=0pt] (u3) at (-1.05, 1.00)
    {\includegraphics[width=.23\textwidth]{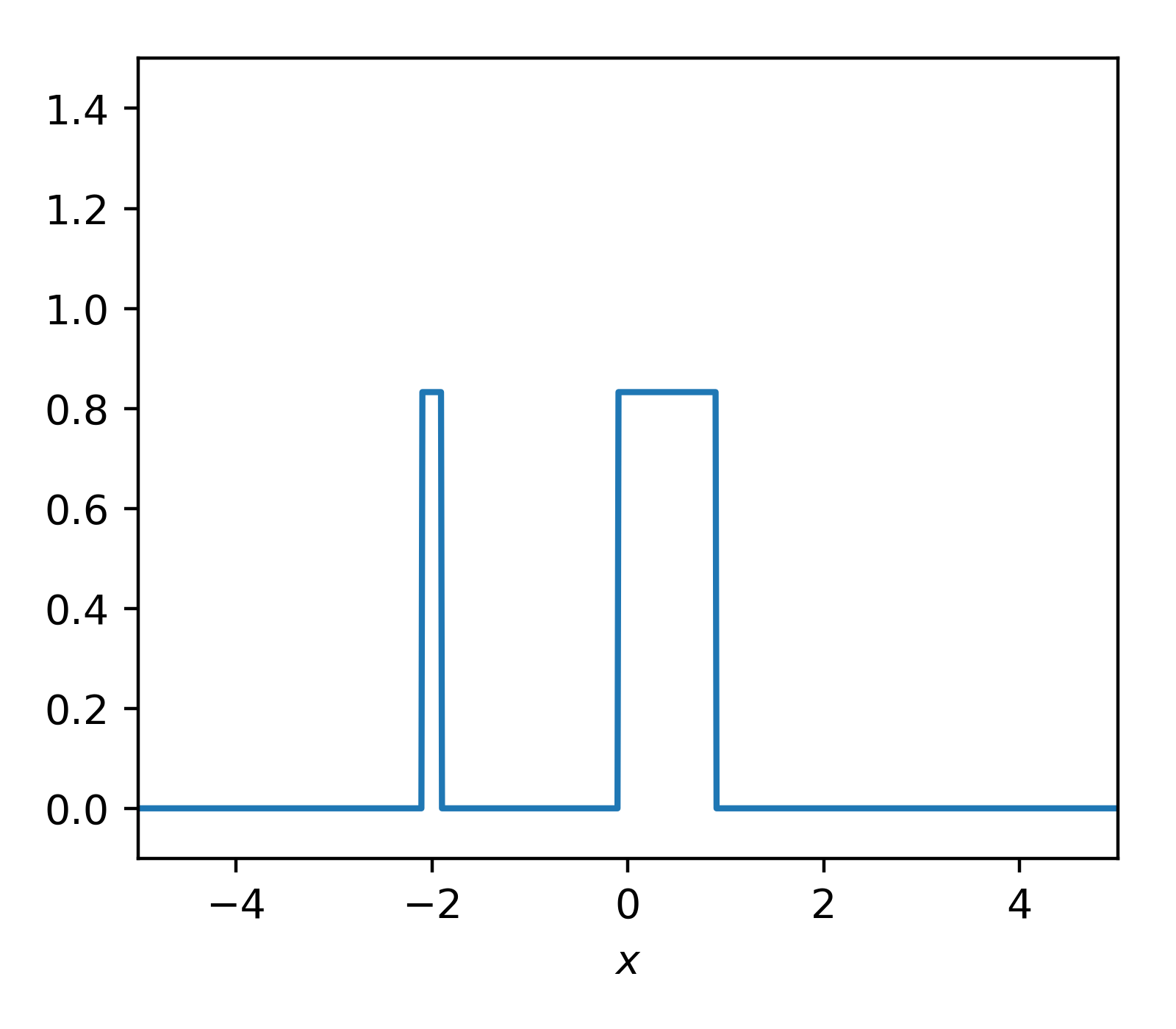}};
    \draw[thick] (u1) to[out= 40,in=-90] (0,0);
    \draw[thick] (u2) to[out=100,in=-90] (1,0);
    \draw[thick] (u3) to[out=  0,in=190] (0,1);
\draw[fill=gray!20] (0,0) -- (1,0) -- (0,1) -- (0,0);
\draw[->] (0,0) -- (0  ,1.25);
\draw[->] (0,0) -- (1.25,0  );
\draw (0,0) node[anchor=east]  {$\balpha_1$};
\draw (1,0) node[anchor=north west]   {$\balpha_2$};
\draw (0,1) node[anchor=south west ]   {$\balpha_3$}; 
\draw[fill=black] (0.  ,0.  ) circle (0.025cm);
\draw[fill=black] (1.  ,0.  ) circle (0.025cm);
\draw[fill=black] (0.  ,1.  ) circle (0.025cm);
\draw[fill=white] (0.25,0.  ) circle (0.025cm) node[anchor=south west] {$1$};
\draw[fill=white] (0.50,0.  ) circle (0.025cm) node[anchor=south west] {$2$};
\draw[fill=white] (0.75,0.  ) circle (0.025cm) node[anchor=south west] {$3$};
\draw[fill=white] (0.  ,0.25) circle (0.025cm) node[anchor=south west] {$9$};
\draw[fill=white] (0.  ,0.50) circle (0.025cm) node[anchor=south west] {$8$};
\draw[fill=white] (0.  ,0.75) circle (0.025cm) node[anchor=south west] {$7$};
\draw[fill=white] (0.25,0.25) circle (0.025cm) node[anchor=south west] {$10$};
\draw[fill=white] (0.50,0.25) circle (0.025cm) node[anchor=south west] {$11$};
\draw[fill=white] (0.25,0.50) circle (0.025cm) node[anchor=south west] {$12$};
\draw[fill=white] (0.75,0.25) circle (0.025cm) node[anchor=south west] {$4$};
\draw[fill=white] (0.50,0.50) circle (0.025cm) node[anchor=south west] {$5$};
\draw[fill=white] (0.25,0.75) circle (0.025cm) node[anchor=south west] {$6$};
\end{tikzpicture} 
\caption{The three functions $u_1=u(\balpha_1)$,$u_2=u(\balpha_2)$,
$u_3=u(\balpha_3)$, corresponding to parameter values $\balpha_1=(0,0)$
$\balpha_2=(1,0)$, $\balpha_3=(0,1)$. The convex hull of $\{\balpha_1,
\balpha_2,\balpha_3\}$ is also displayed along with nodes where the
inteprolants will be computed enumerated 1 to 12.  The computed interpolants
are shown in subsequent \cref{fig:interp2param123,fig:interp2param456,%
fig:interp2param789,fig:interp2param101112}.}
\label{fig:2paramfctns}
\end{figure}
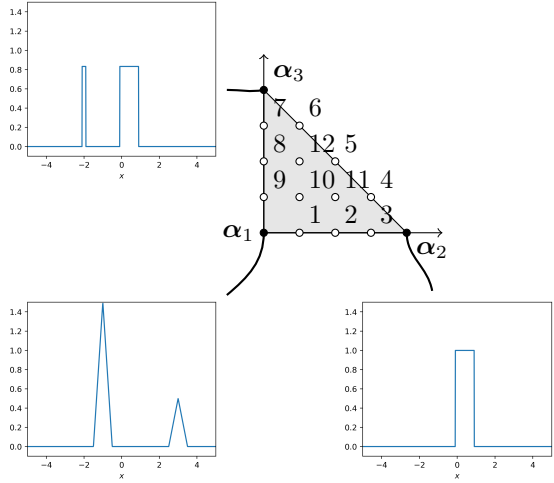

\begin{figure}
\centering
\includegraphics[width=0.80\textwidth]{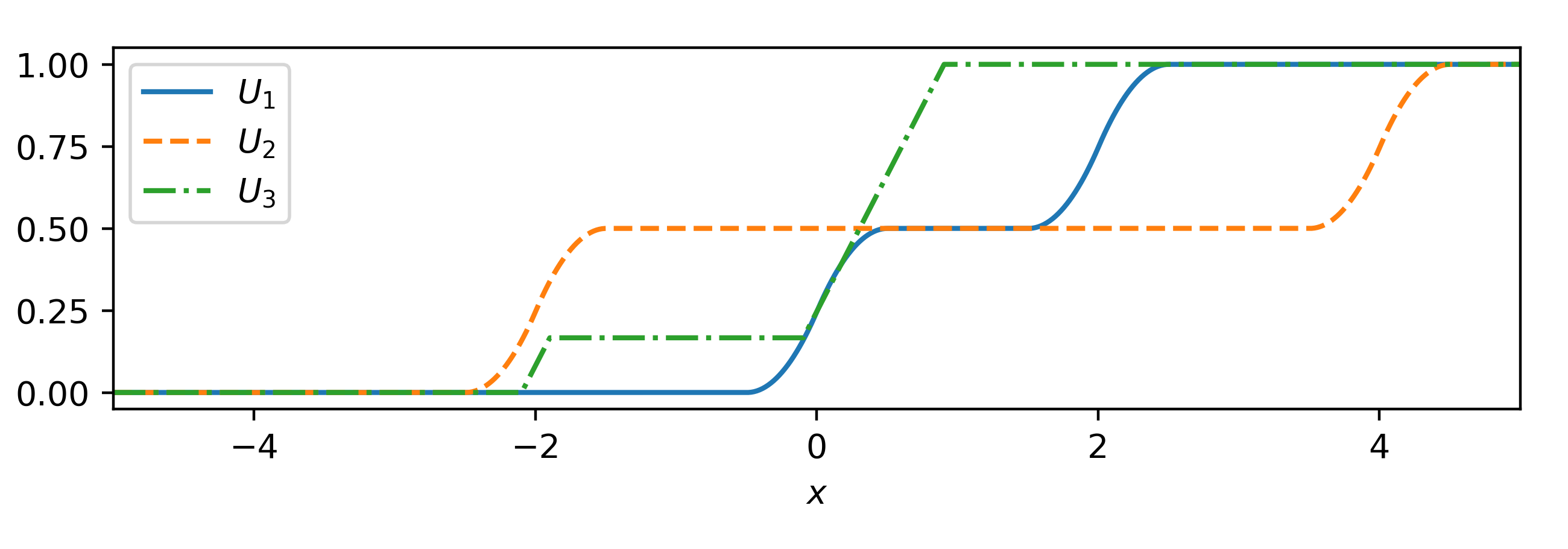}
\caption{The cumulative distribution functions (CDFs) of the functions $u_1$,
$u_2$ and $u_3$, denoted by $U_1, U_2$ and $U_3$, respectively.}
\label{fig:2paramsCDFs}
\end{figure}

In this section, we will illustrate the behavior of the displacement
interpolant defined above by computing the interpolant between three different
functions $u_1, u_2, u_3$ functions, using the barycentric formula. These three
functions are shown in \cref{fig:2paramfctns}.  The function $u_1$ is a sum of
two hat functions $\phi$ \eqref{eq:hatf} with different heights and positions,
and the function $u_2$ is a sharp hump which can be represented by a sum of two
Heaviside jump functions.  The function $u_3$ is a superposition of two
functions in the form of $u_2$.  These functions are neither translates of each
other nor do they have identical number of connected supports. In other words,
if one recalls the definition of $\cY_n$ \eqref{eq:cY}, we have that
$|\cY_1|=|\cY_3|  = 3$ whereas $|\cY_2| = 2$. Also, total mass is distributed
differently between the two superposed profiles in $u_1$ and $u_3$.  These
functions are shown in \cref{fig:2paramfctns}.  The corresponding CDFs are
shown in \cref{fig:2paramsCDFs}.

The displacement interpolant we will compute will be of the form
\eqref{eq:bary_dinterp}, where for the coefficients $c$ we will use the
barycentric coordinates $\lambda$. We will suppose that
\[
u_1 = u(\balpha_1), \quad
u_2 = u(\balpha_2), \quad
u_3 = u(\balpha_3),
\]
where the parameter values are given by
\[
\balpha_1 = (0,0), \quad
\balpha_2 = (1,0), \quad
\balpha_3 = (0,1). 
\]
Then the displacement interpolant is written as
\begin{equation}
 \tilde{u}_{\balpha} = \cI_{\lambda(\balpha)}(\balpha; u_1,u_2,u_3),
    \label{eq:nbarycentric}
\end{equation}
whose pseudo-inverse of the CDF $\widetilde{U}^\dagger$ is computed by
\beq
\begin{aligned}
    \widetilde{U}^\dagger_{\balpha }
    &= \lambda_1(\balpha) U_1^\dagger
      + \lambda_2(\balpha) U_2^\dagger
      + \lambda_3(\balpha) U_3^\dagger \\
    &= (1 - \alpha_1 - \alpha_2) U_1^\dagger
      + \alpha_1 U_2^\dagger
      + \alpha_2 U_3^\dagger .
\end{aligned}
\label{eq:2dquantile}
\eeq

We will evaluate $\tilde{u}_{\balpha}$ at the nodes in the convex hull 
of $\{\balpha_1, \balpha_2, \balpha_3\}$ shown as
white circles in \cref{fig:2paramfctns} with respective enumeration. 
The nodes are,
\begin{equation}
    \begin{aligned}
    \tilde{\balpha}_1 &= (0.25,0), \\
    \tilde{\balpha}_2 &= (0.5,0), \\
    \tilde{\balpha}_3 &= (0.75,0), \\
    \end{aligned}
    \quad
    \begin{aligned}
    \tilde{\balpha}_4 &= (0.75,0.25), \\
    \tilde{\balpha}_5 &= (0.5,0.5), \\
    \tilde{\balpha}_6 &= (0.25,0.75), \\
    \end{aligned}
    \quad
    \begin{aligned}
    \tilde{\balpha}_7 &= (0,0.75), \\
    \tilde{\balpha}_8 &= (0,0.5), \\
    \tilde{\balpha}_9 &= (0,0.25). \\
    \end{aligned}
    \quad
    \begin{aligned}
    \tilde{\balpha}_{10} &= (0.25,0.25), \\
    \tilde{\balpha}_{11} &= (0.5,0.25), \\
    \tilde{\balpha}_{12} &= (0.25,0.5). \\
    \end{aligned}
    \label{eq:talpha}
\end{equation}

\begin{figure}
    \centering
\begin{tikzpicture}[scale=1.9]
\node[inner sep=0pt] (v1) at (-1.10, 1.30-0.75)
    {\includegraphics[width=0.23\textwidth]{exp3_dinterp_00_00.png}};
\node[inner sep=0pt] (v2) at ( 2.50, 1.30-0.75)
    {\includegraphics[width=0.23\textwidth]{exp3_dinterp_10_00.png}};
\node[inner sep=0pt] (u1) at (-1.10, 1.30-2.05)
    {\includegraphics[width=0.23\textwidth]{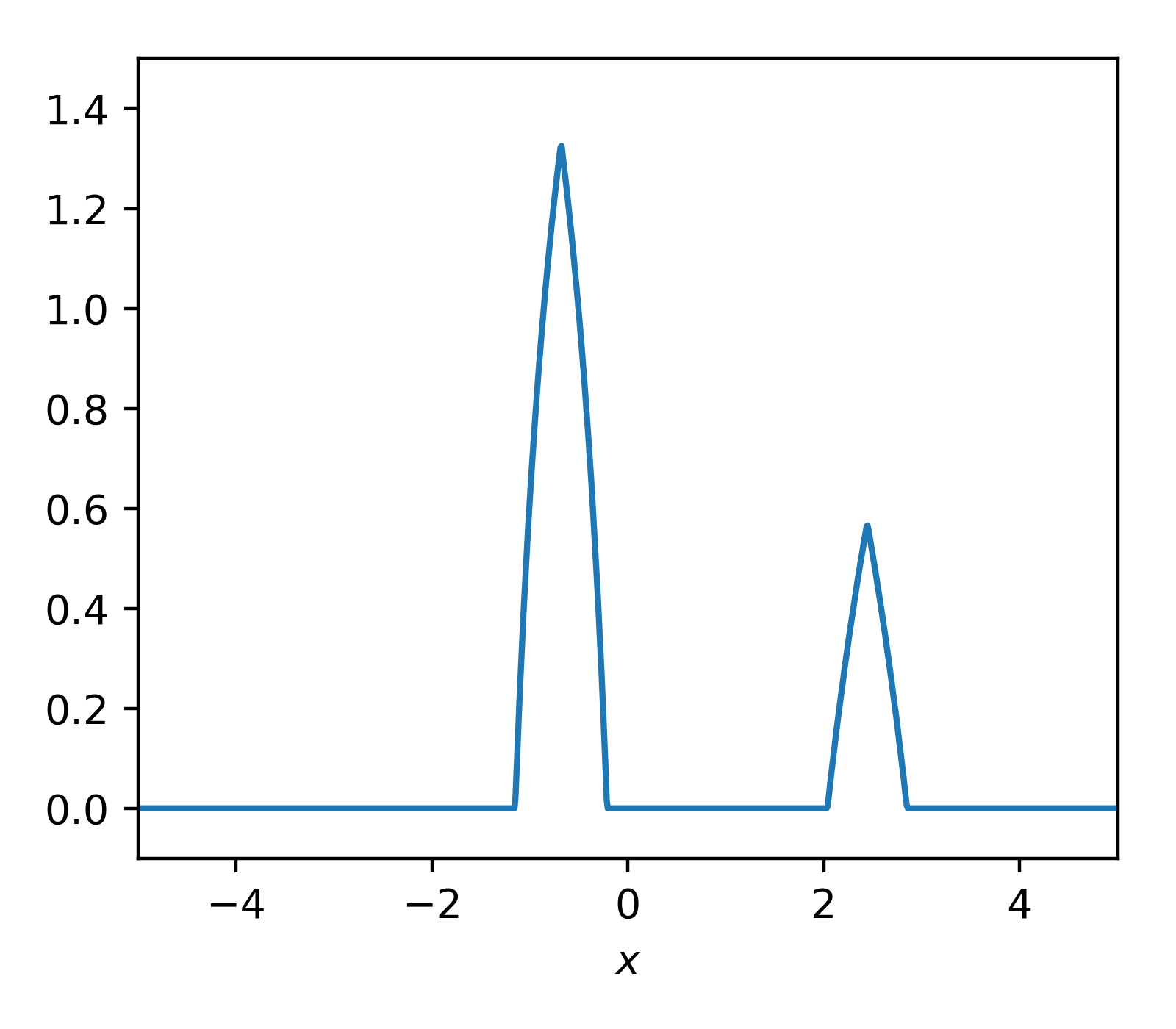}};
\node[inner sep=0pt] (u2) at ( 0.70, 1.30-2.05)
    {\includegraphics[width=0.23\textwidth]{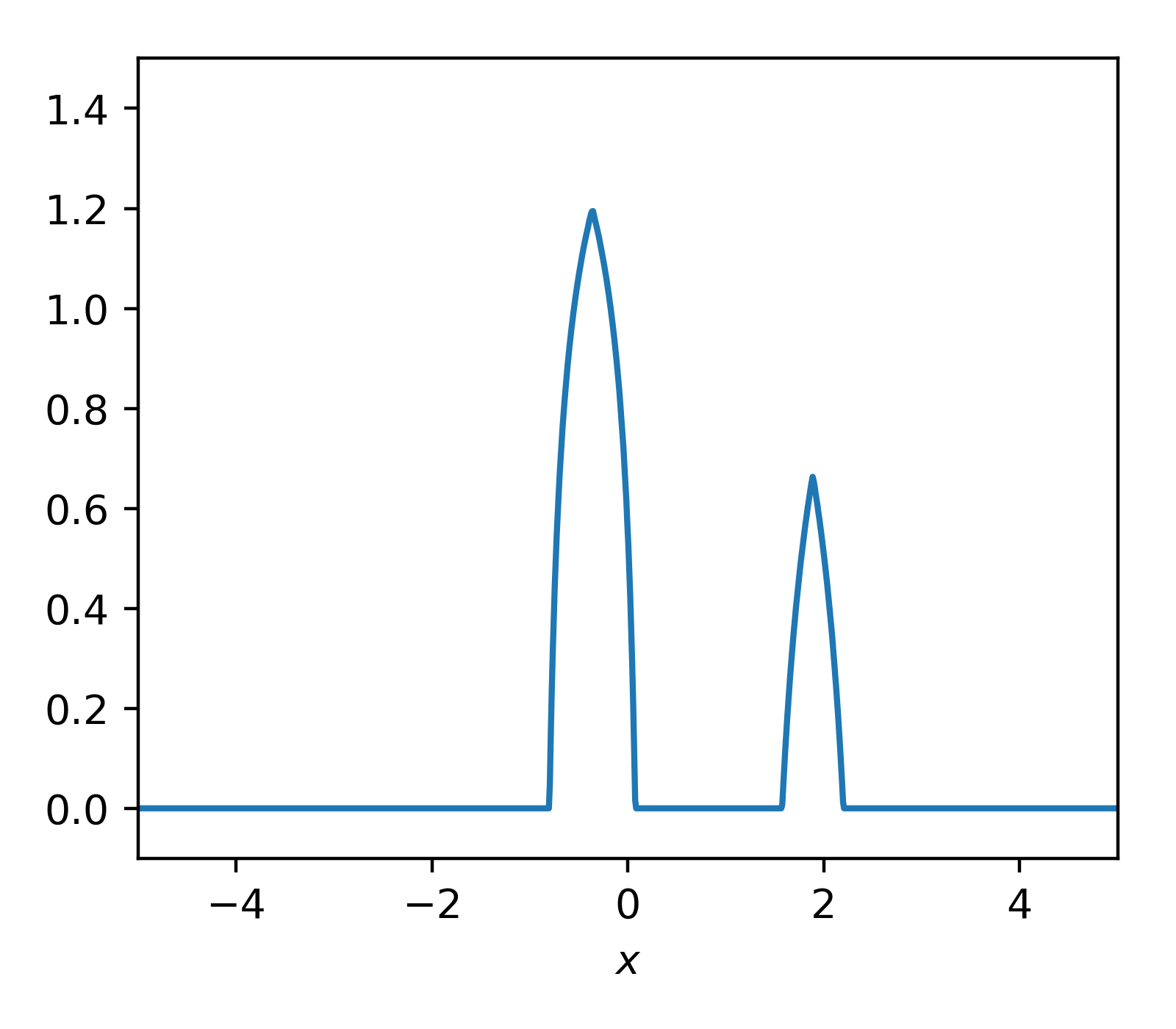}};
\node[inner sep=0pt] (u3) at ( 2.50, 1.30-2.05)
    {\includegraphics[width=0.23\textwidth]{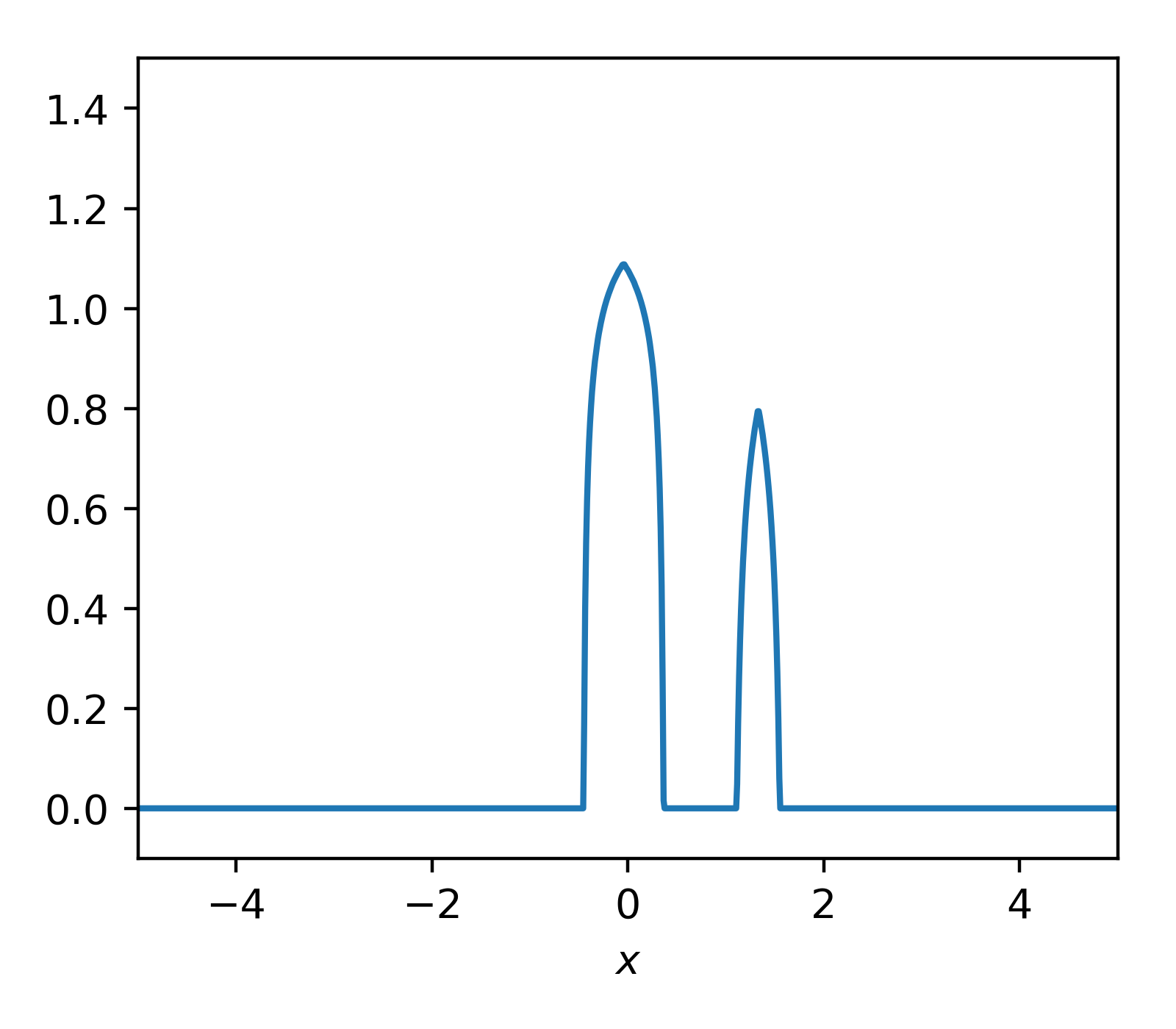}};
\begin{scope}[scale=0.9,shift={(0.20,0.2)}]
\draw[fill=gray!20] (0,0) -- (1,0) -- (0,1) -- (0,0);
\draw[->] (0,0) -- (0  ,1.25);
\draw[->] (0,0) -- (1.25,0  );
\draw[fill=black] (0.  ,0.  ) circle (0.025cm);
\draw[fill=black] (1.  ,0.  ) circle (0.025cm);
\draw[fill=black] (0.  ,1.  ) circle (0.025cm);
\draw[fill=white] (0.25,0.  ) circle (0.025cm) node[anchor=south west]  {$1$};
\draw[fill=white] (0.50,0.  ) circle (0.025cm) node[anchor=south west]  {$2$};
\draw[fill=white] (0.75,0.  ) circle (0.025cm) node[anchor=south west]  {$3$};
\node (b1) at (0.  ,0.  ) {};
\node (b2) at (1.  ,0.  ) {};
\node (a1) at (0.25,0.  ) {};
\node (a2) at (0.50,0.  ) {};
\node (a3) at (0.75,0.  ) {};
\node (a9) at (0.  ,0.25) {};
\node (a8) at (0.  ,0.50) {};
\node (a7) at (0.  ,0.75) {};
\node (a10) at (0.25,0.25) {};
\node (a11) at (0.50,0.25) {};
\node (a12) at (0.25,0.50) {};
\node (a4) at (0.75,0.25) {};
\node (a5) at (0.50,0.50) {};
\node (a6) at (0.25,0.75) {};
\draw (0,0) node[anchor=north]  {$\balpha_1$};
\draw (1,0) node[anchor=north]   {$\balpha_2$};
    \end{scope}
\draw (v1) to[out=   0,in= 180] (b1);
\draw (v2) to[out= 180,in= -40] (b2);
\draw (u1) to[out=  45,in= -90] (a1);
\draw (u2) to[out=  90,in= -90] (a2);
\draw (u3) to[out= 135,in= -90] (a3);
\end{tikzpicture}
    \caption{The displacement interpolants $\tilde{u}_{\balpha}$ for the values 
    $\tilde{\balpha}_1, \tilde{\balpha}_2$, and $\tilde{\balpha}_3$. 
    The specific values are listed in \eqref{eq:talpha}.}
    \label{fig:interp2param123}
\end{figure}

\begin{figure}
    \centering
\begin{tikzpicture}[scale=1.7]
\node[inner sep=0pt] (v1) at ( 2.50, 1.30-1.25)
    {\includegraphics[width=0.23\textwidth]{exp3_dinterp_10_00.png}};
\node[inner sep=0pt] (v2) at (-1.10, 1.30+1.25)
    {\includegraphics[width=0.23\textwidth]{exp3_dinterp_00_10.png}};
\node[inner sep=0pt] (u4) at ( 2.50, 1.30-0.00)
    {\includegraphics[width=0.23\textwidth]{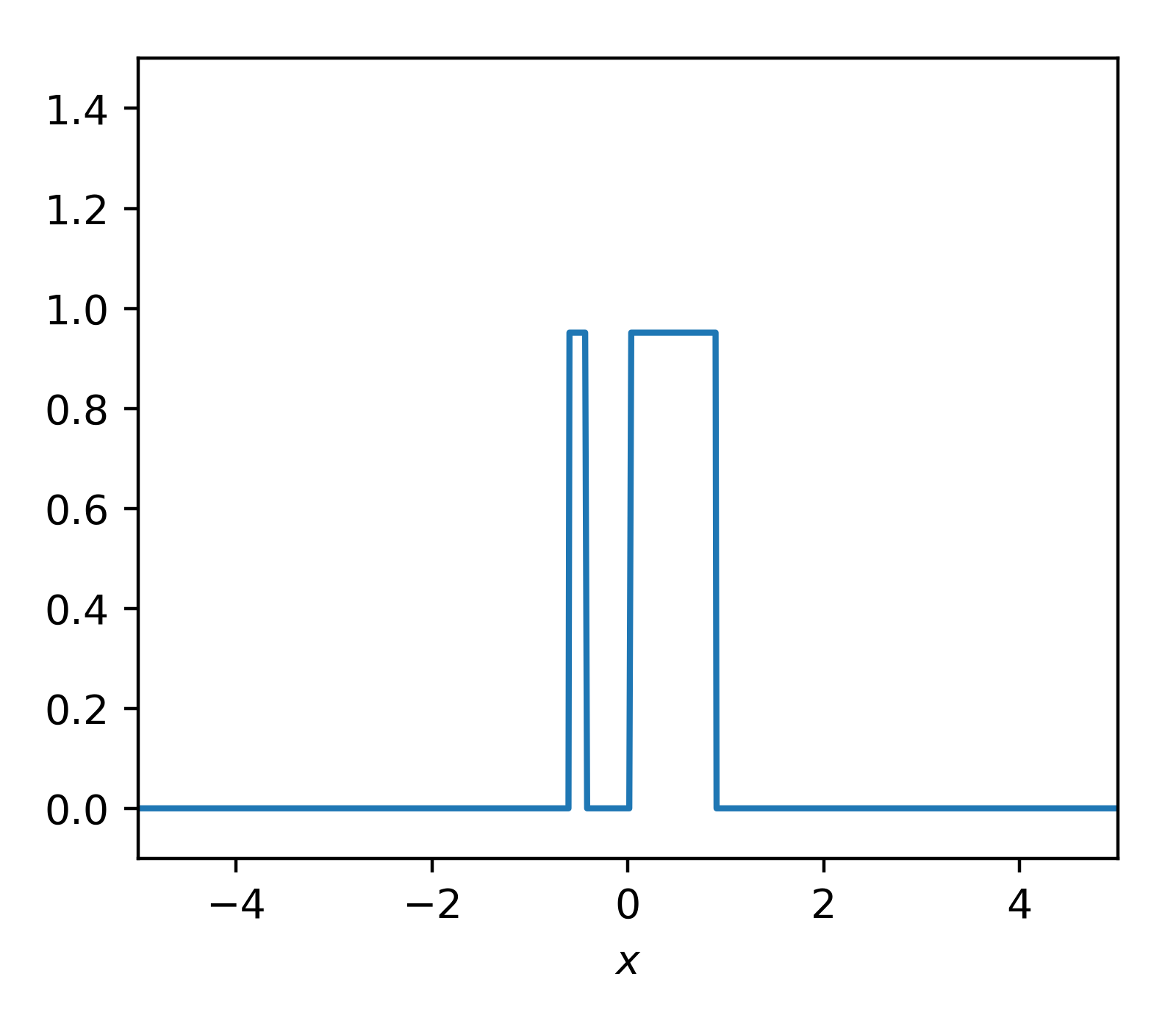}};
\node[inner sep=0pt] (u5) at ( 2.50, 1.30+1.25)
    {\includegraphics[width=0.23\textwidth]{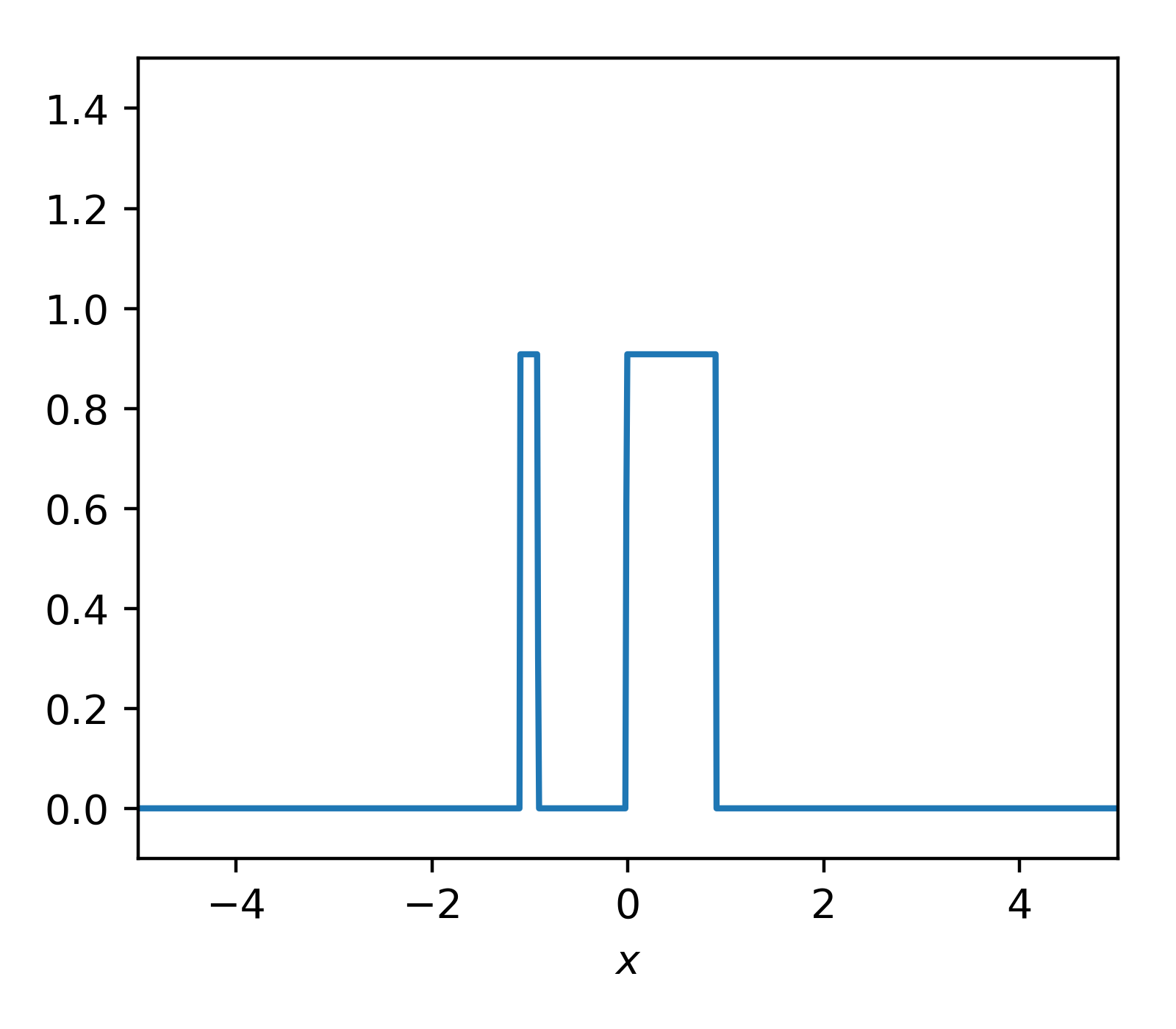}};
\node[inner sep=0pt] (u6) at ( 0.70, 1.30+1.25)
    {\includegraphics[width=0.23\textwidth]{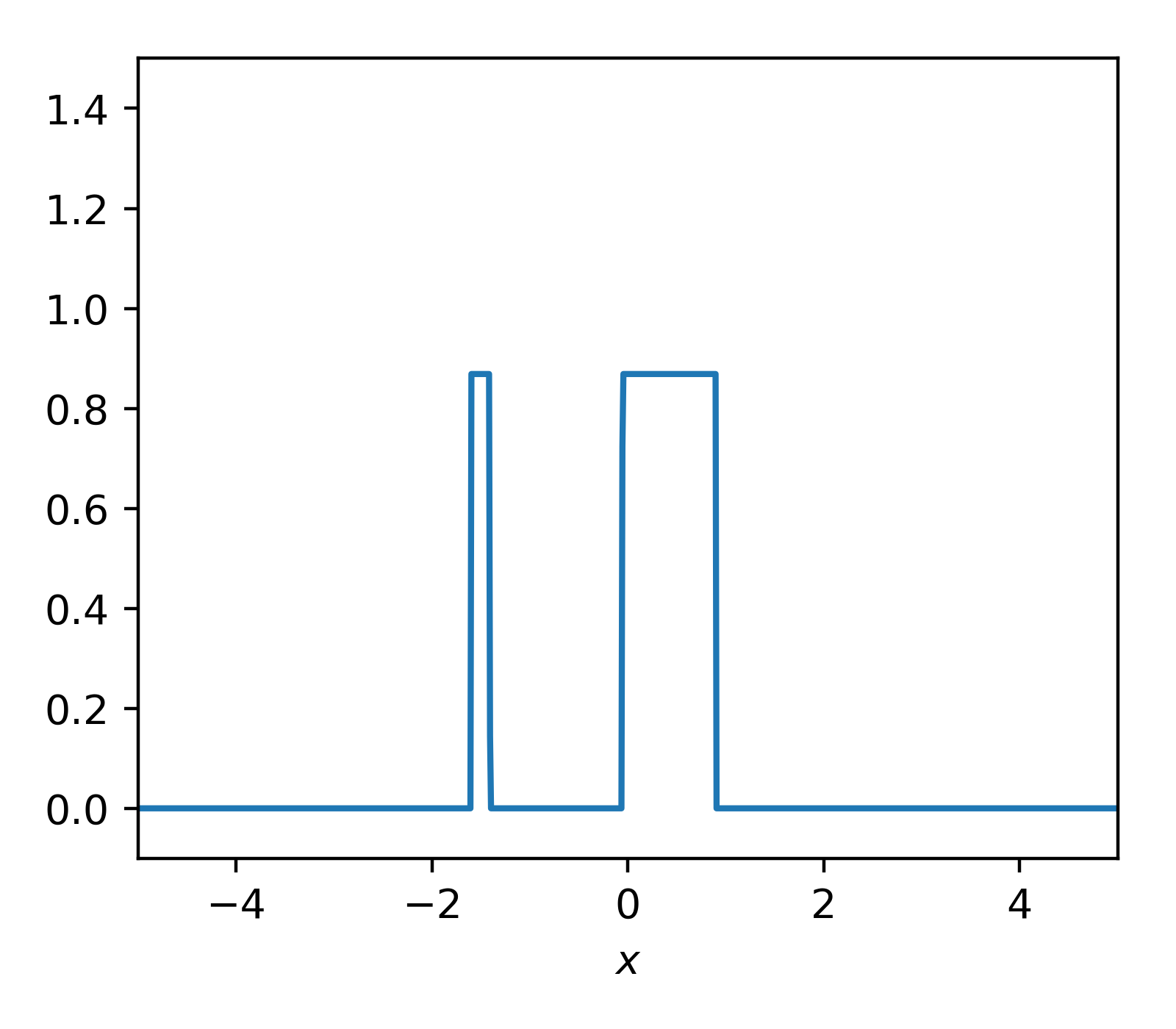}};    
\begin{scope}[scale=0.9,shift={(0.20,0.8)}]
\draw[fill=gray!20] (0,0) -- (1,0) -- (0,1) -- (0,0);
\draw[->] (0,0) -- (0  ,1.25);
\draw[->] (0,0) -- (1.25,0  );
\draw[fill=black] (0.  ,0.  ) circle (0.025cm);
\draw[fill=black] (1.  ,0.  ) circle (0.025cm);
\draw[fill=black] (0.  ,1.  ) circle (0.025cm);
\draw[fill=white] (0.75,0.25) circle (0.025cm) node[anchor=south west]  {$4$};
\draw[fill=white] (0.50,0.50) circle (0.025cm) node[anchor=south west]  {$5$};
\draw[fill=white] (0.25,0.75) circle (0.025cm) node[anchor=south west]  {$6$};
\node (b1) at (1.  ,0.  ) {};
\node (b2) at (0.  ,1.  ) {};
\node (a1) at (0.25,0.  ) {};
\node (a2) at (0.50,0.  ) {};
\node (a3) at (0.75,0.  ) {};
\node (a9) at (0.  ,0.25) {};
\node (a8) at (0.  ,0.50) {};
\node (a7) at (0.  ,0.75) {};
\node (a10) at (0.25,0.25) {};
\node (a11) at (0.50,0.25) {};
\node (a12) at (0.25,0.50) {};
\node (a4) at (0.75,0.25) {};
\node (a5) at (0.50,0.50) {};
\node (a6) at (0.25,0.75) {};
\draw (1,0) node[anchor=north west]   {$\balpha_2$};
\draw (0,1) node[anchor=south west ]   {$\balpha_3$}; 
    \end{scope}
\draw (v1) to[out= 180,in=   0] (b1);
\draw (v2) to[out= -90,in= 180] (b2);
\draw (u4) to[out= 180,in=   0] (a4);
\draw (u5) to[out=-150,in=   0] (a5);
\draw (u6) to[out=- 90,in=   0] (a6);
\end{tikzpicture}
    \caption{The displacement interpolants $\tilde{u}_{\balpha}$ for the values 
    $\tilde{\balpha}_4, \tilde{\balpha}_5$, and $\tilde{\balpha}_6$. 
    The specific values are listed in \eqref{eq:talpha}.}
    \label{fig:interp2param456}
\end{figure}

\begin{figure}
    \centering
\begin{tikzpicture}[scale=1.7]
\node[inner sep=0pt] (v1) at ( 0.8, 1.30+1.75)
    {\includegraphics[width=.23\textwidth]{exp3_dinterp_00_10.png}};
\node[inner sep=0pt] (v2) at ( 0.8, 1.30-0.75)
    {\includegraphics[width=.23\textwidth]{exp3_dinterp_00_00.png}};
\node[inner sep=0pt] (u7) at (-1.1, 1.30+1.75)
    {\includegraphics[width=.23\textwidth]{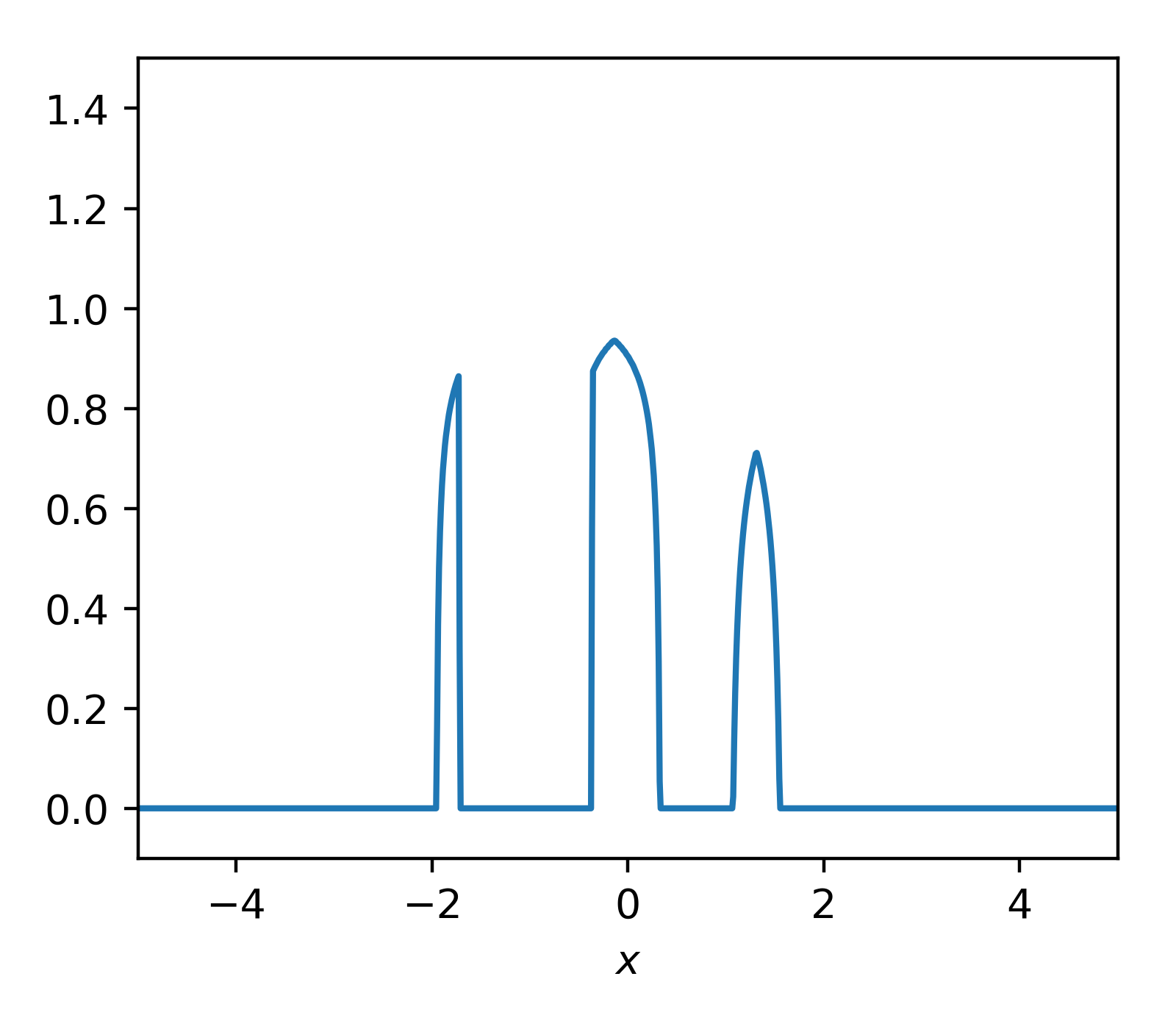}};
\node[inner sep=0pt] (u8) at (-1.1, 1.30+0.50)
    {\includegraphics[width=.23\textwidth]{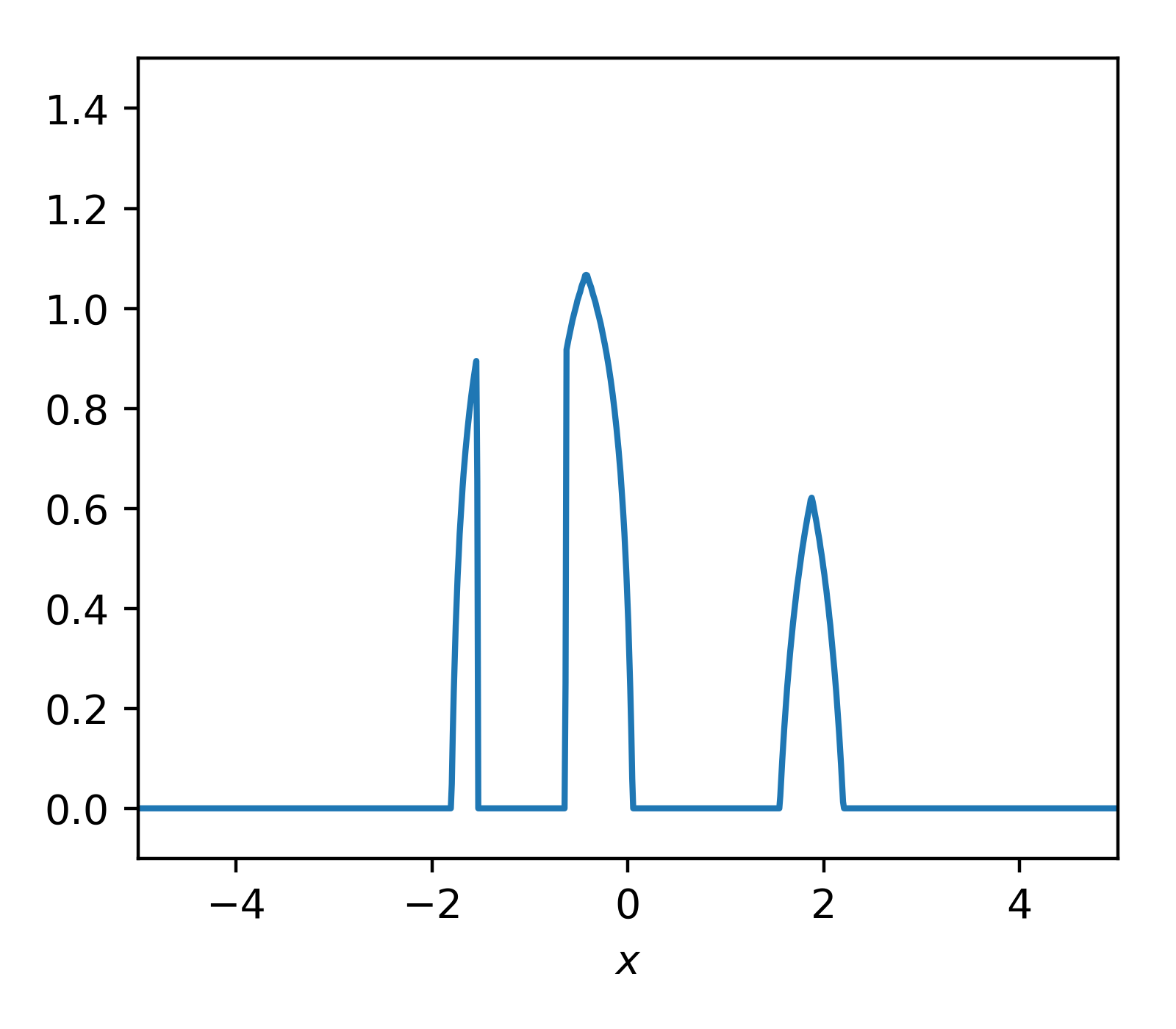}};
\node[inner sep=0pt] (u9) at (-1.1, 1.30-0.75)
    {\includegraphics[width=.23\textwidth]{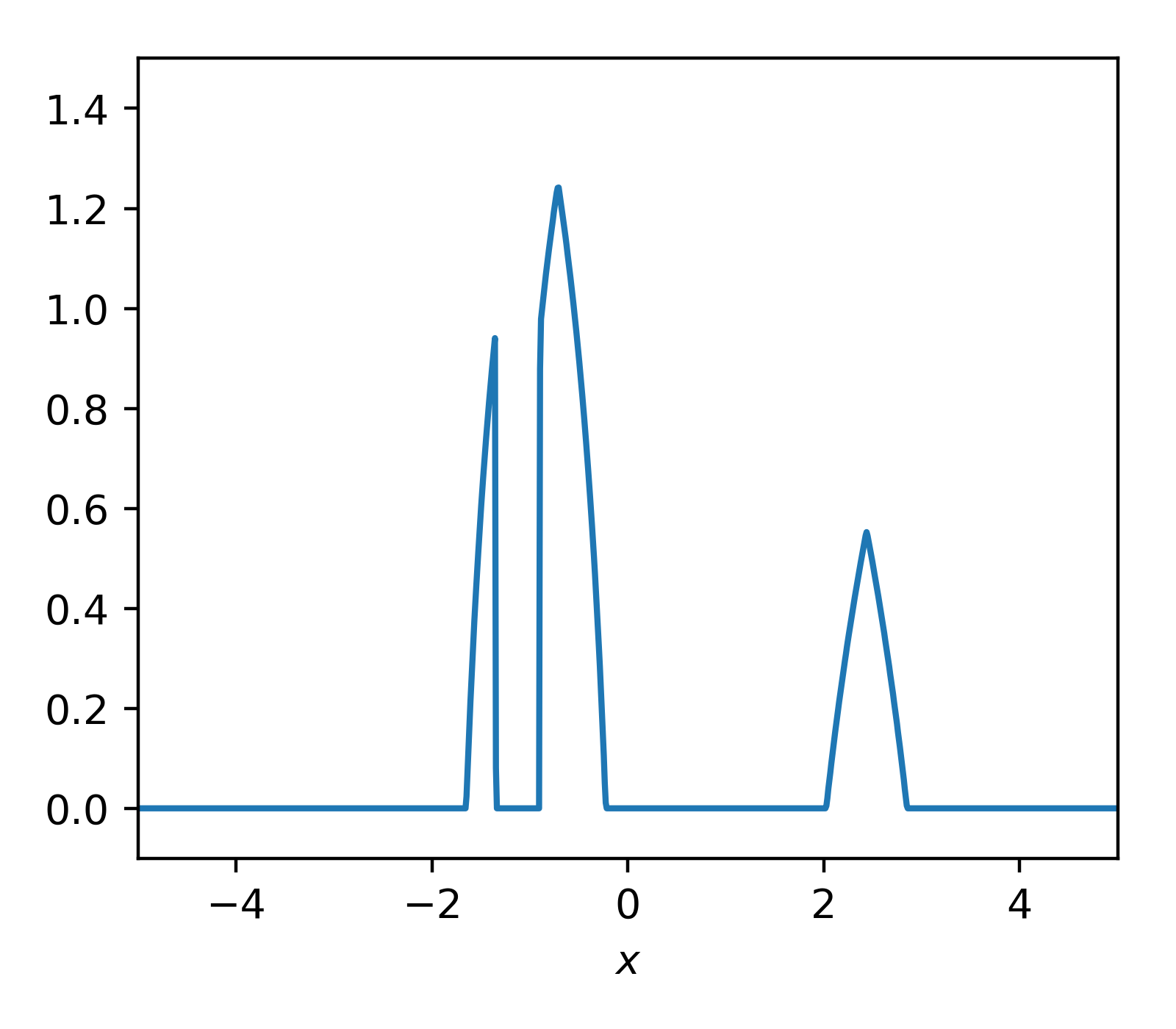}};
\begin{scope}[scale=0.9,shift={(0.20,1.5)}]
\draw[fill=gray!20] (0,0) -- (1,0) -- (0,1) -- (0,0);
\draw[->] (0,0) -- (0  ,1.25);
\draw[->] (0,0) -- (1.25,0  );
\draw[fill=black] (0.  ,0.  ) circle (0.025cm);
\draw[fill=black] (1.  ,0.  ) circle (0.025cm);
\draw[fill=black] (0.  ,1.  ) circle (0.025cm);
\draw[fill=white] (0.  ,0.25) circle (0.025cm) node[anchor=south west]  {$9$};
\draw[fill=white] (0.  ,0.50) circle (0.025cm) node[anchor=south west]  {$8$};
\draw[fill=white] (0.  ,0.75) circle (0.025cm) node[anchor=south west]  {$7$};
\node (b1) at (0.  ,1.  ) {};
\node (b2) at (0.  ,0.  ) {};
\node (a9) at (0.  ,0.25) {};
\node (a8) at (0.  ,0.50) {};
\node (a7) at (0.  ,0.75) {};
\draw (0,0) node[anchor=north east]  {$\balpha_1$};
\draw (0,1) node[anchor=south east ]   {$\balpha_3$}; 
    \end{scope}
\draw (v1) to[out= -90,in=  45] (b1);
\draw (v2) to[out= 130,in= -90] (b2);
\draw (u7) to[out= -30,in= 160] (a7);
\draw (u8) to[out=   0,in= 180] (a8);
\draw (u9) to[out=   0,in= 180] (a9);
\end{tikzpicture}
    \caption{The displacement interpolants $\tilde{u}_{\balpha}$ for the values 
    $\tilde{\balpha}_7, \tilde{\balpha}_8$, and $\tilde{\balpha}_9$. 
    The specific values are listed in \eqref{eq:talpha}.}
    \label{fig:interp2param789}
\end{figure}

\begin{figure}
    \centering
\begin{tikzpicture}[scale=1.9]
\node[inner sep=0pt] (u11) at ( 2.3, 1.30+0.)
    {\includegraphics[width=.23\textwidth]{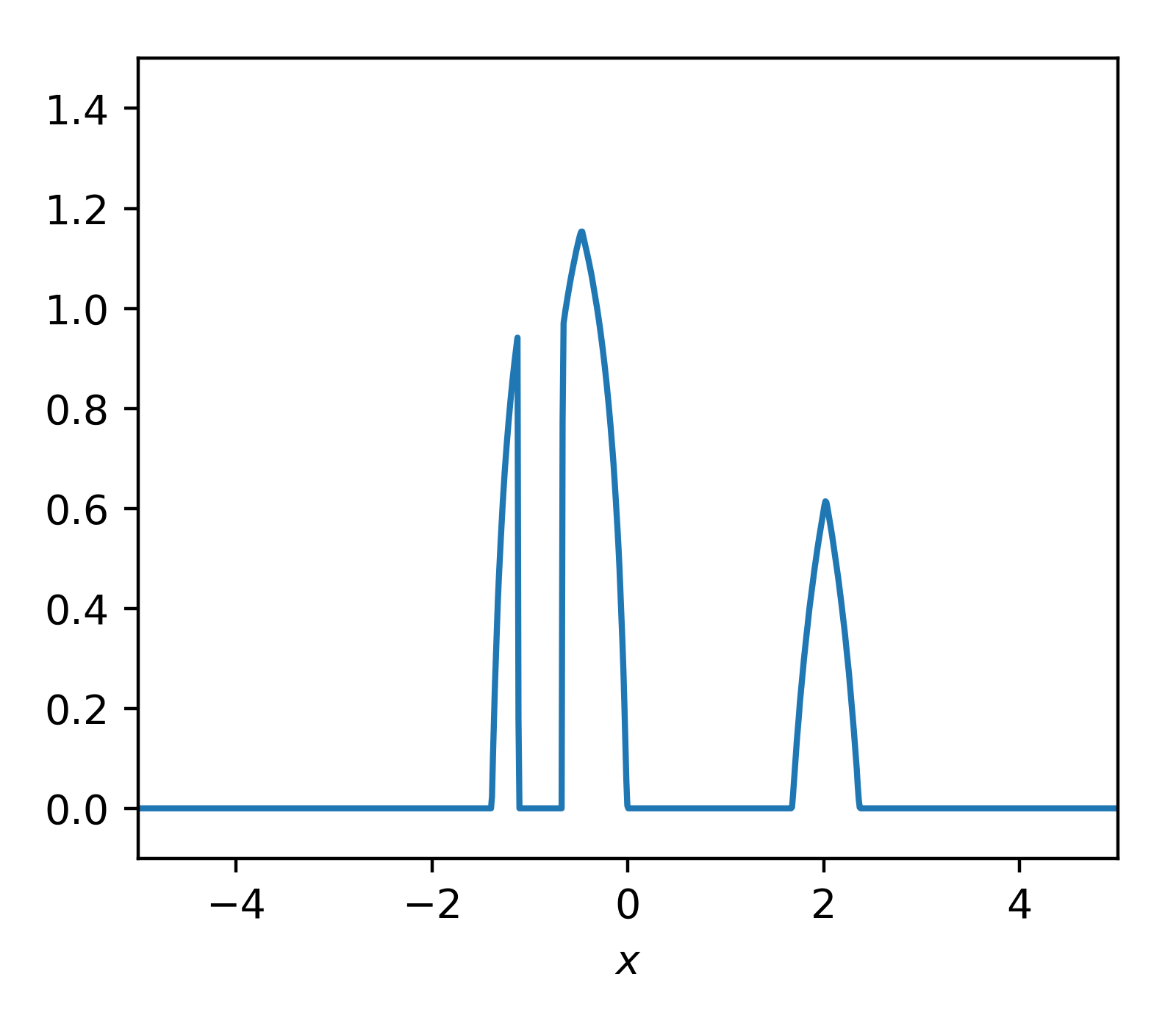}};
\node[inner sep=0pt] (u12) at (-1.1, 1.30+0.80)
    {\includegraphics[width=.23\textwidth]{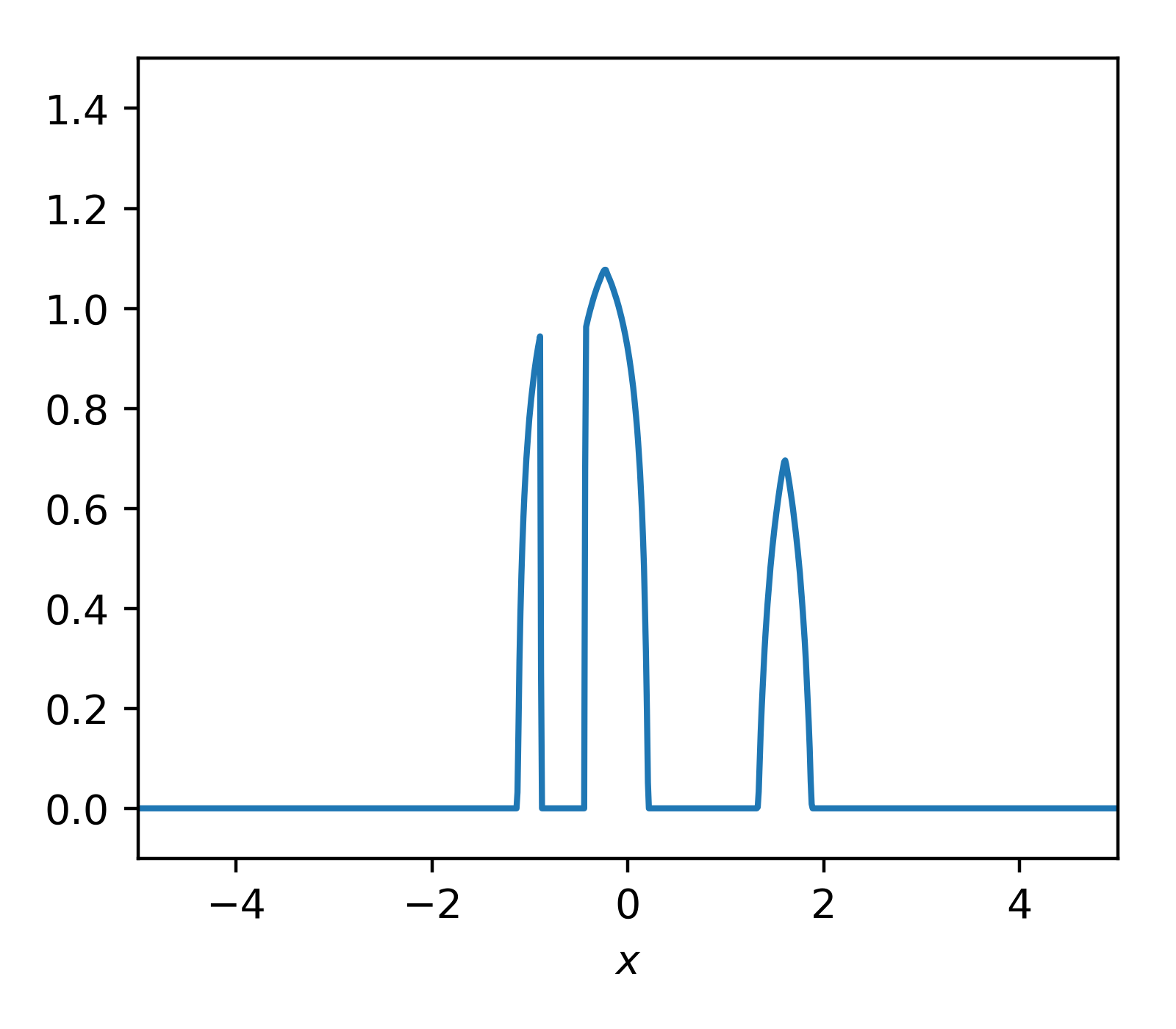}};
\node[inner sep=0pt] (u10) at (-1.1, 1.30-0.55)
    {\includegraphics[width=.23\textwidth]{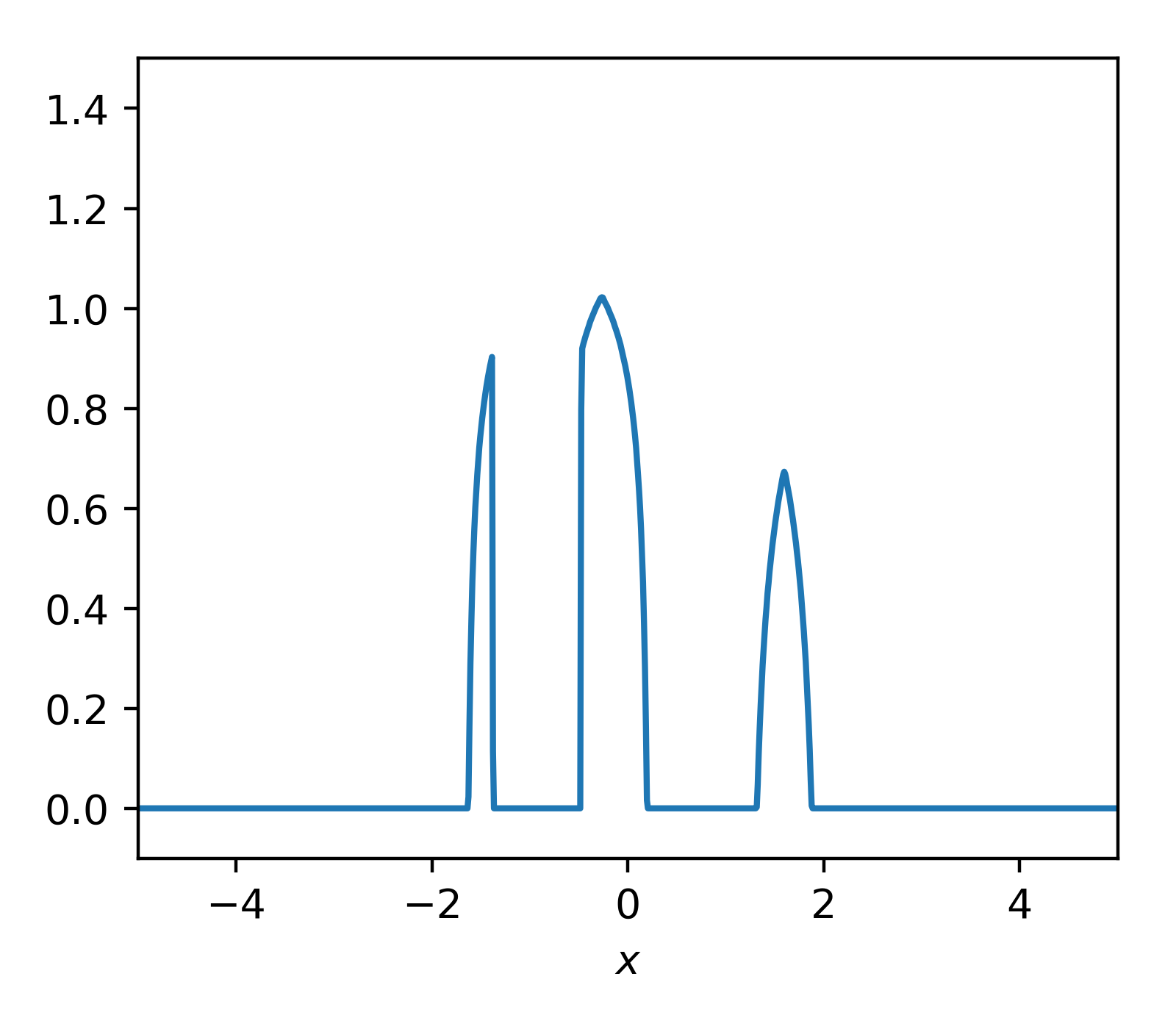}};
\begin{scope}[scale=0.9,shift={(0.20,0.8)}]
\draw[fill=gray!20] (0,0) -- (1,0) -- (0,1) -- (0,0);
\draw[->] (0,0) -- (0  ,1.25);
\draw[->] (0,0) -- (1.25,0  );
\draw[fill=black] (0.  ,0.  ) circle (0.025cm);
\draw[fill=black] (1.  ,0.  ) circle (0.025cm);
\draw[fill=black] (0.  ,1.  ) circle (0.025cm);
\draw[fill=white] (0.25,0.25) circle (0.025cm) node[anchor=south west]  {$10$};
\draw[fill=white] (0.50,0.25) circle (0.025cm) node[anchor=south west]  {$11$};
\draw[fill=white] (0.25,0.50) circle (0.025cm) node[anchor=south west]  {$12$};
\node (a10) at (0.25,0.25) {};
\node (a11) at (0.50,0.25) {};
\node (a12) at (0.25,0.50) {};
\draw (0,0) node[anchor=north east]  {$\balpha_1$};
\draw (1,0) node[anchor=north west]   {$\balpha_2$};
\draw (0,1) node[anchor=south east]   {$\balpha_3$}; 
    \end{scope}
\draw (u10) to[out=  20,in=-150] (a10);
\draw (u11) to[out= 180,in=   0] (a11);
\draw (u12) to[out=   0,in= 180] (a12);
\end{tikzpicture}
    \caption{The displacement interpolants $\tilde{u}_{\balpha}$ for the values 
    $\tilde{\balpha}_{10}, \tilde{\balpha}_{11}$, and $\tilde{\balpha}_{12}$. 
    The specific values are listed in \eqref{eq:talpha}.}
    \label{fig:interp2param101112}
\end{figure}

\cref{fig:interp2param123} displays the displacement interpolant $\tilde{u}$ at
the nodes $\tilde{\balpha}_1, \tilde{\balpha}_2, \tilde{\balpha}_3$, and shows
a gradual deformation of $u_1$ into $u_2$. One notices that the acute angle at
the tip of the two linear hat functions of $u_1$ gradually relax towards a
horizontal line. The two peaks also are transported towards each other, about
to merge into the profile of $u_2$.  In \cref{fig:interp2param456} we see
$\tilde{u}_{\balpha}$ evaluated at the nodes $\tilde{\balpha}_4,
\tilde{\balpha}_5, \tilde{\balpha}_6$.  The one sharp hump in $u_2$ is split
into two sharp humps and is scaled and transported to become $u_3$. The
monotone rearrangement \eqref{eq:mr1d} ensures that the initial split takes
place precisely at the ratio between the mass of two sharp humps in $u_3$.  In
\cref{fig:interp2param789}, $\tilde{u}_{\balpha}$ is evaluated at
$\tilde{\balpha}_7, \tilde{\balpha}_8, \tilde{\balpha}_9$, the two sharp
profiles in $u_3$ are deformed into $u_1$, and since the hump to the left does
not contain enough mass to form the large hat function on the left in $u_1$, a
bulk of mass is taken from the right hump. As parameter approaches that of
$\balpha_1$ the profiles sharpen and form a peak.  The left two peaks will
merge to form one hat function.  In \cref{fig:interp2param101112} the
displacement interpolant at $\tilde{\balpha}_{10}, \tilde{\balpha}_{11},
\tilde{\balpha}_{12}$ are shown. These corresponds to interior points, and we
observe interesting behavior of three hump functions of various shapes.

The behavior is predictable and natural, and the location of the profiles as
well as their shapes are adjusted smoothly according to the paramters.  This is
in contrast to linear interpolation, where the basis functions are largely
stationary.

\begin{figure}
    \centering
    \begin{tabular}{cc}
        \includegraphics[width=0.35\textwidth]{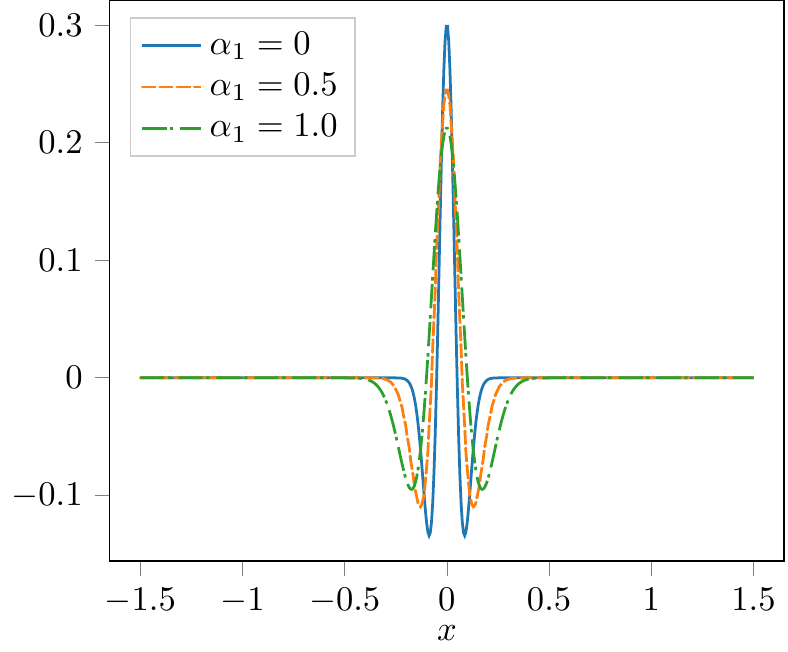}
        &
        \includegraphics[width=0.35\textwidth]{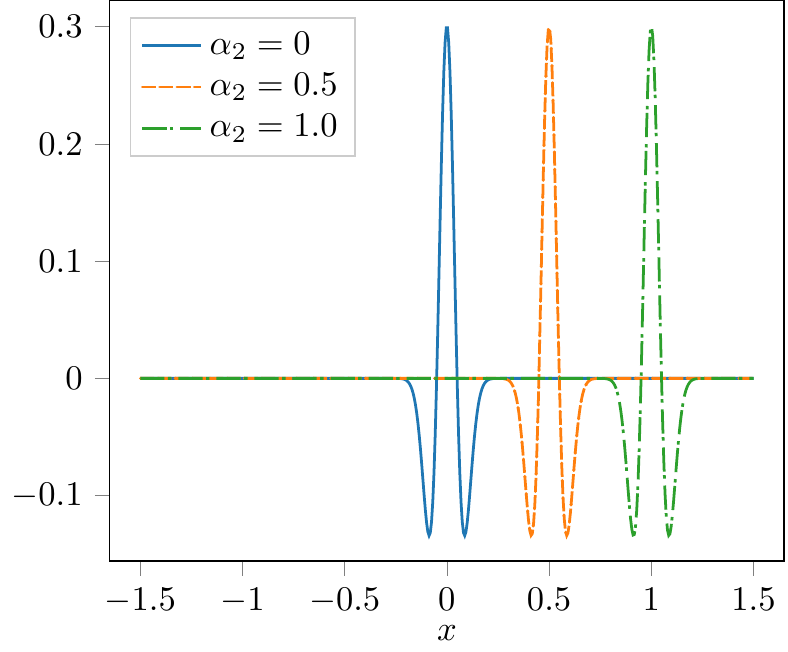}\\
    \end{tabular}
    \caption{\added{Plot of the Mexican hat wavelet function $\psi(x)$ and its 
       displacement interpolants $\zeta(x;\alpha_1,\alpha_2)$ defined in 
        \cref{eq:wavelet_dinterp} for parameter weights 
       $\alpha_1 = 0,0.5,1$ resulting in scaling (left),
       and $\alpha_2 = 0, 0.5, 1$ resulting in translation (right).
      }}
    \label{fig:wavelet_plot}
\end{figure}

\added{
\subsection{Wavelets} \label{sec:wavelets}
It is straightforward to see the continuous wavelet basis
\cite{daubechies} as a special case of a two-parameter version of
\cref{eq:bary_dinterp} (or \cref{eq:2dquantile}), with $u_1(x)$ taken as a
\emph{wavelet function} $\psi(x)$. We perform the following displacement
interpolation,
\beq
\zeta(x;\alpha_1,\alpha_2) := \cI(\alpha_1,\alpha_2;u_1,u_2,u_3), 
\quad \text{ where } 
\left\{
\begin{aligned}
u_1(x) &= \psi(x),\\
u_2(x) &= \psi(x/2)/\sqrt{2},\\
u_3(x) &= \psi(x-1).\\
\end{aligned}
\right. 
\label{eq:wavelet_dinterp}
\eeq
The interpolants are shown in \cref{fig:wavelet_plot}.  The continuous wavelet
basis $\{\psi_{a,b}(x)\}$ is constructed by scaling and translating 
$\psi(x)$, which is related to the above interpolation $\cI(\alpha_1,\alpha_2)$
with 
\beq
    \psi_{a,b}(x) = |a|^{-\frac{1}{2}}\psi\left( \frac{x - b}{a} \right) 
                  =  C(a,b) \zeta \left(x; \, a, \frac{b}{a}\right).
\eeq
Therefore the wavelet basis results from $\zeta$ in \cref{eq:wavelet_dinterp}
up to a scalar multiple $C$.

In this sense, the displacement interpolation produces a generalization of the
wavelet basis, beyond simple operations like dilation or translation. This
observation can be expressed as an expansion in the characteristic variables,
see \cite{iollo14,rim18}.
}

\added{
\subsection{Multiple hat functions}\label{sec:hats} 
Here we will revisit the example of multiple hat functions displayed in
\cref{fig:hats}. Recall the definition of hat functions $\phi(x;w)$
\cref{eq:hatf}. Let us define a collection of these hat functions with two
parameters, the translate $t$ and the width $w$,
\beq
\cU = \{\phi(x-t,w): 0.1 < t < 0.9, 0 < w < 0.1\}.
\label{eq:randomhats}
\eeq
We will draw functions from $\cU$ at random, by selecting random values of $t$
and $w$ with the uniform distribution over their allowed intervals to form
$\cU_n = \{u_n\}$.  Then we compute the SVD of the snapshot matrix $A$ whose
$(n-1)$-th column is the difference of normalized pseudo-inverses $U_n^\dagger -
U_1^\dagger$
\eqref{eq:bary_dinterp}, 
\beq
A = \begin{bmatrix} \cdots, & 
  \frac{U_n^\dagger - U_1^\dagger}
       {\Norm{U_n^\dagger - U_1^\dagger}{2}}, & \cdots \end{bmatrix}.
\label{eq:randomhats_snapshot}
\eeq
The matrix $A$ is just a normalization of the gradient of
$\widetilde{U}_{\balpha}^\dagger$ with respect to $\alpha_j$ in
\cref{eq:bary_dinterp}.  We see that there are two singular values that are the
most significant as shown in \cref{fig:hats_svd}.  The singular vectors
corresponding to these singular values coincide roughly with translation and
dilation, as also shown in the same figure. One may compare these singular
vectors with the case when only translation is present, as displayed in
\cref{fig:hats_CDFs}.}

\added{
\begin{figure}
    \centering
    \includegraphics[width=0.80\textwidth]{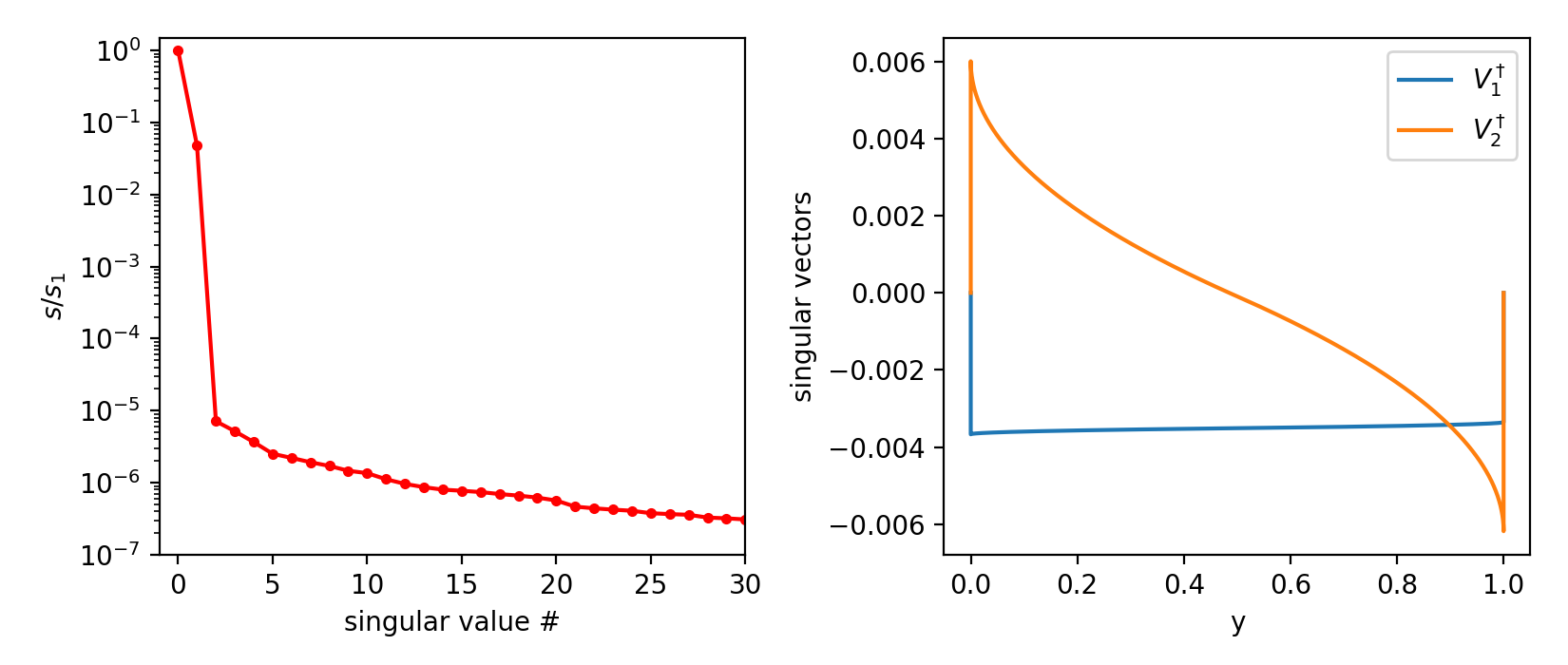}
    \caption{\added{Semi-log plot of scaled singular values $s/s_1$ of the
snapshot matrix \cref{eq:randomhats_snapshot} (left) and the two leading
singular vectors (right). $V_1^\dagger$ and $V_2^\dagger$ correspond roughly to
translation and dilation; see \cref{fig:hats_CDFs} for comparison.}} 
    \label{fig:hats_svd}.  
\end{figure}
}

\subsection{1D acoustics equations}\label{sec:1dacoustics}

\begin{figure}
    \centering
    \begin{tabular}{ll}
    (a) & (b) \\
    \includegraphics[width=0.35\textwidth]{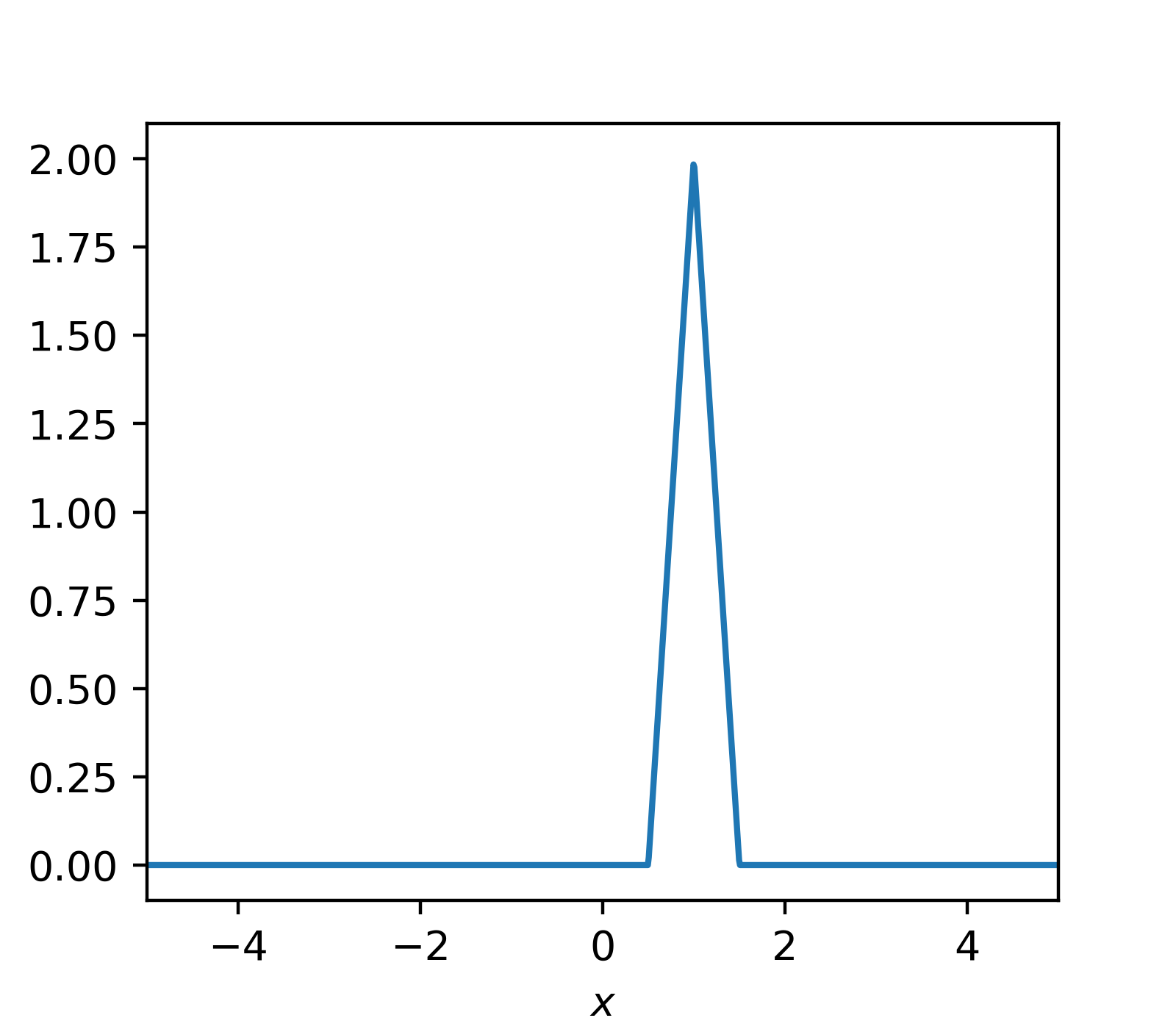}
        &
    \includegraphics[width=0.35\textwidth]{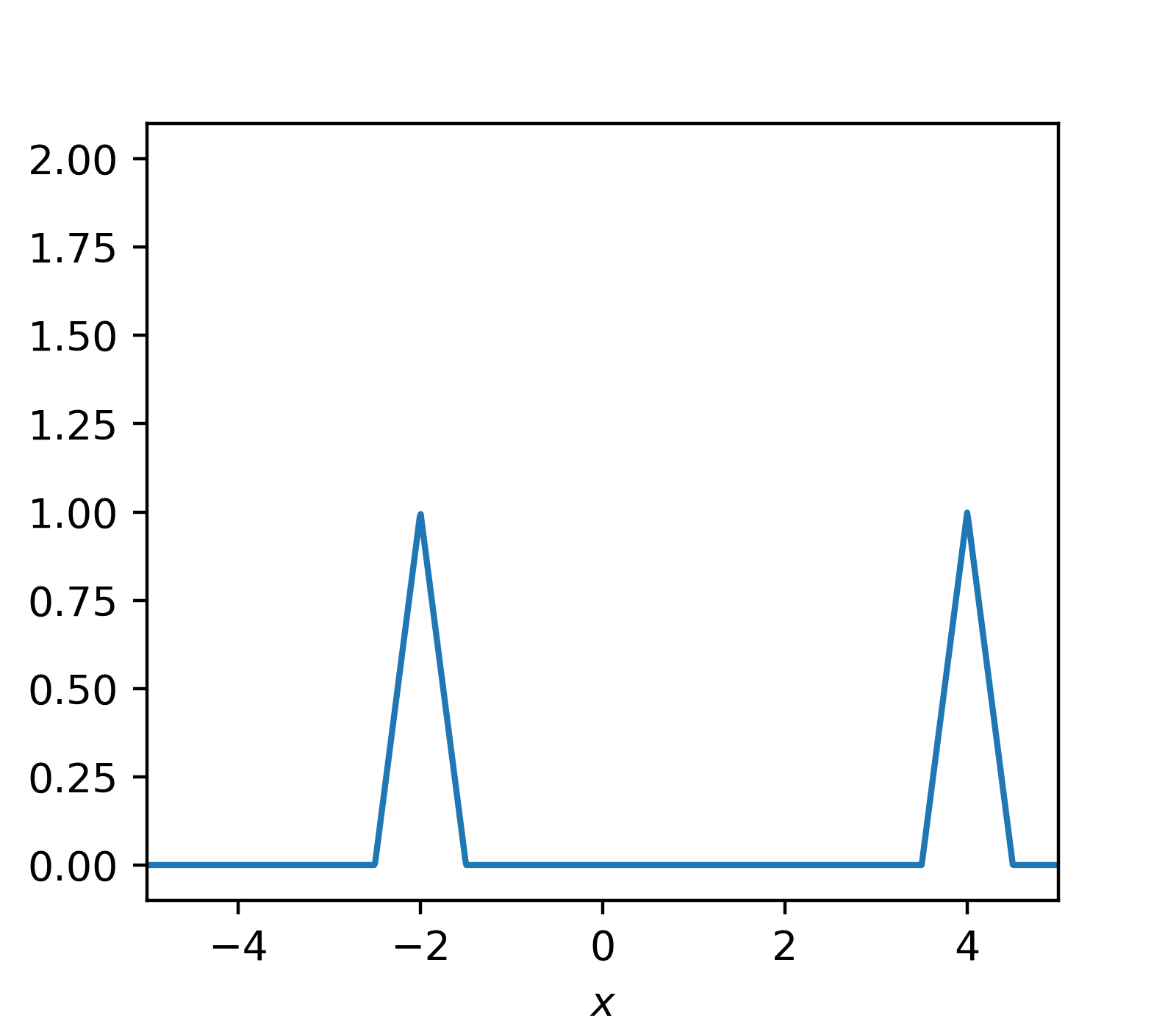}
        \\
    (c) & (d) \\
    \includegraphics[width=0.35\textwidth]{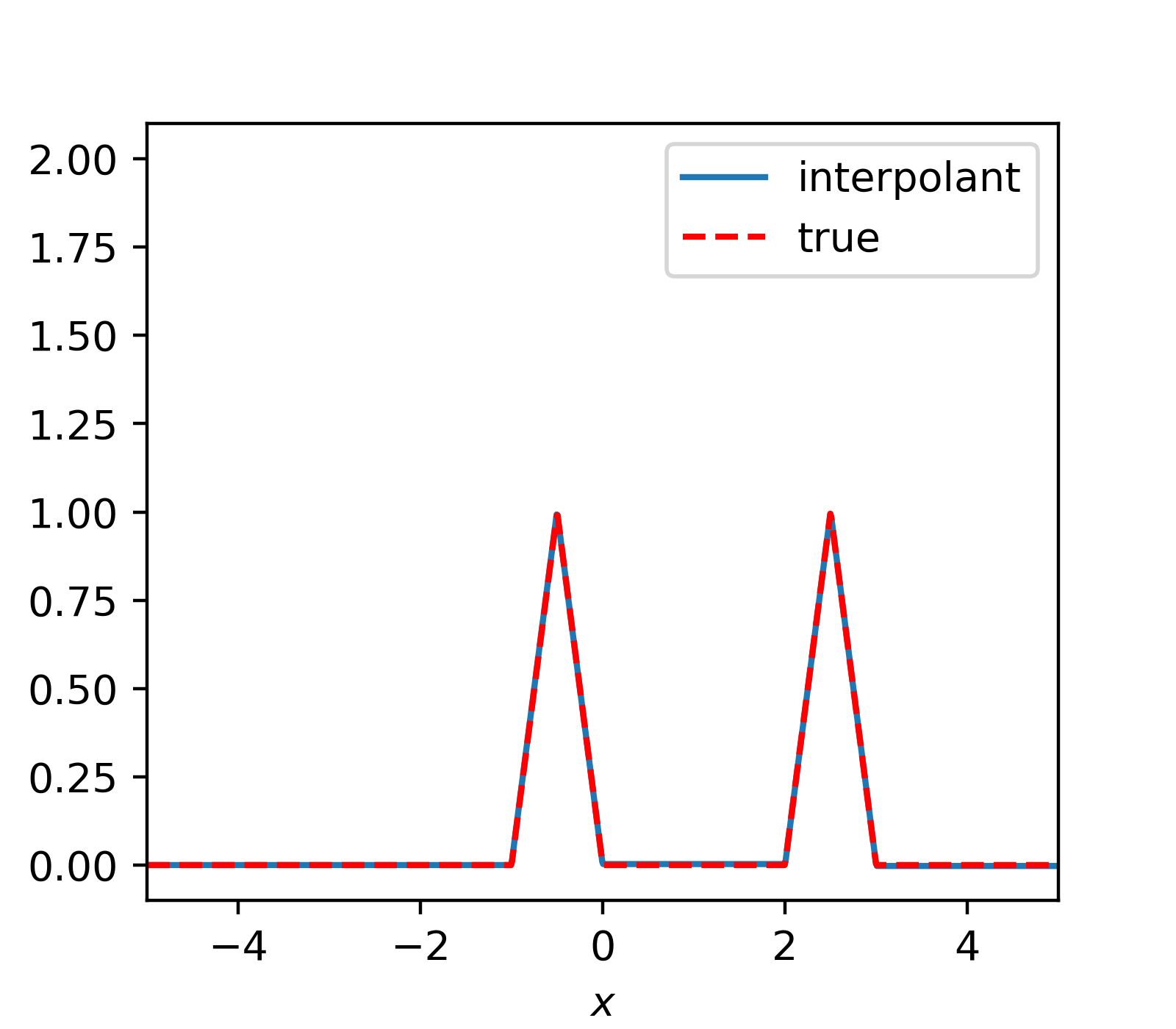}
        &
    \includegraphics[width=0.35\textwidth]{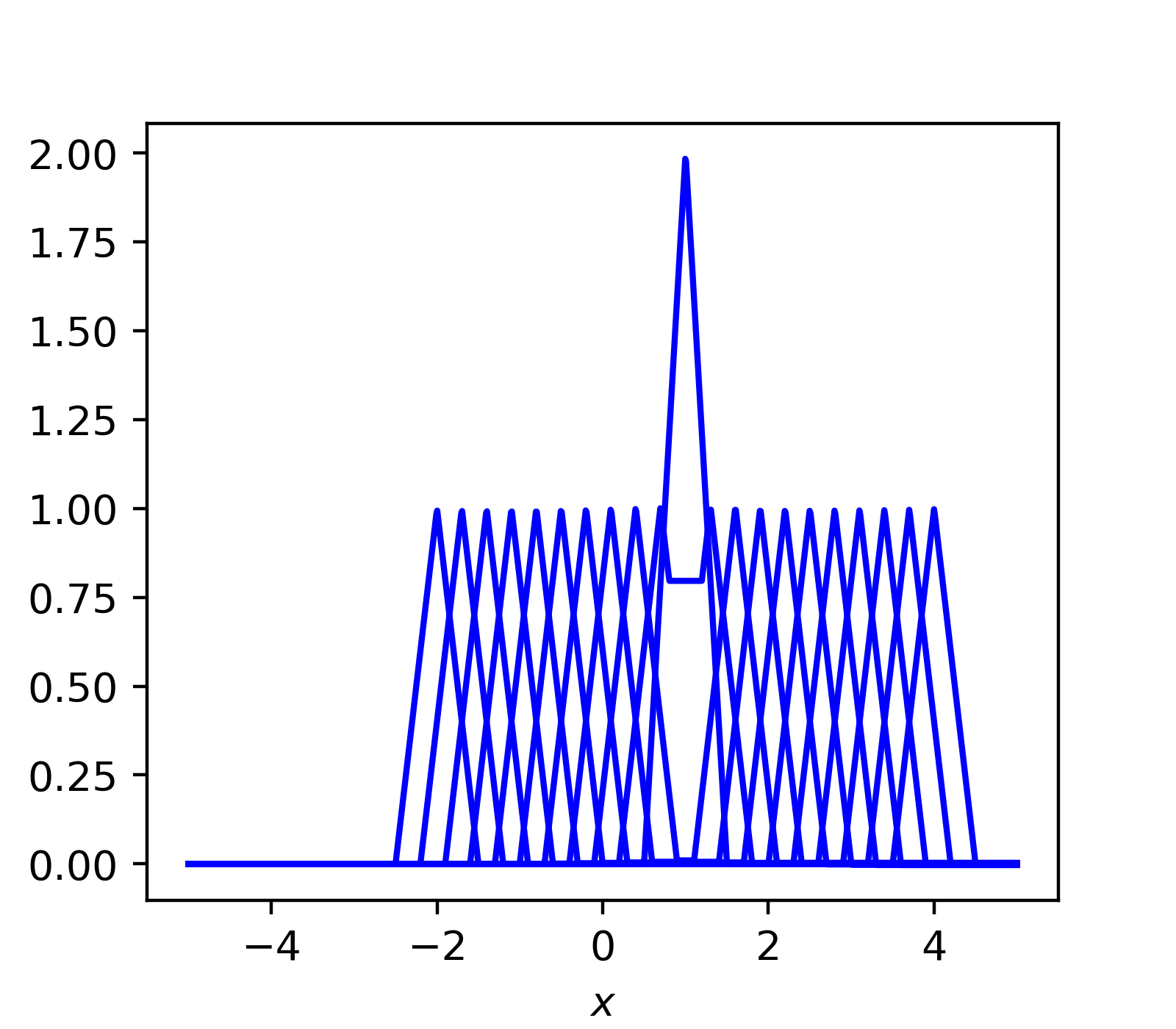}
    \end{tabular}
    \caption{(a) Solution for the pressure variable $p$ to the 
    1D acoustics equations\eqref{eq:acoustic1d} at times $t_1 = 0$  
    and (b) $t_2 = 3$. (c) Displacement interpolant 
    $\tilde{p}_t$ \eqref{eq:pressureI} to these two
    at $t = 1.5$ or $\lambda(t) = 0.5$, and (d) $\tilde{p}_t$
    at times $t = 0.2n$ for $n = 0,1,2, ... , 10$.}
    \label{fig:acoustic}.
\end{figure}

Consider the 1D acoustics equations for pressure $p$ and velocity $u$,
\begin{equation}
    \begin{bmatrix} p \\ u \end{bmatrix}_t
  + \begin{bmatrix} 0 & K_0 \\ 1/\rho_0 & 0 \end{bmatrix}
    \begin{bmatrix} p \\ u \end{bmatrix}_x = 0,
      \label{eq:acoustic1d}
\end{equation}
where $K_0$ is the bulk modulus and $\rho_0$ is the density.  Let us suppose
the initial condition for the pressure $p(x,0)$ is given by the hat function
\eqref{eq:hatf} that is translated $\phi(x - 2w)$ with $w=0.05$, and $u(x,0) =
0$. Then the pressure at time $t$ is given by
\begin{equation}
    p(x,t) = \frac{1}{2} [\phi(x-2h-ct;h) + \phi(x-2h+ct;h)] 
    \quad \text{ where } c = \sqrt{K_0/\rho_0}.
    \label{eq:dalembert}
\end{equation}
For simplicity, let us assume that $K_0=\rho_0 = 1$, so that $c = 1$.  Suppose
we are interested in the behavior of the pressure $p$ with respect to the
paramter $\alpha=t$, the time variable. We are  given the solution at time $t_1
= 0$ and $t_2 = 3$, as shown in the top row of \cref{fig:acoustic}, and we will
use the displacement interpolant to approximate the solution times between
$t_1$ or $t_2$. 

Let us make an observation regarding the example in \cref{sec:2param}.  When
the number of connected supports of the two given functions are different, the
interpolation allocates the probability mass correspondingly and then divides
the support.  An example of this is shown in \cref{fig:interp2param123}.
However, this behavior is different from that of the d'Alembert solution given
in \eqref{eq:dalembert}, where two copies of smaller amplitudes are produced
then are transported in opposite directions.  In what follows we will show
that, for the given initial condition above, the wave-like behavior of the
d'Alembert solution can be reproduced by applying the interpolation procedure
to the derivatives of the given functions instead. 

First, let us define $p_n(x) = p(x,t_n)$ for $n=1,2$ then compute
\begin{equation}
\begin{aligned}
p_n^{d+d+} &:= %
        \left[ \frac{d}{dx} \left( \frac{d p_n}{d x} \right)^+\right]^+,\\
p_n^{d+d-} &:= %
        \left[ \frac{d}{dx} \left( \frac{d p_n}{d x} \right)^+\right]^-,\\
\end{aligned}
\qquad
\begin{aligned}
p_n^{d-d+} &:= %
        \left[ \frac{d}{dx} \left( \frac{d p_n}{d x} \right)^-\right]^+,\\
p_n^{d-d-} &:= %
        \left[ \frac{d}{dx} \left( \frac{d p_n}{d x} \right)^-\right]^-.\\
\end{aligned}
\end{equation}
In the continuous setting, these functions will be linear combinations
of delta functions, but in a discretized setting these will be piece-wise
constant functions whose support is contained in one grid cell. 

\begin{figure}
    \centering
    \begin{tabular}{c}
    \includegraphics[width=0.80\textwidth]{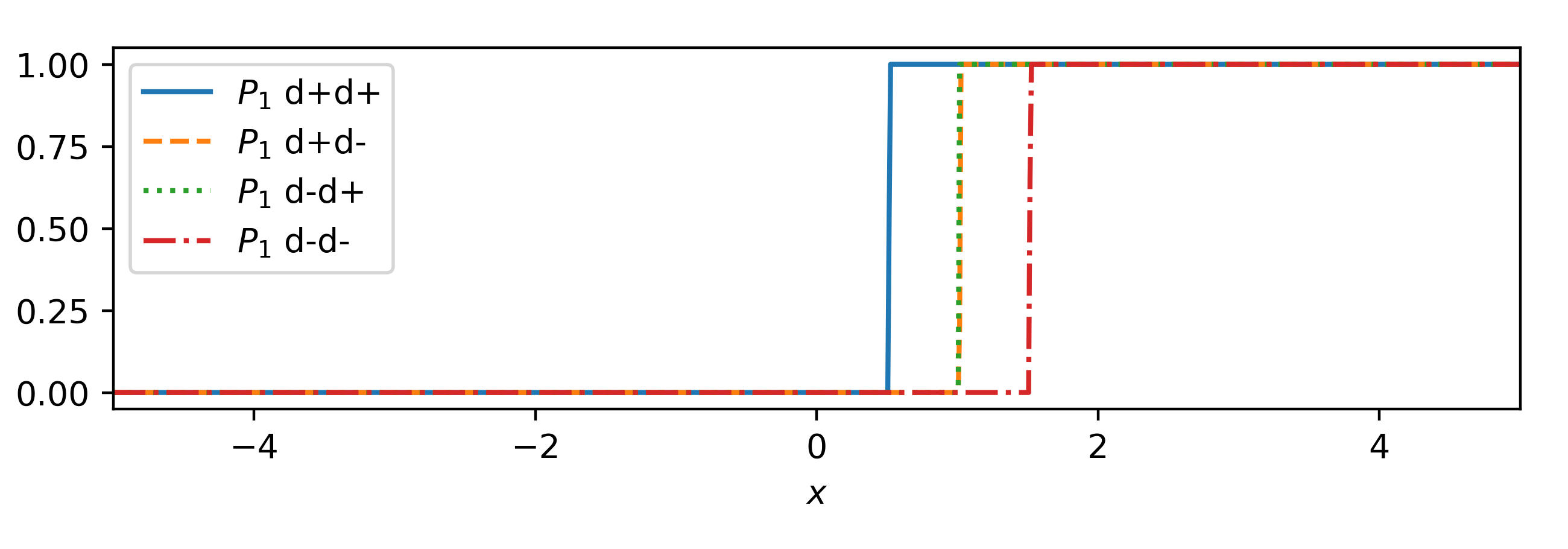} \\
    \includegraphics[width=0.80\textwidth]{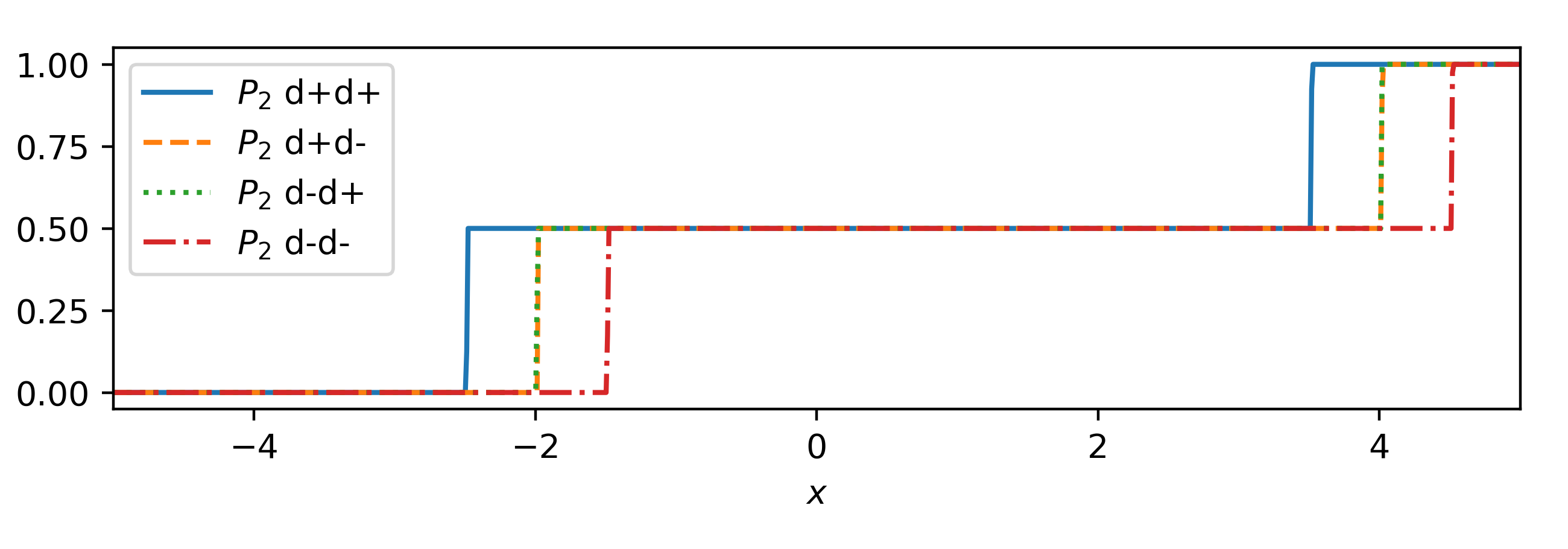} 
    \end{tabular}
    \caption{The CDFs  \eqref{eq:acousticCDF} 
    of the functions $p_1^{d \pm d\pm}$ (top) and 
    $p_2^{d \pm d\pm}$ (bottom). }
    \label{fig:acousticCDF}.
\end{figure}

Then we compute the interpolants
\begin{equation}
    \begin{aligned}
    \tilde{p}^{d+d+}_t &:= \cI(t;p_1^{d+d+},p_2^{d+d+}),\\ 
    \tilde{p}^{d+d-}_t &:= \cI(t;p_1^{d+d-},p_2^{d+d-}),
    \end{aligned}
    \quad
    \begin{aligned}
    \tilde{p}^{d-d+}_t &:= \cI(t;p_1^{d-d+},p_2^{d-d+}),\\ 
    \tilde{p}^{d-d-}_t &:= \cI(t;p_1^{d-d-},p_2^{d-d-}).
    \end{aligned}
    \label{eq:pressureI}
\end{equation}
To compute the interpolant we first obtain the corresponding CDFs 
\eqref{eq:CDFs} for $p_n^{d\pm d\pm}$ for $n=1,2$,
\begin{equation}
    P_n^{d\pm d\pm}(x) := \int_{-\infty}^x p_n^{d\pm d\pm}(y) \, dy,
\label{eq:acousticCDF}
\end{equation}
which are plotted in \cref{fig:acousticCDF}. 

The interpolant is given by summing and integrating,
\begin{equation}
\tilde{p}_t := 
   \int_{-\infty}^x \int_{-\infty}^y   
    \left(\tilde{p}^{d+d+}_t(z) - \tilde{p}^{d+d-}_t(z) 
        + \tilde{p}^{d-d+}_t(z) - \tilde{p}^{d-d-}_t(z)\right) \, dz \, dy.
\end{equation}
The interpolant $\tilde{p}_t$ at $t=1.5$ is shown in in \cref{fig:acoustic}(c).
It agrees exactly with the true pressure $p(x,1.5)$, successfully translating
the two hat functions with two opposite speeds. It is not difficult to verify
that $\tilde{p}_t$ is equal to the exact solution for all values $0< t <3$.
The plot of $\tilde{p}_t$ for 11 different values $t_n =  0.2n$ for $n = 0,1,2,
... , 10$ is shown in the lower right plot in \cref{fig:acoustic}(d).

This demonstrates that for functions whose connected supports are being
translated together with some speed, monotone rearrangement can be used to
exploit the low-dimensional structure by utilizing the linear subspace formed
by the pseudo-inverses of the CDFs computed in \eqref{eq:acousticCDF}.  This
also implies that applying the monotone rearrangement to derivatives of a
function results in a different transport map, and this property could be used
to construct another monotone rearrangement map. Depending on the context,
other maps could be more useful than the solution to the Monge-Kantorovich
problem itself. 

\subsection{1D Burgers shock} \label{sec:burgers}
\begin{figure}
\centering
\begin{tabular}{ll}
    (a) & (b)  \\
    \includegraphics[width=0.35\textwidth]{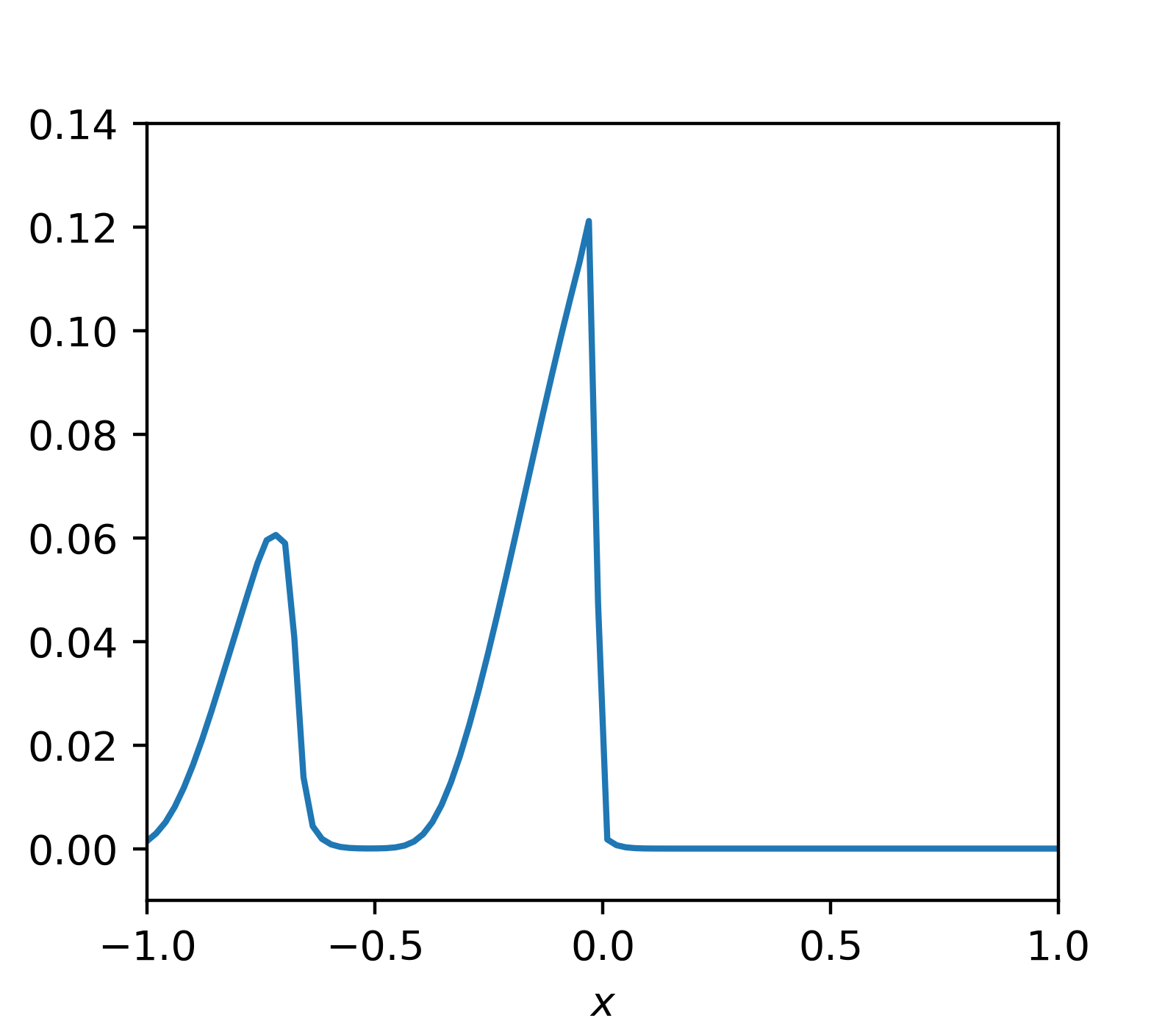}
        &
    \includegraphics[width=0.35\textwidth]{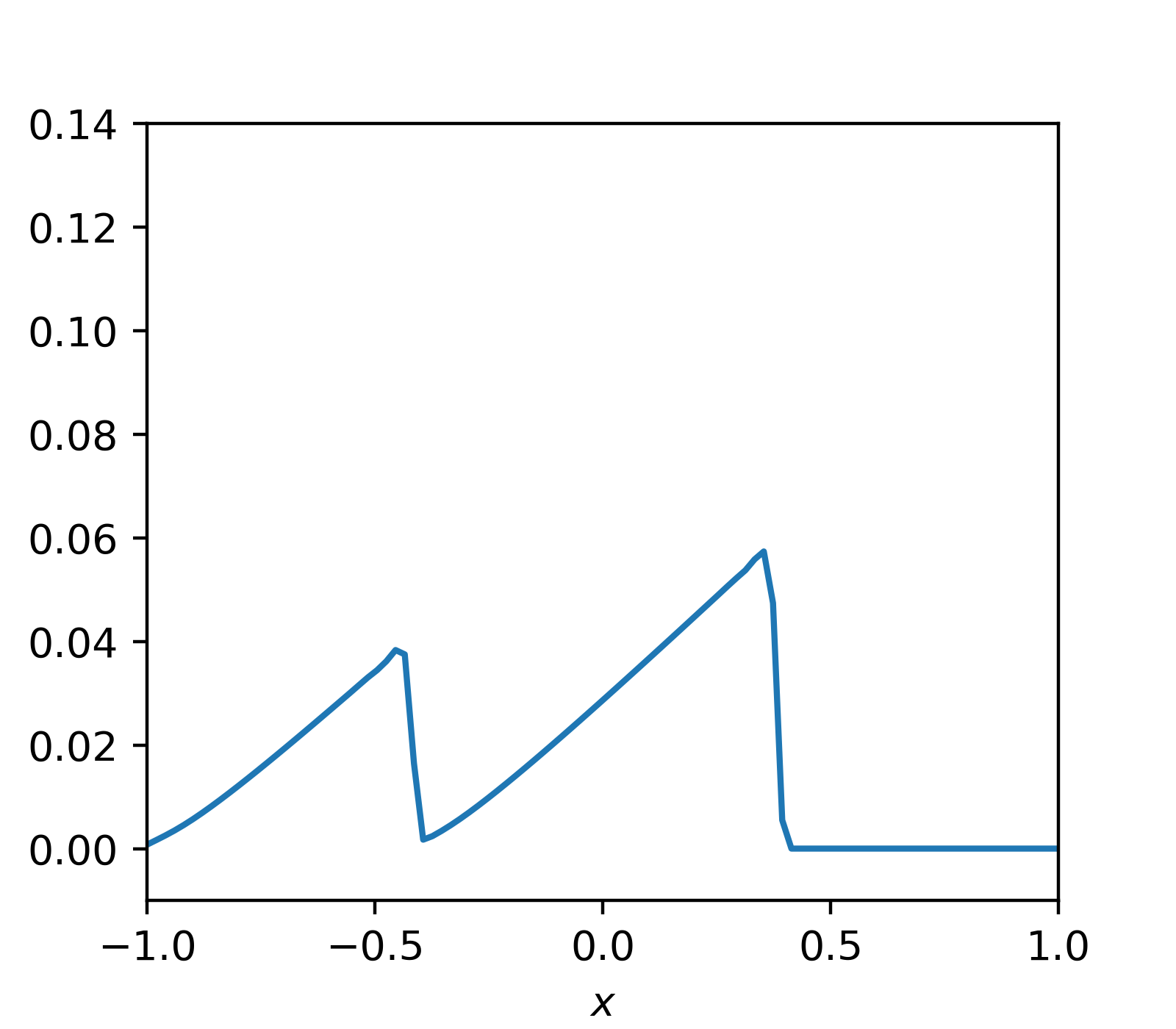}
        \\
    (c) & (d) \\
    \includegraphics[width=0.35\textwidth]{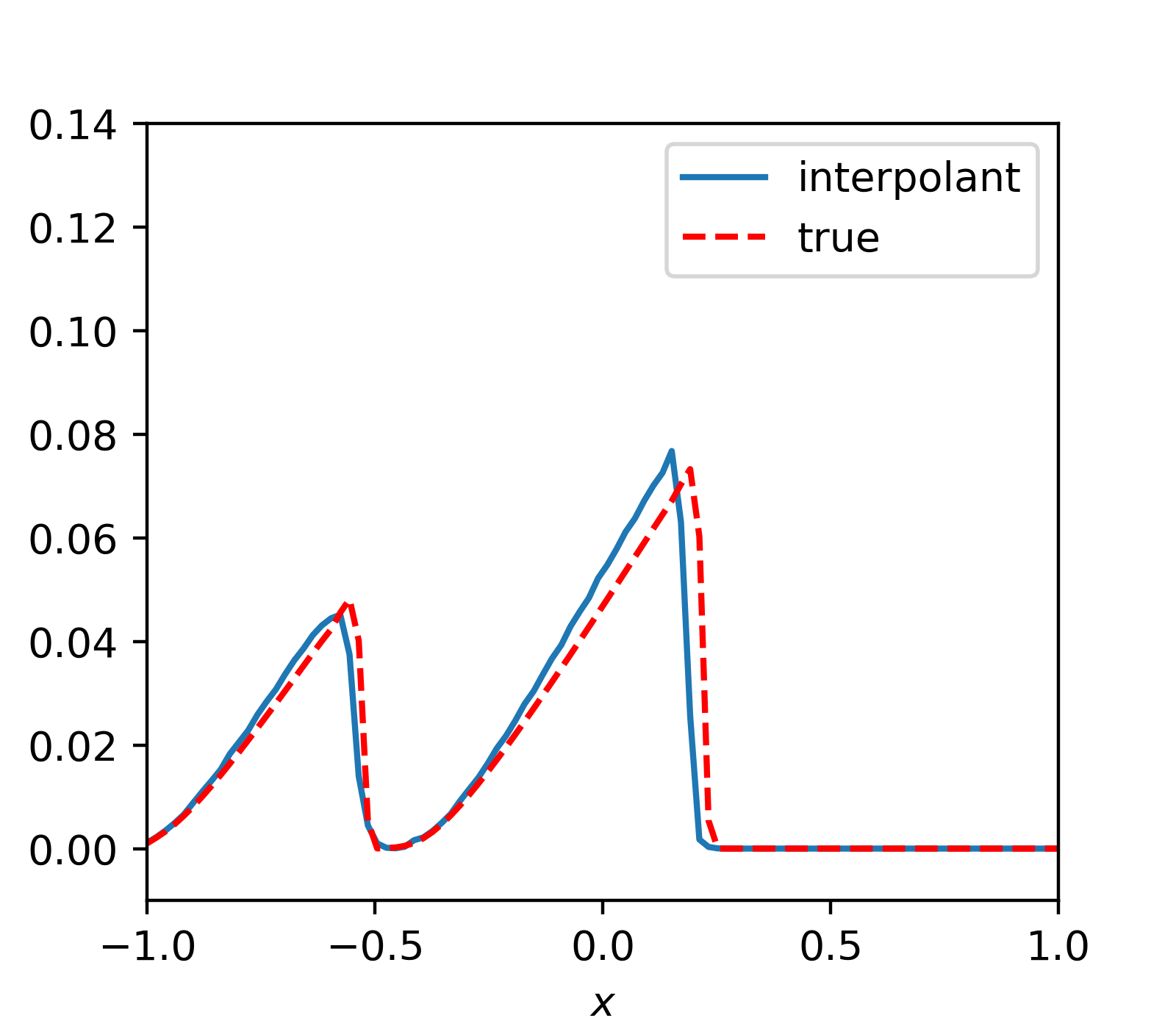}
        &
    \includegraphics[width=0.35\textwidth]{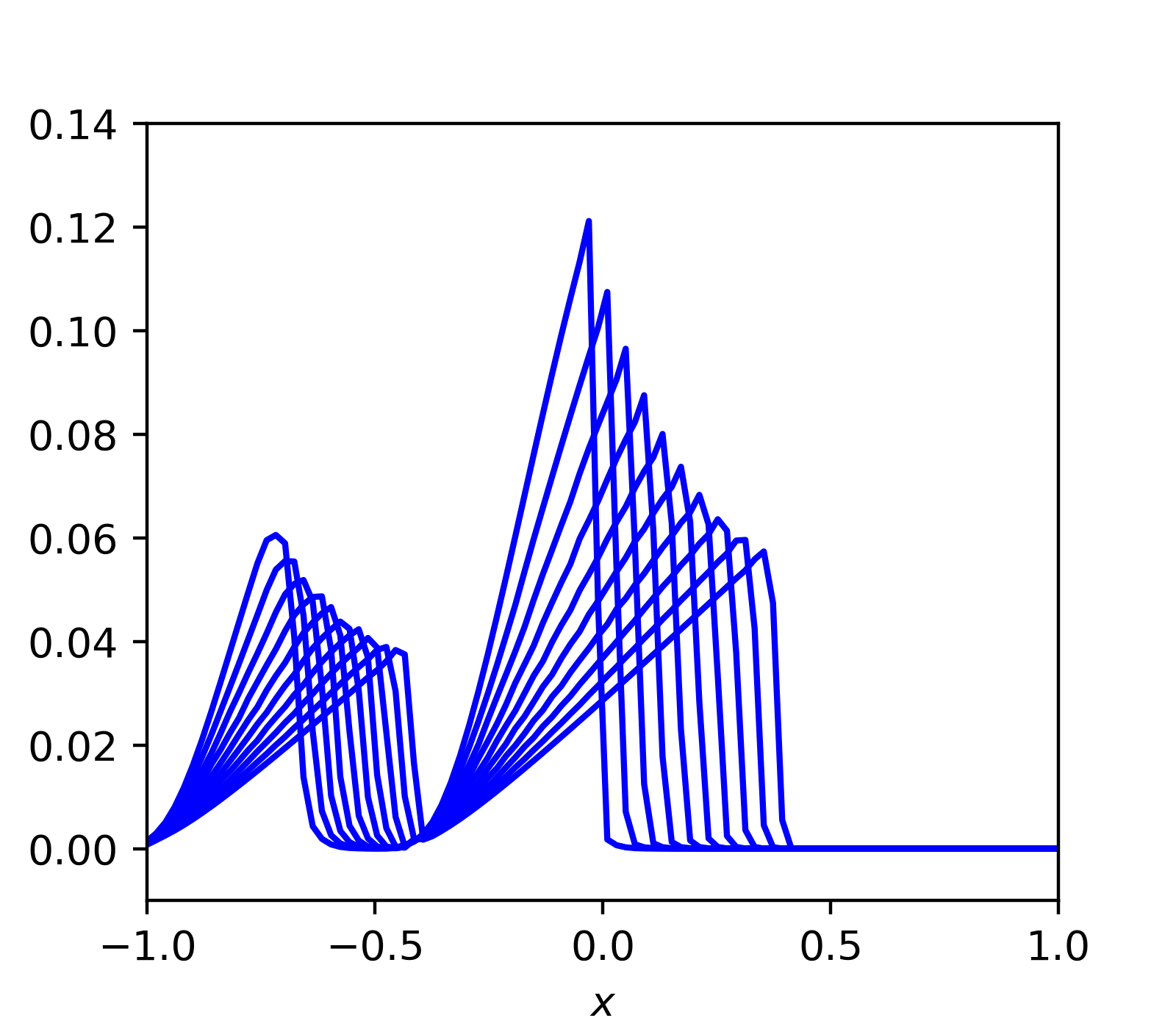}
\end{tabular}
\caption{(a) Solution $u$ to the Burger's equation \eqref{eq:acoustic1d} at
times $t_1 = 1$  and (b) $t_2 = 3$. (c) Displacement interpolant $\tilde{u}_t$
\eqref{eq:pressureI} to these two at $t = 2$ or $\lambda(t) = 0.5$, and (d)
$\tilde{u}_t$ at times $t = 1 + 0.2n$  for $n = 1,2, ... ,10$ (bottom right).}
\label{fig:burgers}
\end{figure}

\begin{figure}
\centering
\includegraphics[width=0.80\textwidth]{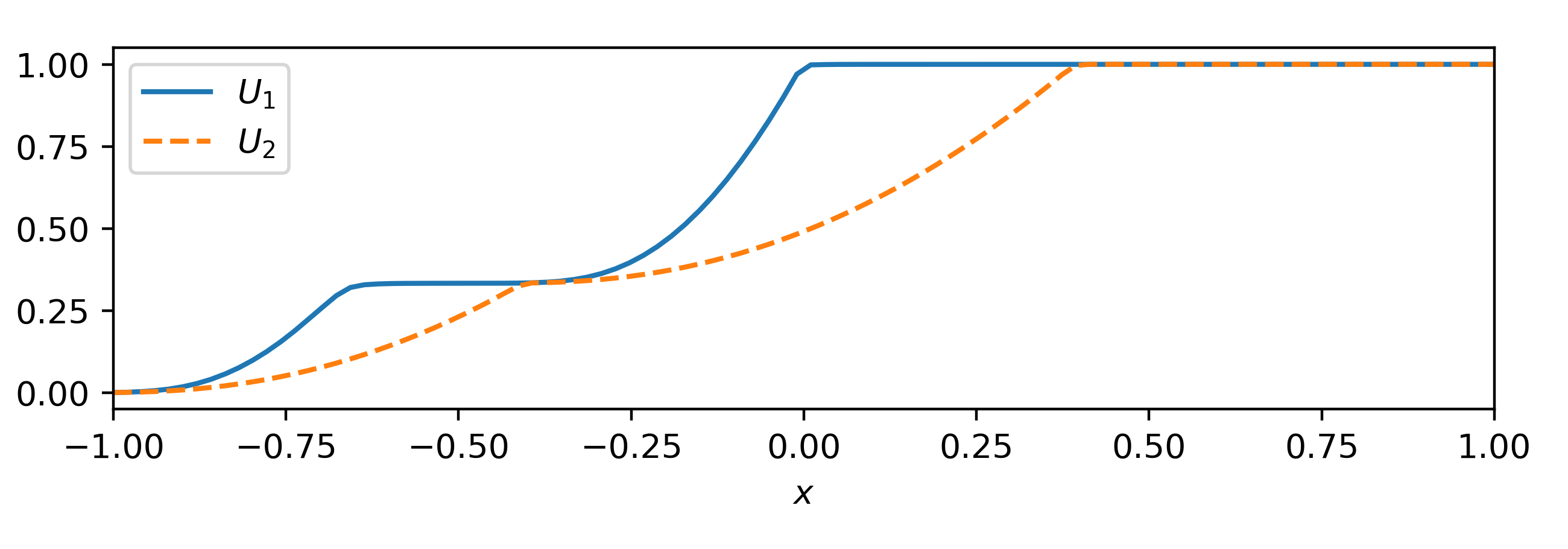}
\caption{The cumulative distribution functions (CDFs) of the solution $u$ to
the Burger's equation \eqref{eq:burgers} at times $t_1$ and $t_2$, denoted by
$U_1$ and $U_2$.}
\label{fig:burgersCDFs}
\end{figure}

In this example, we consider the interpolation of the sotluion to the 1D
Burger's equation,
\begin{equation}
u_t + u u_x = 0.
\label{eq:burgers}
\end{equation}
As in the previous section, we are interested in the change in the behavior of
the solution with respect to the paramter $\alpha = t$, the time variable.
Suppose we are given the solution at two times $t_1$ and $t_2$, ($u_1 =
u(x,t_1)$, $u_2 = u(x,t_2)$) shown in \cref{fig:burgers}(a),(c).  The
corresponding CDFs \eqref{eq:CDFs} are shown in \cref{fig:burgersCDFs}.  The
solution $u$ to the Burgers' equation develops shocks, and once a shock is
formed, the speed at which the shock travels is determined by the
Rankine-Hugoniot jump condition \cite{fvmbook}. Although we do not expect the
interpolant $\tilde{u}_t$ will agree with the solution $u$ at time $t$, unlike
in the previous example, we aim to recover a reasonable deformation map that
will provide a suitable approximation to the evolution of the shock.

The solutions at two different times $u_1$ and $u_2$ are shown in the top row
of \cref{fig:burgers}. The interpolant $\tilde{u}_t$ where $t = 2$ ($\lambda(t)
= 0.5$) is shown in the bottom left plot of the same figure. The dotted line
represents the true solution  $u(x,t)$.  The plot of $\tilde{u}_t$ for 10
different values $t=1 + 0.2n$ for $n = 1,2, ... , 10$.  is shown in the lower
right plot in the same figure.  It is easy to see that the speeds of the
propagation of the shock is linear with respect to time for the interpolant.
This is not true for the true solution, therefore the interpolant is only an
approximation whose quality depends on the distance to the nearest data point
in the paramter space of $\alpha = t$. However, this yields a far superior
approximation than the corresponding linear interpolant of the two.

\subsection{2D displacement interpolation}
In this section, we present an example of the interpolant \eqref{eq:ndinterp}
with spatial dimension $d=2$. The two functions $u_1(x,y)$ and $u_2(x,y)$ are
given as follows,
\begin{align*}
    u_1(x,y) &= \frac{3}{2} 
    \exp\left[
        - \frac{ (|x| + |y| - 0.75)^2}{2 \cdot (0.2)^2}
        - \frac{ (|x| + |y| - 0.5)^2}{2 \cdot (0.2)^2}
    \right], \\
    u_2(x,y) &= \frac{3}{2} 
    \exp\left[
        - \frac{(|x-0.5| + |y-0.25|)^2}{2 \cdot (0.05)^2}
    \right].
\label{eq:2d2fctns}
\end{align*}
The first function $u_1$ is concentrated near a diamond-shaped Gaussian hump
centered at the origin, where as $u_2$ is a much narrower diamond-shaped
Gaussian hump centered at the coordinates $(x,y) = (0.5, 0.25)$. Let us suppose
$u_1$ is the function at parameter $\alpha = 0$ whereas $u_2$ is a function at
parameter $\alpha = 1$.

We compute the interpolant \eqref{eq:ndinterp},
\begin{equation}
    \tilde{u}_\alpha = \cI_2(\alpha; u_1, u_2)
\end{equation}
for values $\alpha = 0.25, 0.5, 0.75$.

\begin{figure}
\centering
\begin{tabular}{rcc}
    \rotatebox{90}{$\qquad \quad \alpha = 0$} 
    &
    \includegraphics[width=0.3\textwidth]{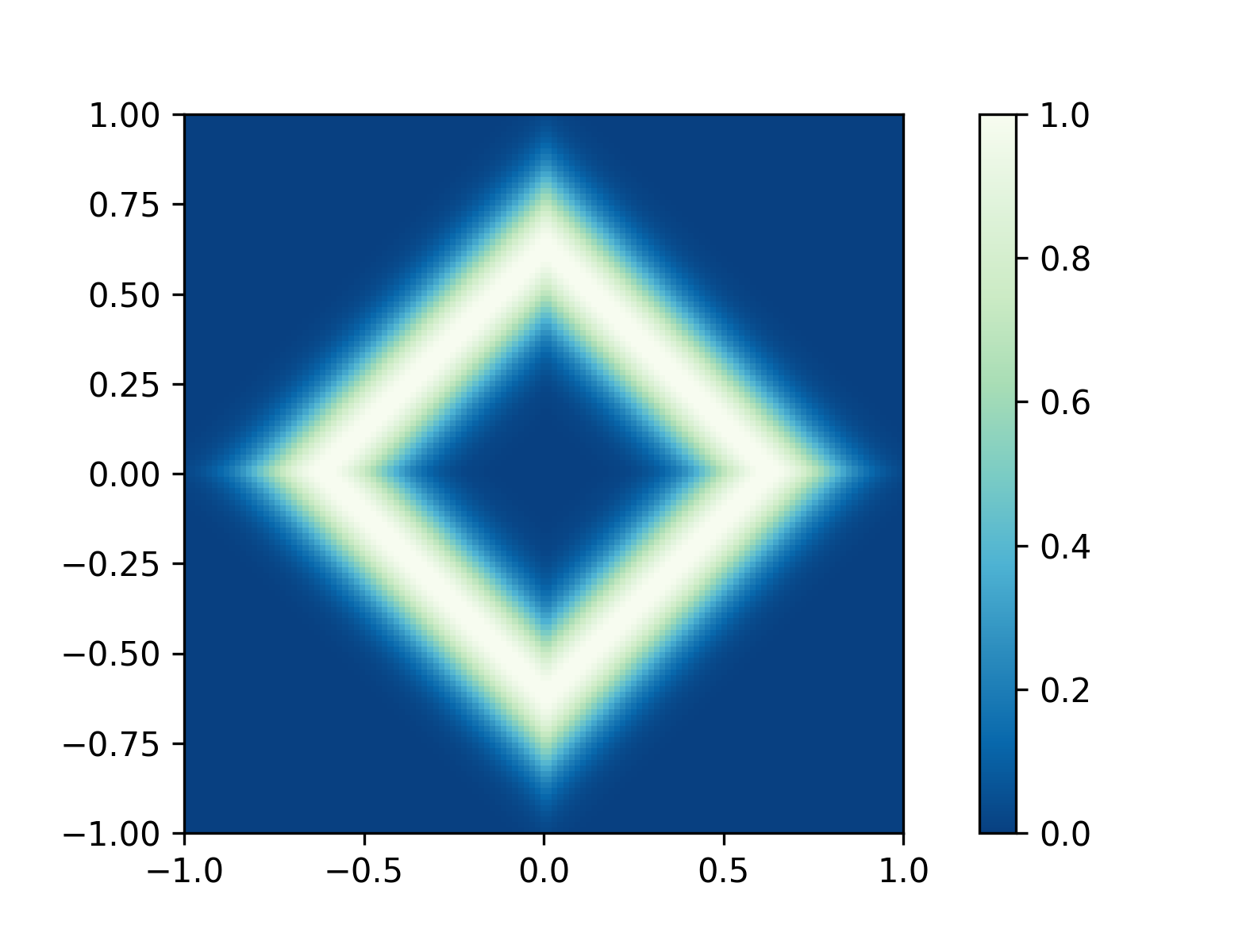}
    &
    \includegraphics[width=0.4\textwidth]{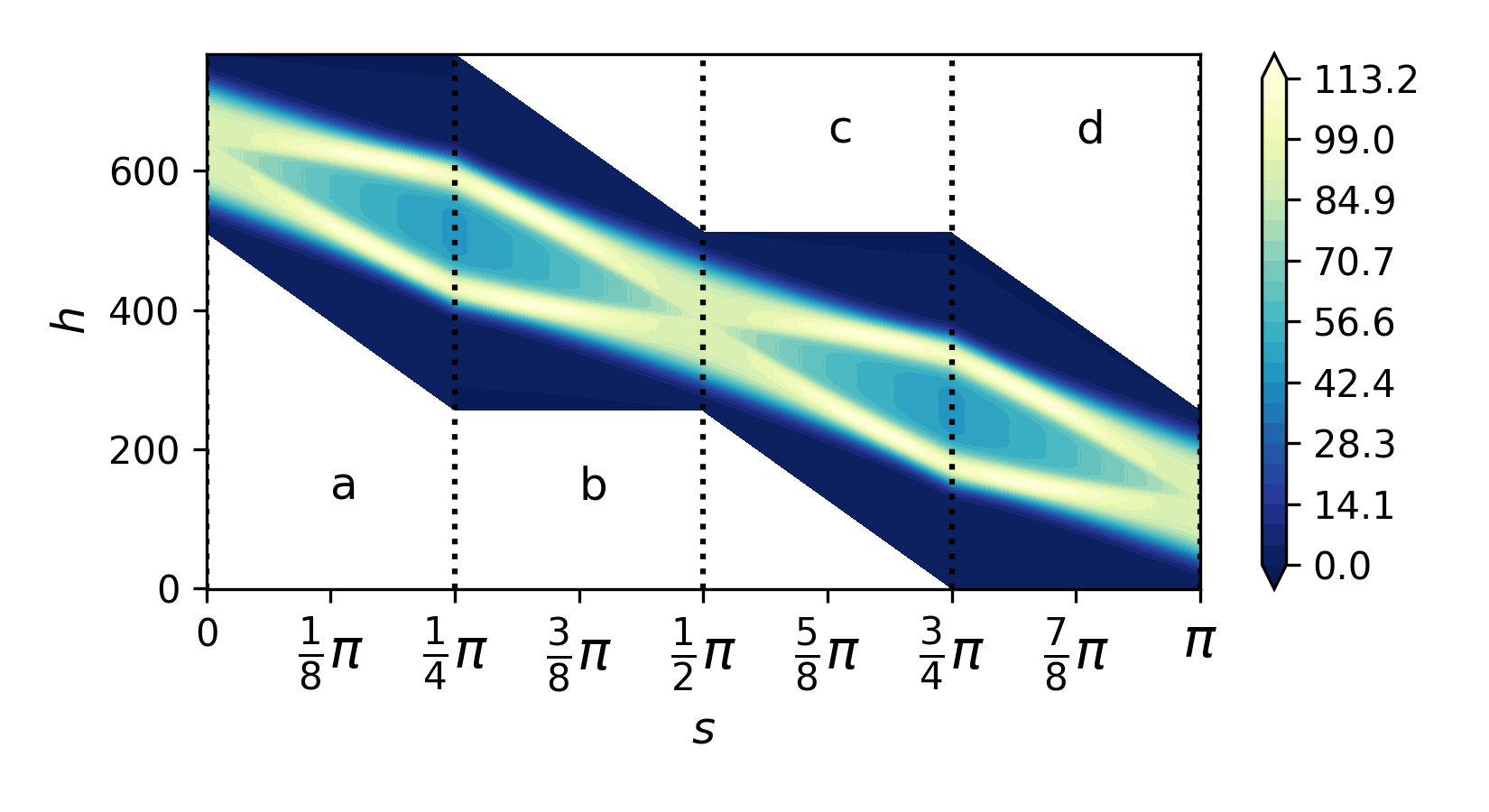}
    \\
    \rotatebox{90}{$\qquad \quad \alpha = 0.25$} 
    &
    \includegraphics[width=0.3\textwidth]{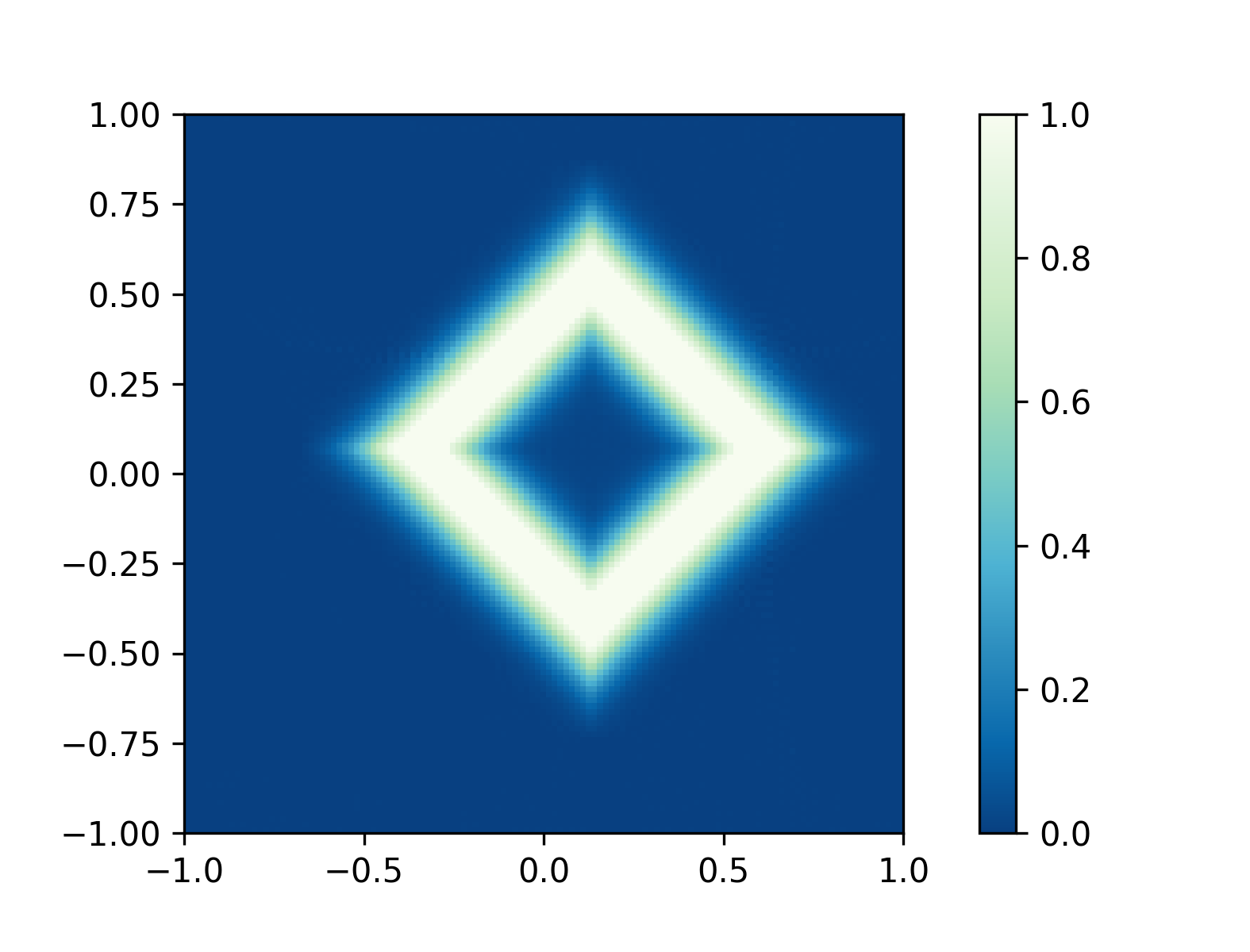}
    &
    \includegraphics[width=0.4\textwidth]{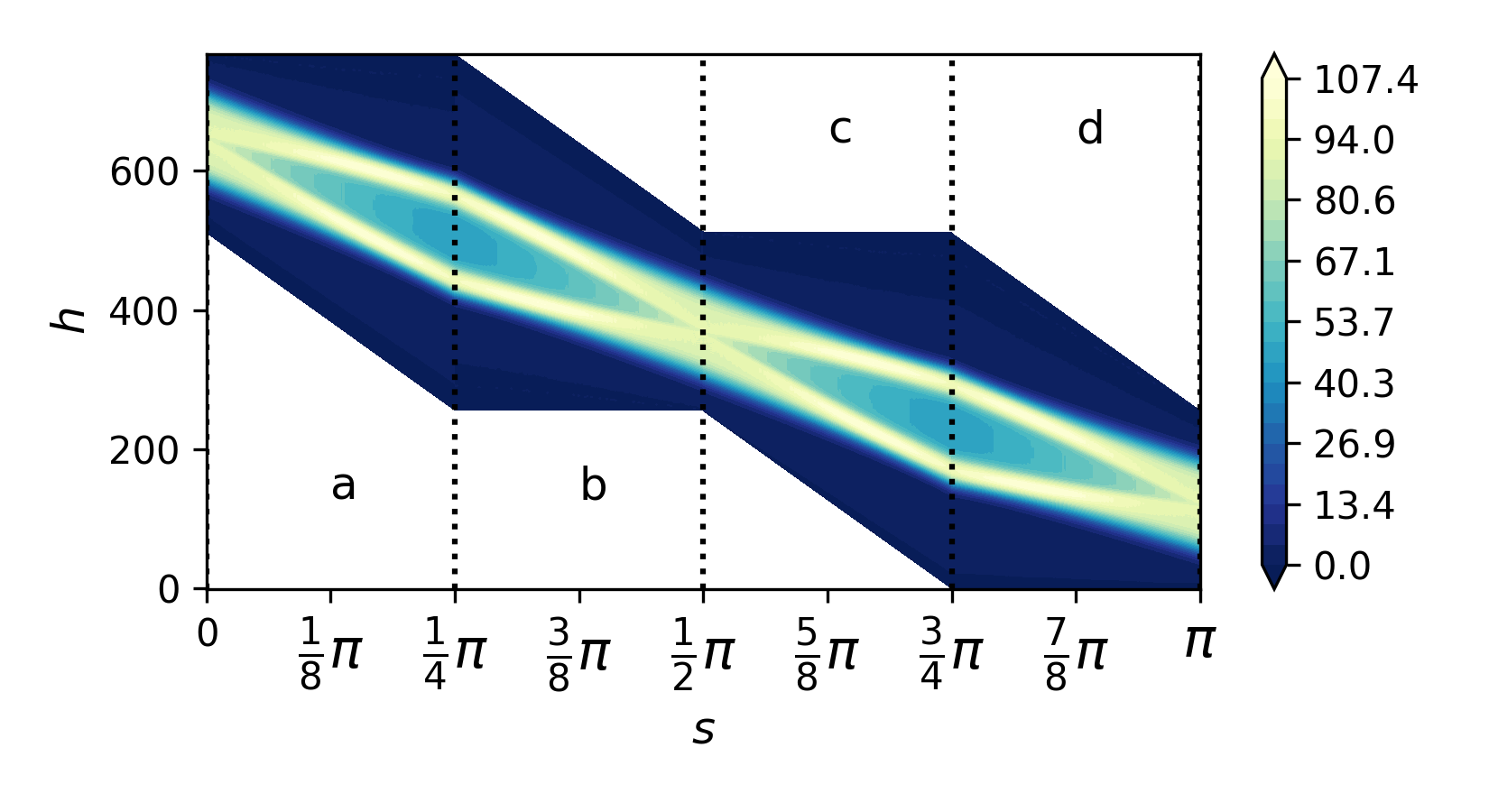}
    \\
    \rotatebox{90}{$\qquad \quad \alpha = 0.5$} 
    &
    \includegraphics[width=0.3\textwidth]{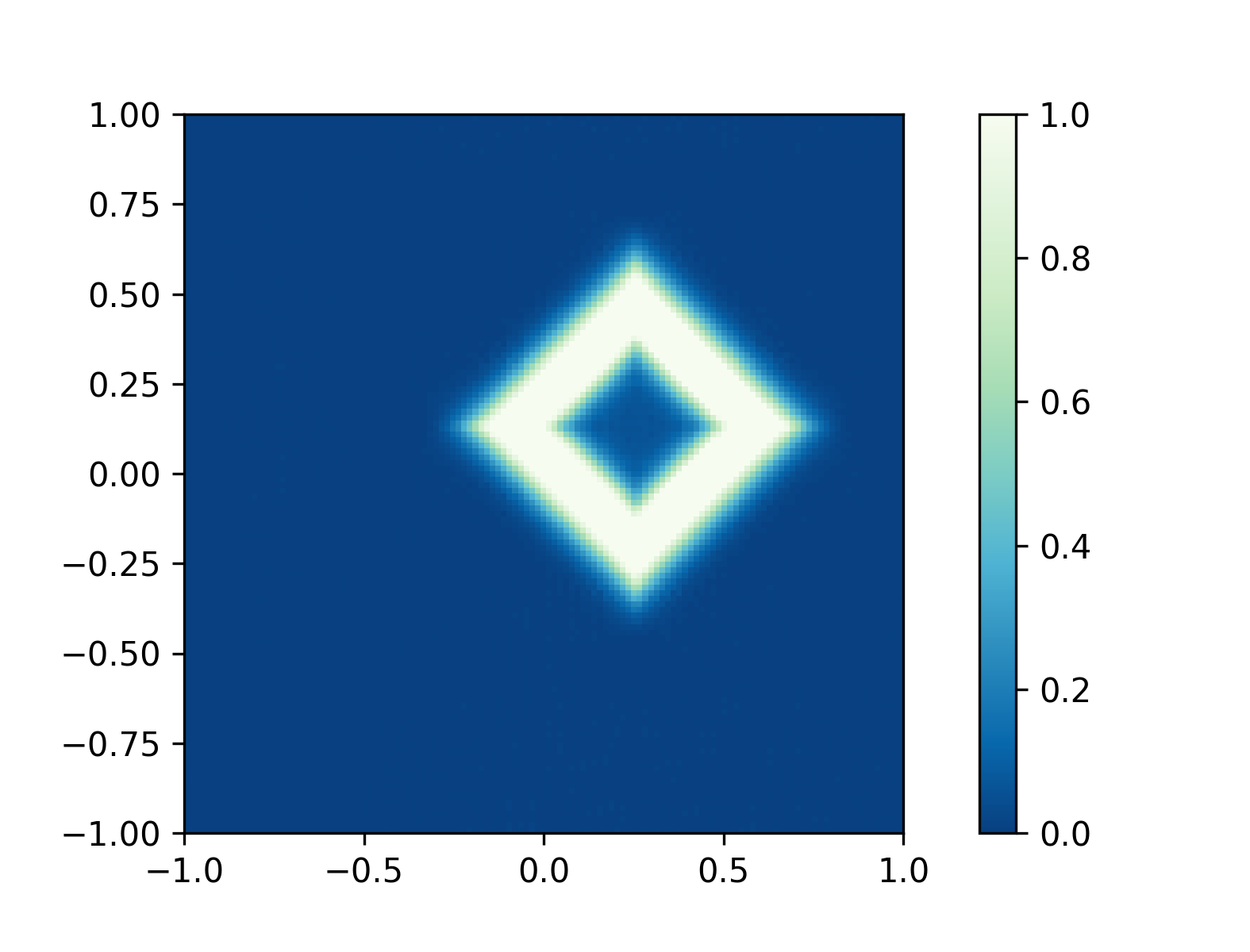}
    &
    \includegraphics[width=0.4\textwidth]{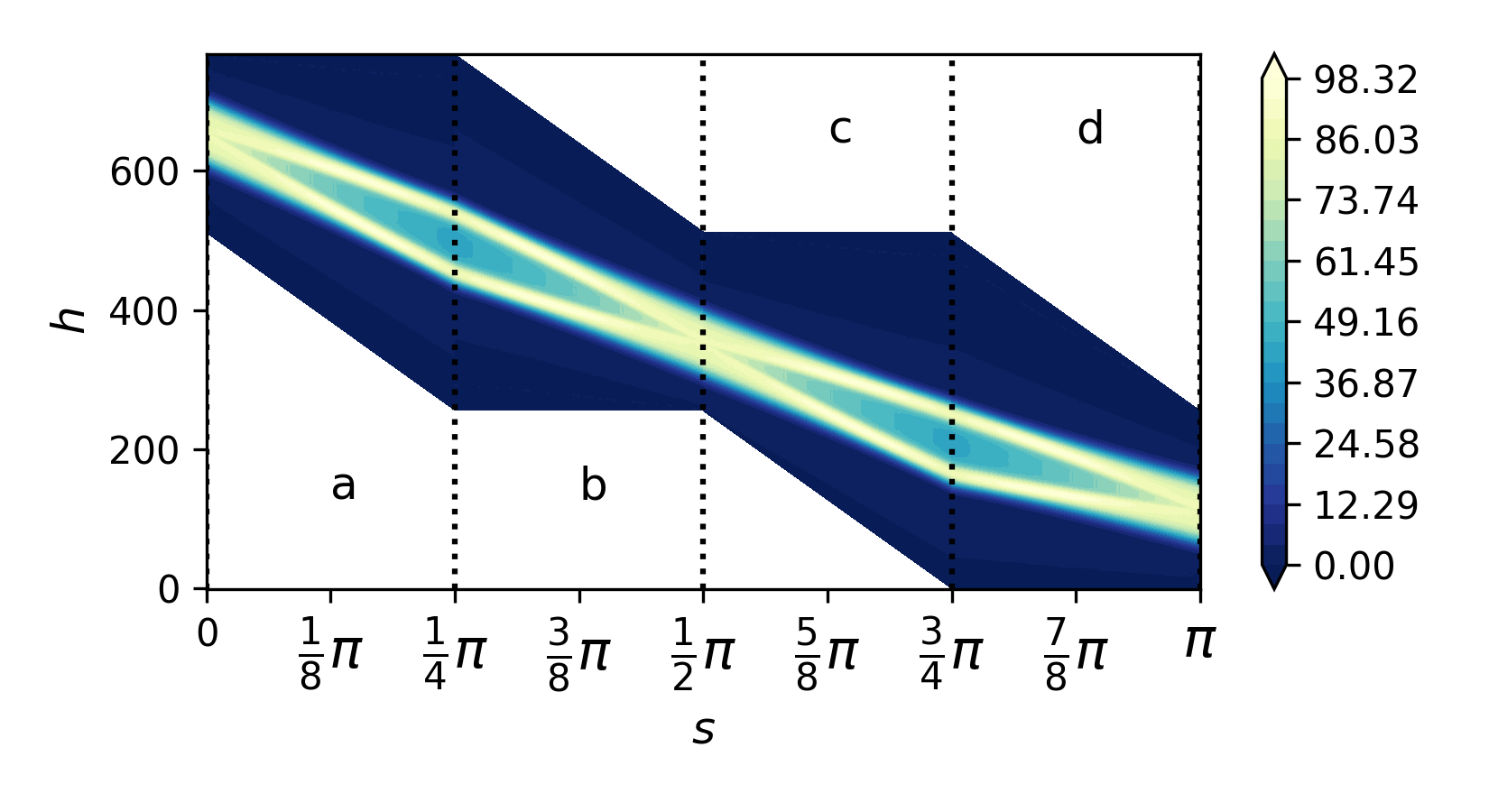}
    \\
    \rotatebox{90}{$\qquad \quad \alpha = 0.75$} 
    &
    \includegraphics[width=0.3\textwidth]{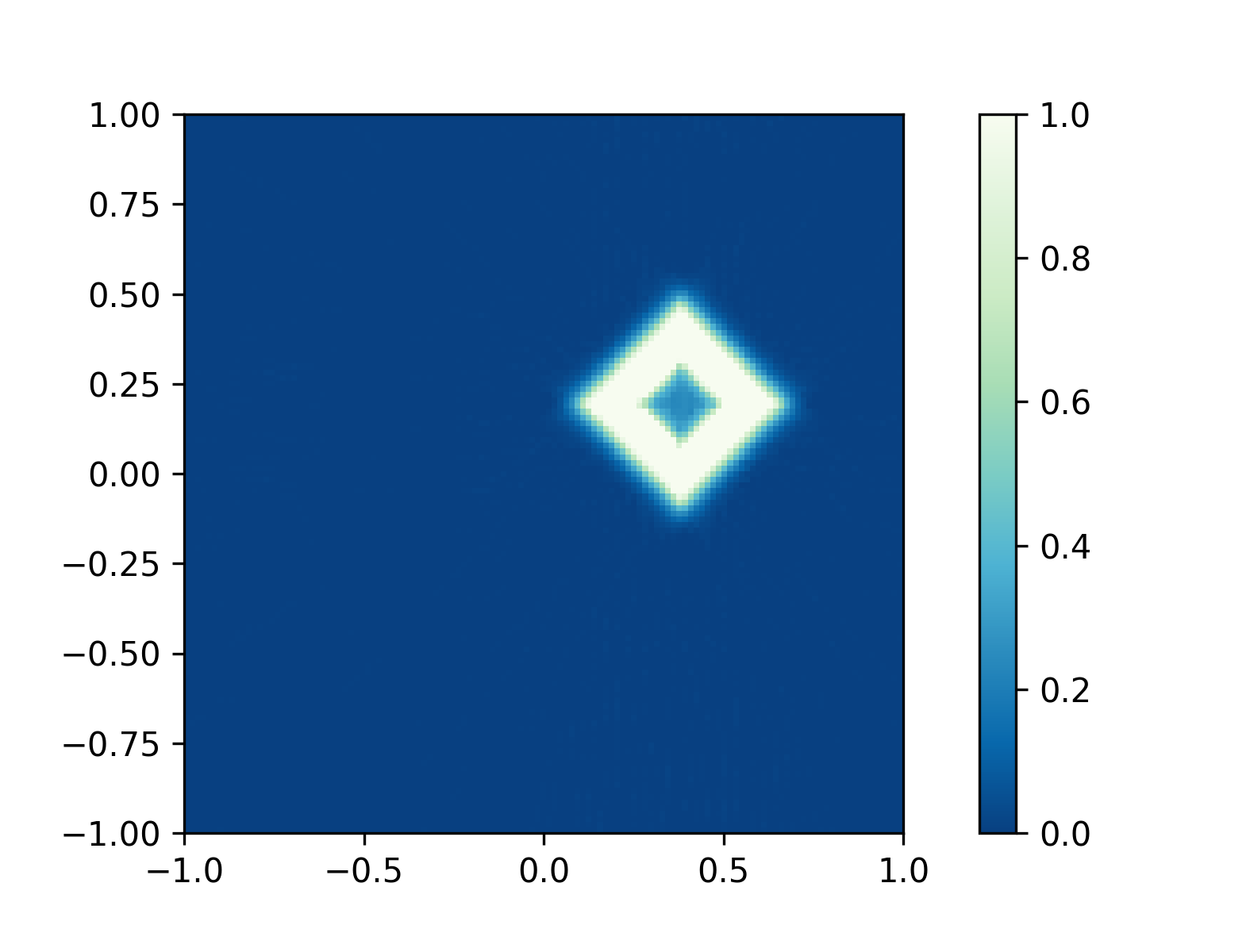}
    &
    \includegraphics[width=0.4\textwidth]{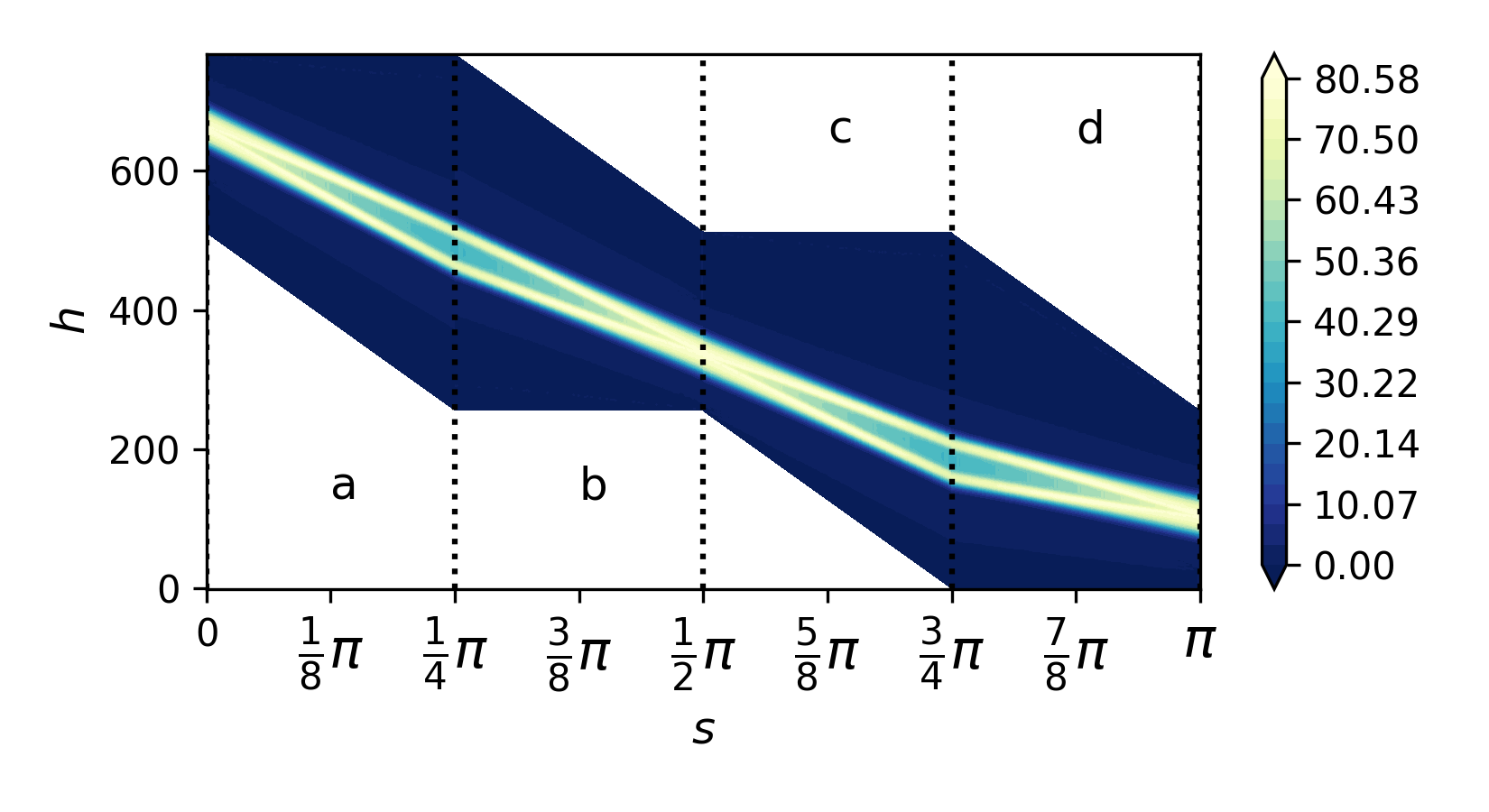}
    \\
    \rotatebox{90}{$\qquad \quad \alpha = 1.0$} 
    &
    \includegraphics[width=0.3\textwidth]{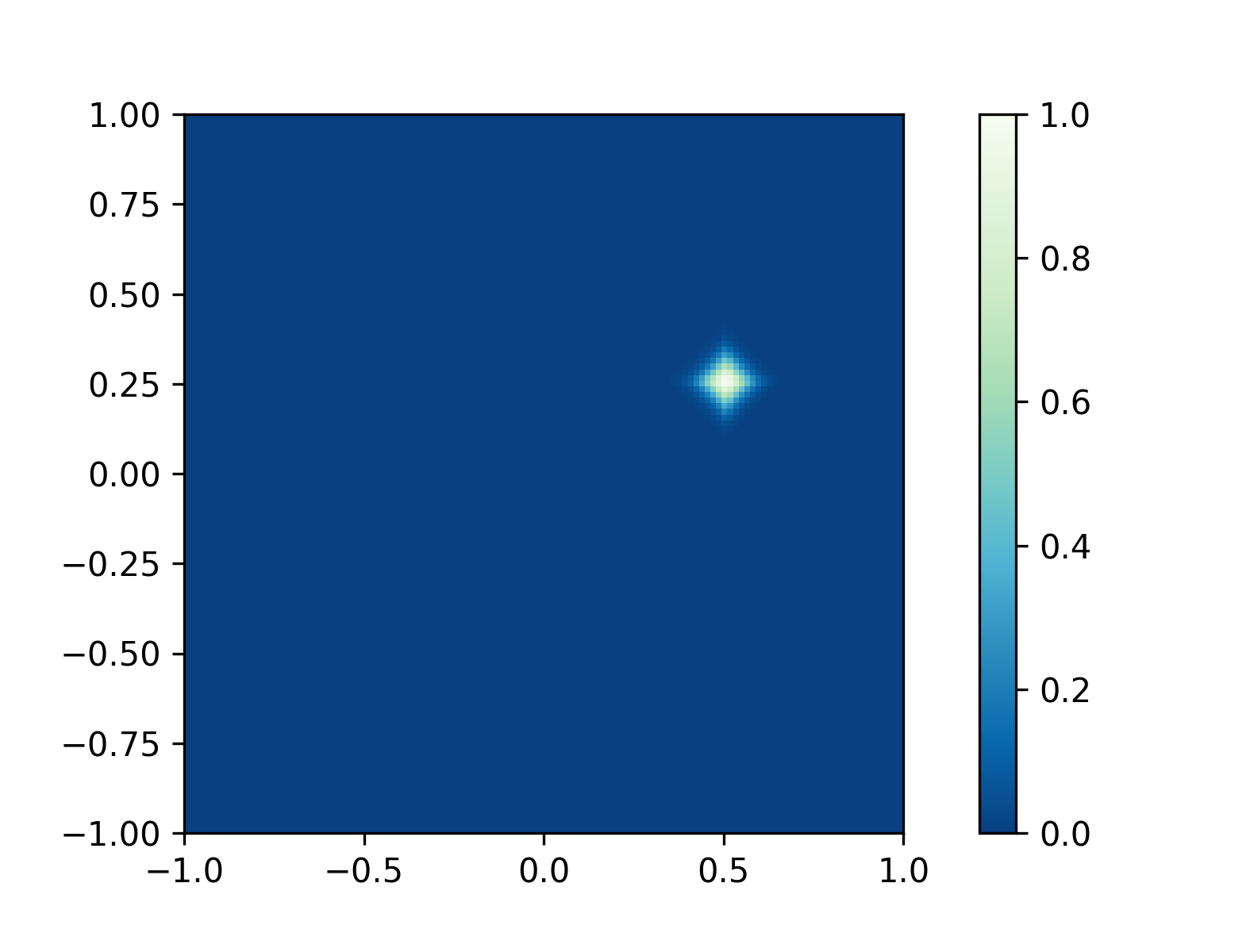}
    &
    \includegraphics[width=0.4\textwidth]{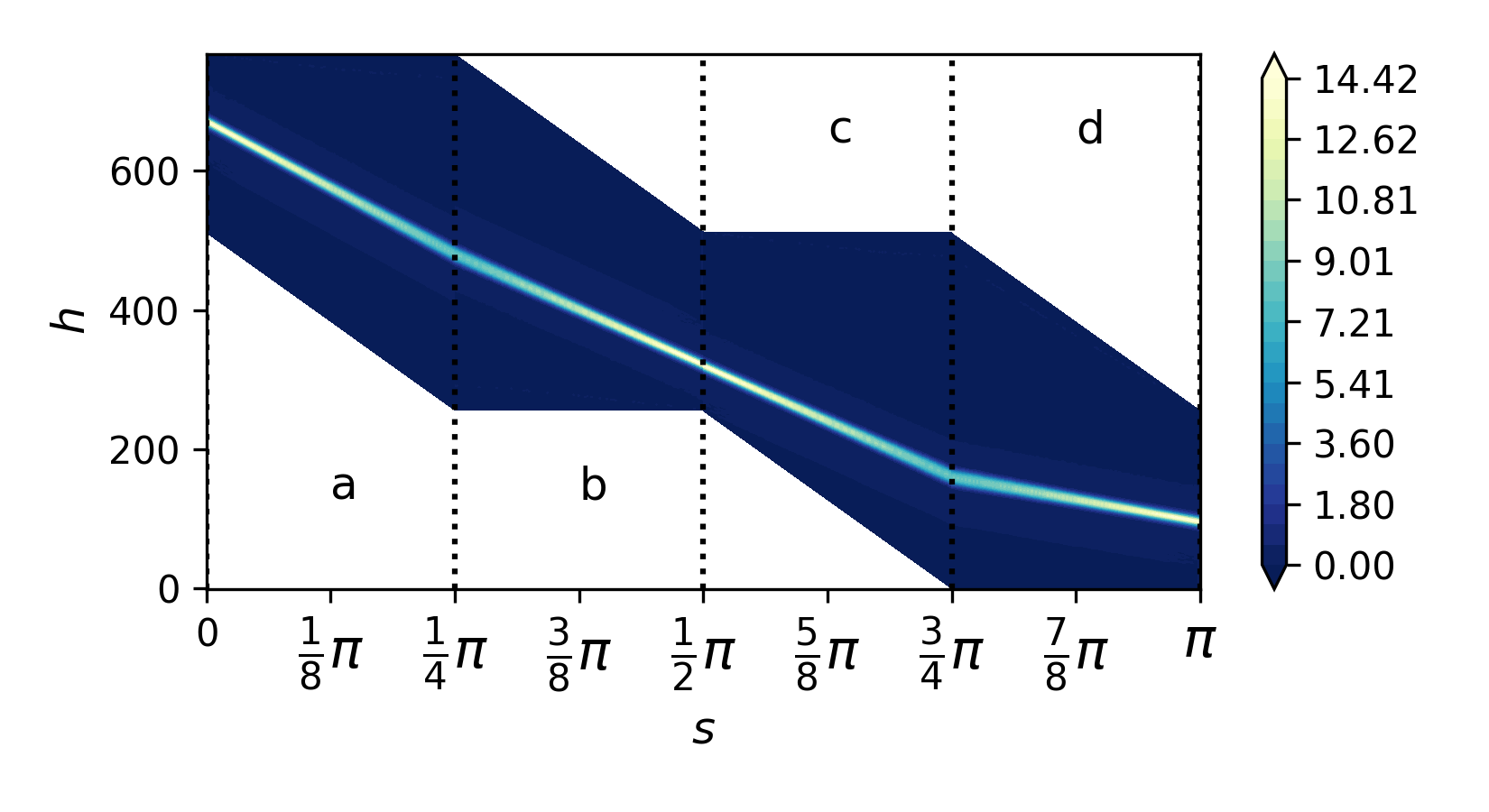}
\end{tabular}
\caption{2-dimensional displacement interpolant between two functions $u_1$ and
$u_2$ given in \eqref{eq:2d2fctns}, shown together with their discrete Radon
transform. The functions themselves are given in the first row $\alpha = 0$ and
the last row $\alpha =1$.}
\label{fig:2dinterp}
\end{figure}

The two functions and the interpolants are shown in \cref{fig:2dinterp}.  The
interpolant discovers the translation map, that shifts the center of mass of
$u_1$ towards that of $u_2$, while shrinking the region of concentrated mass
towards the concentrated hump. In the transformed variables $(h,s)$, this
transport map is reduced to a 1D monotone rearrangement, and the smooth
displacement interpolant we have obtained is a low-dimensional for each slice
of the transformed variables for fixed $\boldsymbol{\omega}$.  This is clear
from \eqref{eq:dinterp-slice}.

\begin{figure}
\centering
\begin{tabular}{c}
    \includegraphics[width=0.9\textwidth]{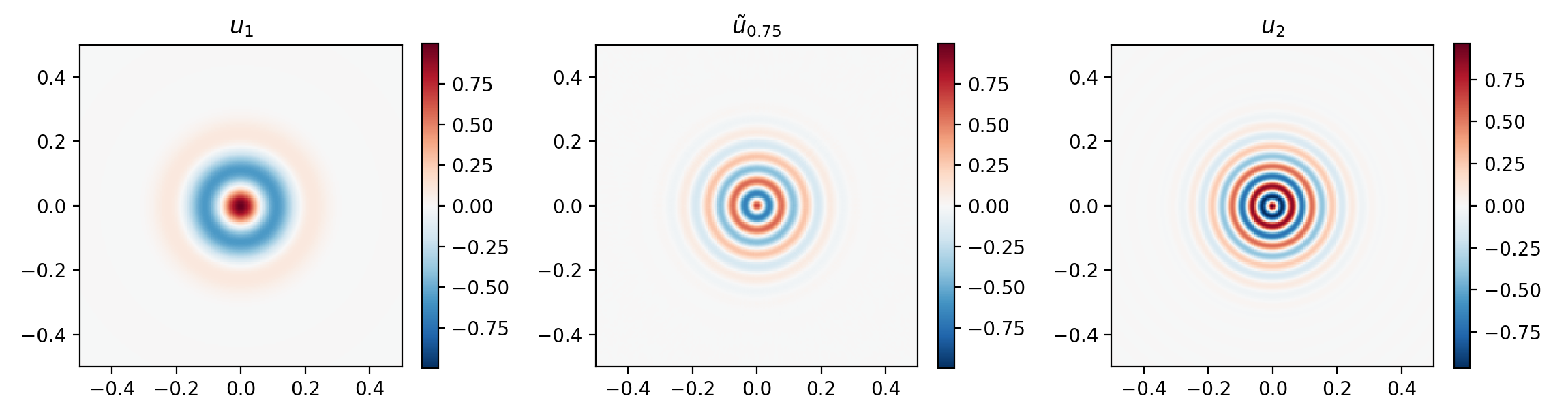}
    \\
    \includegraphics[width=0.9\textwidth]{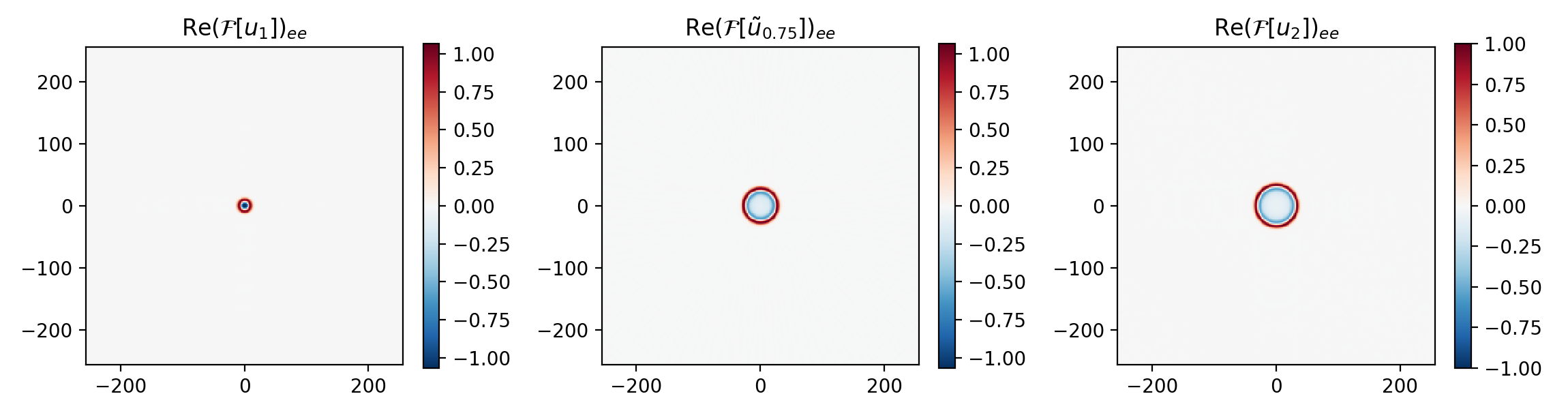}
\end{tabular}
\caption{ \added{2D displacement interpolant between two oscillating functions
$u_1$ and $u_2$ given in \eqref{eq:osc_2d2fctns} and the interpolant
$\tilde{u}_{0.75} = \cI_\cT(0.75;u_1,u_2)$ where $\cT$ is given by
\cref{eq:oscT} (top row), shown together with the corresponding even-even
component of $\Real(\cP \cF[u])$ \cref{eq:subindices} (bottom row). The plots
have been normalized with respect to the maximum value for ease of
comparison.}}
\label{fig:osc_2d}
\end{figure}

\begin{figure}
\centering
\includegraphics[width=0.9\textwidth]{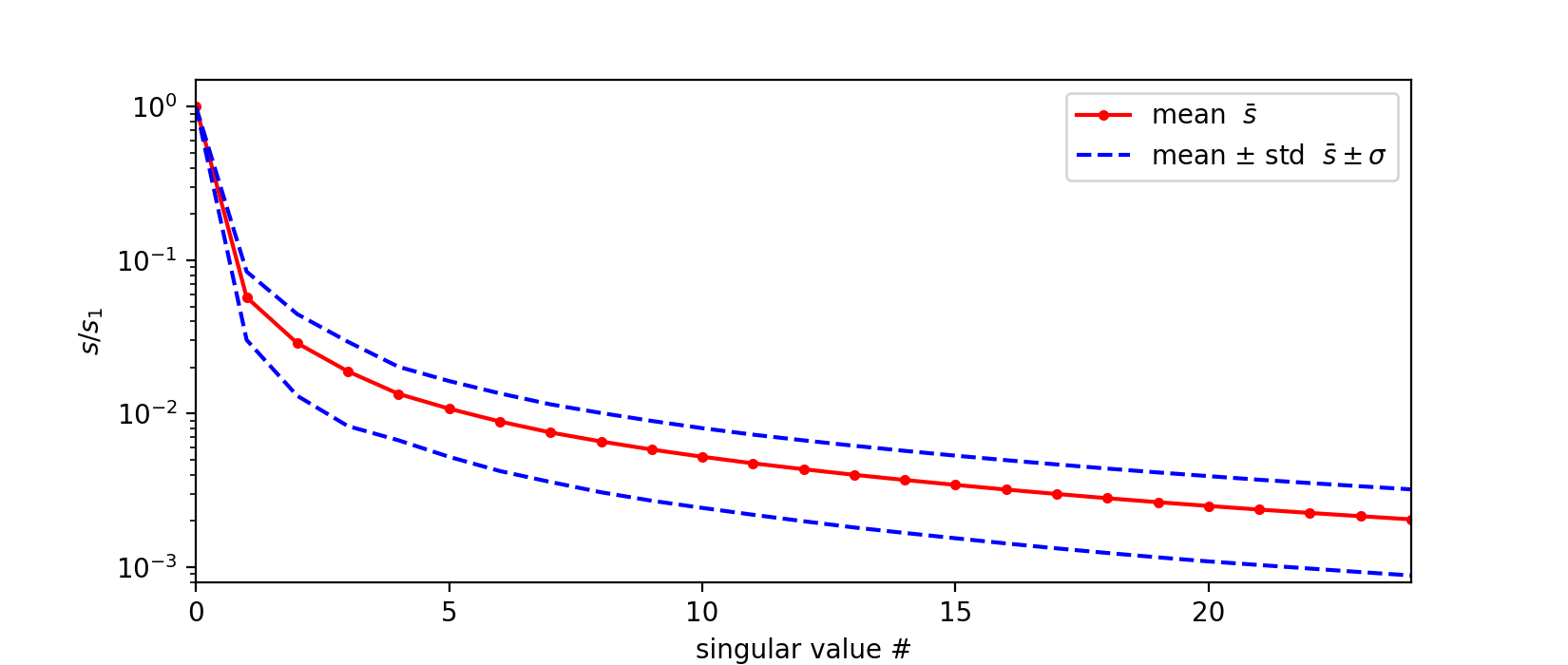}
\caption{\added{Mean and standard deviation of the singular values of the
transport maps, over all angles $\boldsymbol{\omega} \in S^1$ 
\cref{eq:meanvar_singvals}.}}
\label{fig:osc_2d_singvals}
\end{figure}

\added{
\subsection{2D oscillating function} \label{sec:osc}
In this example, we will interpolate between 2D oscillatory functions.  Let us
consider the radial function with a free parameter $k$,
\beq
u(r;k) = 
\exp\left[ - \left( \frac{r^2}{2\sigma^2}\right) \right] 
\cos\left( k \pi  r  \right),
\quad
\text{ where }
r = \sqrt{x^2 + y^2}, 
\quad
0\le x,y \le 1.
\label{eq:osc2d_fctn}
\eeq
Let us set $\sigma^2$ as $0.0125$ and the two functions $u_1$ and $u_2$ as
\beq
u_1(x,y) := u(r;k_1), \quad k_1 = 8
\quad \text{ and } \quad
u_2(x,y) := u(r;k_2), \quad k_2 = 16.
\label{eq:osc_2d2fctns}
\eeq
We will apply the displacement interpolation of the type $\cI_\cT$
\cref{eq:dinterp_transform} in which we will choose the transform $\cT$ as,
\beq
\cT = \cR \otimes \cP \cF,
\label{eq:oscT}
\eeq
where $\cF$ is the Fourier transform $\cP$ rearranges the Fourier transform,
\beq
\begin{aligned}
\cP \cF[u] = &\left[
               \Real{(\cF[u])}_{\textrm{oo}},
               \Real{(\cF[u])}_{\textrm{eo}},
               \Real{(\cF[u])}_{\textrm{oe}}, 
               \Real{(\cF[u])}_{\textrm{ee}},\right. \\
      &\left.   \Imag{(\cF[u])}_{\textrm{oo}}, 
               \Imag{(\cF[u])}_{\textrm{eo}}, 
               \Imag{(\cF[u])}_{\textrm{oe}}, 
               \Imag{(\cF[u])}_{\textrm{ee}} \right] 
\end{aligned}
\label{eq:subindices}
\eeq
where the subscript indices $\{\textrm{o},\textrm{e}\}$ indicate the
restriction to either odd or even sub-indices of the 2D arrays,
and $\cR \otimes$ acts on the individual components of \cref{eq:subindices}.
The interpolation $\cI_\otimes$ appearing in $\cI_\cT = \cT^{-1} \cI_\otimes
\cT$ acts on the $s$-variable of the Radon transform, as before.  

The functions $u_1,u_2$ as well as the resulting interpolant $\tilde{u}_\alpha$
with $\alpha =0.75$, along with the even-even component
$\Real{(\cF[\cdot])}_{\textrm{ee}}$ are shown in \cref{fig:osc_2d}. The
interpolant is able to capture the oscillatory nature of the functions
\cref{eq:osc2d_fctn}. Note that their Fourier transforms are wave-like and that
the location of high amplitudes in the frequency space correspond to the
frequency of the oscillation.

Next, let us examine the decay of the singular values of the transport maps
computed for different angles of $\cT[\cdot]$ for a larger sample of functions
of the form \cref{eq:osc2d_fctn} with more diverse values of $k$. To be more
clear, we consider the SVD of the transport maps between the 1D functions
\beq
\begin{aligned}
\cT[v](\cdot,\boldsymbol{\omega})  
    &= \cR \otimes \cP \cF [v](\cdot,\boldsymbol{\omega}) \\
&\quad \text{ where } 
v \in \cU := \{u(r;k ): k=8 + .5j, \, j=0,1,..., 49\}.
\end{aligned}
\eeq
We display the sample mean and the standard deviation over all angles
$\boldsymbol{\omega}$ of the singular values $s_j(\boldsymbol{\omega})$
computed from these maps in \cref{fig:osc_2d_singvals},  that is,
\beq
\bar{s}_j := \mathbb{E}_{\boldsymbol{\omega}} 
                    \left[ s_j(\boldsymbol{\omega})\right]
\quad \text{ and } \quad
\sigma_j := \sqrt {\textrm{Var}_{\boldsymbol{\omega}}
                        \left[ s_j(\boldsymbol{\omega}) \right]}.
\label{eq:meanvar_singvals}
\eeq
There is now a clear decay of singular values, whereas it is straightforward to
show that the original functions are nearly orthogonal to each other.
}

\added{
\subsection{Displacement interpolation for oscillatory functions in 2D}}

\section{Conclusion and future work} 
We have introduced a displacement interpolation scheme based on monotone
rearrangement solution to the Monge-Kantorovich problem in 1D. We extended this
to the case when the given functions have arbitrary sign, as well as when
multiple functions are given.  We then showed that this scheme can be naturally
generalized to multiple spatial dimensions through the use of the Radon
transform, resulting in a generalization of the Lax-Philips transform.

The interpolation \added{allowed the authors to achieve the dimensionality
reduction} of wave-like phenomena \added{in 1D scalar conservation laws
\cite{rim18}}, and we will investiage its application in the model reduction of
hyperbolic PDEs in a multi-dimensional setting \added{in a future work.}
\added{The interpolation method is by no means restricted to hyperbolic PDEs,
however, and may well complement existing dimensionality reduction methods in
various applications where standard methods are unsuccessful.}

\bibliographystyle{siamplain}

\end{document}

%% file: macros.tex
\newcommand{\cF}{\mathcal{F}}

\newcommand{\cO}{\mathcal{O}}
\newcommand{\cT}{\mathcal{T}}

\newcommand{\cP}{\mathcal{P}}
\newcommand{\cR}{\mathcal{R}}

\newcommand{\cU}{\mathcal{U}}

\newcommand{\cY}{\mathcal{Y}}

\newcommand{\beq}{\begin{equation}}
\newcommand{\eeq}{\end{equation}}
\def\bals#1\eals{\begin{align*} #1 \end{align*}}
\def\bal#1\eal{\begin{align} #1 \end{align}}

\newcommand\Dom\Omega

\newcommand\RR{\mathbb{R}}

\newcommand\Lap\Delta

\def\bpde#1\epde{\[\left\{\begin{aligned}#1\end{aligned}\right. \]}
\def\inbpde#1\inepde{\left\{\begin{aligned}#1\end{aligned}\right.}
\def\binpde#1\einpde{\left\{\begin{aligned}#1\end{aligned}\right.}

\newcommand\Norm[2]{\left\lVert { #1 } \right\rVert_{#2}}

\def\cI{\mathcal{I}}
\def\cL{\mathcal{L}}
\def\cU{\mathcal{U}}
\def\cV{\mathcal{V}}

\def\cR{\mathcal{R}}

\def\bu{\mathbf{u}}

\def\cU{\mathcal{U}}

\def\bfU{\mathbf{U}}

\def\bfU{\mathbf{U}}

\def\b0{\mathbf{0}}

\def\bbmat{\begin{bmatrix}[r]}
\def\ebmat{\end{bmatrix}}

\def\balpha{\boldsymbol\alpha}

\newcommand{\barr}{\begin{array}}
\newcommand{\ea}{\end{array}}
\newcommand{\bea}{\begin{eqnarray}}
\newcommand{\eea}{\end{eqnarray}}
\newcommand{\bt}{\begin{table}}
\newcommand{\et}{\end{table}}

\DeclareMathOperator\supp{supp}

\DeclareMathOperator\Real{Re}
\DeclareMathOperator\Imag{Im}

\theoremstyle{plain}
\theoremstyle{definition}

\numberwithin{equation}{section}

\newcommand{\TheTitle}{Displacement interpolation using monotone rearrangement}

\newcommand{\remove}[1]{#1}%
\newcommand{\added}[1]{#1}%

%% file: mr.bbl
\begin{thebibliography}{10}

\bibitem{amsallem}
{\sc R.~Abgrall and D.~Amsallem}, {\em Robust model reduction by ${L}^1$-norm
  minimization and approximation via dictionaries: Application to linear and
  nonlinear hyperbolic problems}, Advanced Modeling and Simulation in
  Engineering Sciences, 3 (2016), p.~1.

\bibitem{agueh11}
{\sc M.~Agueh and G.~Carlier}, {\em Barycenters in the {W}asserstein space},
  SIAM Journal on Mathematical Analysis, 43 (2011), pp.~904--924.

\bibitem{benamou15}
{\sc J.~Benamou, G.~Carlier, M.~Cuturi, L.~Nenna, and G.~Peyré}, {\em
  Iterative {B}regman projections for regularized transportation problems},
  SIAM Journal on Scientific Computing, 37 (2015), pp.~A1111--A1138.

\bibitem{Benamou}
{\sc J.-D. Benamou and Y.~Brenier}, {\em A computational fluid mechanics
  solution to the monge-kantorovich mass transfer problem}, Numerische
  Mathematik, 84 (2000), pp.~375--393.

\bibitem{bolley05}
{\sc F.~Bolley, Y.~Brenier, and G.~Loeper}, {\em Contractive metrics for scalar
  conservation laws}, Journal of Hyperbolic Differential Equations, 02 (2005),
  pp.~91--107.

\bibitem{bonneel15}
{\sc N.~Bonneel, J.~Rabin, G.~Peyr{\'e}, and H.~Pfister}, {\em Sliced and
  {R}adon {W}asserstein barycenters of measures}, Journal of Mathematical
  Imaging and Vision, 51 (2015), pp.~22--45.

\bibitem{bonneel11}
{\sc N.~Bonneel, M.~van~de Panne, S.~Paris, and W.~Heidrich}, {\em Displacement
  interpolation using lagrangian mass transport}, ACM Trans. Graph., 30 (2011),
  pp.~158:1--158:12.

\bibitem{drt}
{\sc M.~L. Brady}, {\em A fast discrete approximation algorithm for the {R}adon
  transform}, SIAM Journal on Computing, 27 (1998), pp.~107--119.

\bibitem{cuturi13}
{\sc M.~Cuturi}, {\em Sinkhorn distances: Lightspeed computation of optimal
  transport}, in Advances in Neural Information Processing Systems 26, C.~J.~C.
  Burges, L.~Bottou, M.~Welling, Z.~Ghahramani, and K.~Q. Weinberger, eds.,
  Curran Associates, Inc., 2013, pp.~2292--2300.

\bibitem{daubechies}
{\sc I.~Daubechies}, {\em Ten Lectures on Wavelets}, Society for Industrial and
  Applied Mathematics, 1992.

\bibitem{dewet02}
{\sc T.~de~Wet}, {\em Goodnes-of-fit tests for location and scale families
  based on a weighted {$L_2$}-{W}asserstein distance measure}, Test, 11 (2002),
  pp.~89--107.

\bibitem{barrio99}
{\sc E.~del Barrio, J.~A. Cuesta-Albertos, C.~Matran, and J.~M.
  Rodriguez-Rodriguez}, {\em Tests of goodness of fit based on the
  {$L_2$}-{W}asserstein distance}, The Annals of Statistics, 27 (1999),
  pp.~1230--1239.

\bibitem{delon04}
{\sc J.~Delon}, {\em Midway image equalization}, Journal of Mathematical
  Imaging and Vision, 21 (2004), pp.~119--134.

\bibitem{dobrushin70}
{\sc R.~Dobrushin}, {\em Prescribing a system of random variables by
  conditional distributions}, Theory of Probability \& Its Applications, 15
  (1970), pp.~458--486.

\bibitem{engquist18}
{\sc B.~Engquist and H.~Zhao}, {\em Approximate separability of the {G}reen's
  function of the {H}elmholtz equation in the high frequency limit},
  Communications on Pure and Applied Mathematics, 0.

\bibitem{frogner15}
{\sc C.~Frogner, C.~Zhang, H.~Mobahi, M.~Araya, and T.~A. Poggio}, {\em
  Learning with a {W}asserstein loss}, in Advances in Neural Information
  Processing Systems 28, C.~Cortes, N.~D. Lawrence, D.~D. Lee, M.~Sugiyama, and
  R.~Garnett, eds., Curran Associates, Inc., 2015, pp.~2053--2061.

\bibitem{gangbo96}
{\sc W.~Gangbo and R.~J. McCann}, {\em The geometry of optimal transportation},
  Acta Mathematica, 177 (1996), pp.~113--161.

\bibitem{greenbaum}
{\sc A.~Greenbaum}, {\em Iterative Methods for Solving Linear Systems}, Society
  for Industrial and Applied Mathematics, 1997.

\bibitem{Helgason2011}
{\sc S.~Helgason}, {\em Integral Geometry and {R}adon Transforms}, Springer New
  York, New York, NY, 2011.

\bibitem{crb-book}
{\sc J.~S. Hesthaven, G.~Rozza, and B.~Stamm}, Springer Cham, Cham,
  Switzerland, 2016.

\bibitem{iollo14}
{\sc A.~Iollo and D.~Lombardi}, {\em Advection modes by optimal mass transfer},
  Phys. Rev. E, 89 (2014), p.~022923.

\bibitem{Kantorovich}
{\sc L.~V. Kantorovich}, {\em On a problem of {M}onge}, Uspekhi Mat. Nauk, 3
  (1948), pp.~225--226.

\bibitem{laxphilips}
{\sc P.~D. Lax and R.~S. Phillips}, {\em Scattering theory}, Bull. Amer. Math.
  Soc., 70 (1964), pp.~130--142.

\bibitem{fvmbook}
{\sc R.~J. LeVeque}, {\em Finite Volume Methods for Hyperbolic Problems},
  Cambridge University Press, Cambridge, 1st~ed., 2002.

\bibitem{mallows72}
{\sc C.~L. Mallows}, {\em A note on asymptotic joint normality}, The Annals of
  Mathematical Statistics, 43 (1972), pp.~508--515.

\bibitem{mccann97}
{\sc R.~J. McCann}, {\em A convexity principle for interacting gases}, Advances
  in Mathematics, 128 (1997), pp.~153 -- 179.

\bibitem{munk02}
{\sc A.~Munk and C.~Czado}, {\em Nonparametric validation of similar
  distributions and assessment of goodness of fit}, Journal of the Royal
  Statistical Society: Series B (Statistical Methodology), 60, pp.~223--241.

\bibitem{metivier16}
{\sc L.~Métivier, R.~Brossier, Q.~Mérigot, E.~Oudet, and J.~Virieux}, {\em An
  optimal transport approach for seismic tomography: application to 3d full
  waveform inversion}, Inverse Problems, 32 (2016), p.~115008.

\bibitem{natterer}
{\sc F.~Natterer}, {\em The Mathematics of Computerized Tomography}, Society
  for Industrial and Applied Mathematics, 2001.

\bibitem{siamrev-survey}
{\sc S.~G. P.~Benner and K.~Willcox}, {\em A survey of projection-based model
  reduction methods for parametric dynamical systems}, SIAM Rev., 57 (2015),
  pp.~483--531.

\bibitem{pitie07}
{\sc F.~Pitié, A.~C. Kokaram, and R.~Dahyot}, {\em Automated colour grading
  using colour distribution transfer}, Computer Vision and Image Understanding,
  107 (2007), pp.~123 -- 137.

\bibitem{pressdrt}
{\sc W.~H. Press}, {\em Discrete {R}adon transform has an exact, fast inverse
  and generalizes to operations other than sums along lines}, Proceedings of
  the National Academy of Sciences, 103 (2006), pp.~19249--19254.

\bibitem{rabin12}
{\sc J.~Rabin, G.~Peyr{\'e}, J.~Delon, and M.~Bernot}, {\em {W}asserstein
  barycenter and its application to texture mixing}, in Scale Space and
  Variational Methods in Computer Vision, A.~M. Bruckstein, B.~M. ter
  Haar~Romeny, A.~M. Bronstein, and M.~M. Bronstein, eds., Berlin, Heidelberg,
  2012, Springer Berlin Heidelberg, pp.~435--446.

\bibitem{Schulze15}
{\sc J.~Reiss, P.~Schulze, J.~Sesterhenn, and V.~Mehrmann}, {\em The shifted
  proper orthogonal decomposition: A mode decomposition for multiple transport
  phenomena},  (2015),
  \href{http://arxiv.org/abs/1512.01985}{arXiv:1512.01985}.

\bibitem{radonsplit}
{\sc D.~Rim}, {\em Dimensional splitting of hyperbolic partial differential
  equations using the {R}adon transform},  (2017),
  \href{http://arxiv.org/abs/1705.03609}{arXiv:1705.03609}.

\bibitem{rim17thesis}
{\sc D.~Rim}, {\em Uncertainty quantification problems in tsunami modeling and
  reduced order models for hyperbolic partial differential equations}, Ph.D.
  Thesis, University of Washington,  (2017).

\bibitem{dinterp}
{\sc D.~Rim}, {\em Code archive}, 2018.
\newblock \href{https://zenodo.org/record/1405576}{doi:10.5281/zenodo.1405576}.

\bibitem{rim18}
{\sc D.~Rim and K.~Mandli}, {\em Model reduction of a parametrized scalar
  hyperbolic conservation law using displacement interpolation}, Preprint,
  (2018), \href{http://arxiv.org/abs/1805.05938}{arXiv:1805.05938}.

\bibitem{rim18reversal}
{\sc D.~Rim, S.~Moe, and R.~LeVeque}, {\em Transport reversal for model
  reduction of hyperbolic partial differential equations}, SIAM/ASA Journal on
  Uncertainty Quantification, 6 (2018), pp.~118--150.

\bibitem{marsden1}
{\sc C.~W. Rowley and J.~E. Marsden}, {\em Reconstruction equations and the
  {Karhunen-Lo\`eve} expansion for systems with symmetry}, Physica D,  (2000),
  pp.~1--19.

\bibitem{solomon15}
{\sc J.~Solomon, F.~de~Goes, G.~Peyr{\'e}, M.~Cuturi, A.~Butscher, A.~Nguyen,
  T.~Du, and L.~Guibas}, {\em Convolutional {W}asserstein distances: Efficient
  optimal transportation on geometric domains}, ACM Trans. Graph., 34 (2015),
  pp.~66:1--66:11.

\bibitem{szego}
{\sc G.~Szeg\"o}, {\em Orthogonal Polynomials}, American Mathematical Society,
  New York, NY, USA, 1939.

\bibitem{trefethen}
{\sc L.~N. Trefethen}, {\em Approximation Theory and Approximation Practice},
  Society for Industrial and Applied Mathematics, Philadelphia, PA, USA, 2012.

\bibitem{villani03}
{\sc C.~Villani}, {\em Topics in Optimal Transportation}, American Mathematical
  Society, Providence, RI, 2003.

\bibitem{villani2008}
{\sc C.~Villani}, {\em Optimal transport: old and new}, vol.~338, Springer
  Science \& Business Media, 2008.

\bibitem{Welper17p}
{\sc G.~Welper}, {\em $h$ and $hp$-adaptive {I}nterpolation by {T}ransformed
  {S}napshots for {P}arametric and {S}tochastic {H}yperbolic {PDE}s}, {\tt
  arXiv:1710.11481 [math.NA]} (2017),
  \href{http://arxiv.org/abs/1710.11481}{arXiv:1710.11481}.

\bibitem{welper17}
{\sc G.~Welper}, {\em Interpolation of functions with parameter dependent jumps
  by transformed snapshots}, SIAM Journal on Scientific Computing, 39 (2017),
  pp.~A1225--A1250.

\bibitem{xiu10}
{\sc D.~Xiu}, {\em Numerical Methods for Stochastic Computations: A Spectral
  Method Approach}, Princeton University Press, Princeton, NJ, USA, 2010.

\bibitem{yang18}
{\sc Y.~Yang, B.~Engquist, J.~Sun, and B.~F. Hamfeldt}, {\em Application of
  optimal transport and the quadratic {W}asserstein metric to full-waveform
  inversion}, Geophysics, 83 (2018), pp.~R43--R62.

\bibitem{zygmund_fefferman}
{\sc A.~Zygmund and R.~Fefferman}, {\em Trigonometric Series}, Cambridge
  Mathematical Library, Cambridge University Press, 3~ed., 2003.

\end{thebibliography}
